\documentclass[graybox]{svmult}

\usepackage{type1cm}        
\usepackage{graphicx}        
\usepackage{multicol}        
\usepackage{multirow} 
\newcommand{\minitab}[2][l]{\begin{tabular}{#1}#2\end{tabular}}
\usepackage[bottom]{footmisc}
\usepackage{mathtools}

\usepackage{newtxtext}       %
\usepackage[varvw]{newtxmath}       

\newcommand{\topint}{\operatorname{int}}
 
\newcommand{\bc}{\mathcal{B}}
\newcommand{\pow}{\operatorname{pow}}
\newcommand{\size}{\operatorname{size}}

\newcommand{\cL}{\mathcal{L}}

\newcommand{\one}{\operatorname{1\hspace*{-0.55ex}I}}

\newcommand{\N}{\ensuremath{\mathbb{N}}} 
\newcommand{\R}{\ensuremath{\mathbb{R}}} 

\renewcommand{\phi}{\varphi}

\newcommand{\B}{\ensuremath{\mathcal{B}}} 

\renewcommand{\E}{\ensuremath{\mathbf{E}}} 
\renewcommand{\P}{\ensuremath{\mathbf{P}}} 
\newcommand{\Fc}{\ensuremath{\mathcal{F}}} 
\newcommand{\Ac}{\ensuremath{\mathcal{A}}} 
\newcommand{\Cc}{\ensuremath{\mathcal{C}}} 
\newcommand{\ind}{\ensuremath{\mathbf{1}}} 
\newcommand{\Ns}{\ensuremath{\mathsf{N}}} 
\newcommand{\Nf}{\ensuremath{\mathfrak{N}}}

\newcommand{\Tb}{\ensuremath{\mathbb{T}}}

\newcommand{\Lc}{\ensuremath{\mathcal{L}}} 
\newcommand{\Rc}{\ensuremath{\mathcal{R}}} 

\newcommand{\cH}{\ensuremath{\mathcal{H}}}

\newcommand{\Sc}{\ensuremath{\mathcal{S}}}

\newcommand{\Tc}{\ensuremath{\mathcal{T}}}

\newcommand{\conv}{{\rm conv}\,}

\newcommand{\kommentar}[1]{}

\makeindex             

\begin{document}

\title*{Random Tessellations - An Overview of Models}
\author{Claudia Redenbach, Christian Jung}
\institute{Claudia Redenbach \at RPTU University Kaiserslautern-Landau, Kaiserslautern, \email{claudia.redenbach@rptu.de}
\and Christian Jung \at RPTU University Kaiserslautern-Landau, Kaiserslautern, \email{christian.jung@rptu.de}}
%
%
\maketitle

\abstract{Random tessellations are a prominent class of models in stochastic geometry. In this chapter, we give an overview of mechanisms that have been used to formulate random tessellation models. First, the notion of a random tessellation and basic geometric
characteristics of random tessellations are introduced. Then, several model classes are presented. This includes, but is not limited to, Voronoi tessellations and their weighted generalizations, hyperplane tessellations, and STIT tessellations. Simulation of the tessellation models and approaches for model fitting are also discussed.}

\section{Introduction}
\label{sec:1}
A tessellation is a division of space
into cells that only intersect in their boundaries.
In practice, such structures are observed in biological cells, honeycombs, the cells of a foam or a polycrystalline material as well as crack patterns in soil. On a very different scale, also the road network in a city or influence zones of supermarkets or restaurants can be modelled as tessellations. 

Taking the geometry of these structures and their generation mechanisms into account, various approaches for tessellation construction can be formulated. In stochastic geometry, models that randomize these constructions have been proposed. 
For fitting such models to observed structures, characteristics derived from the distributions of cell size or shape as well as topological characteristics of the tessellation can be used.  
Some models are analytically well tractable, such that explicit results for such characteristics are available. In contrast, other tessellation models can only be studied by Monte-Carlo simulation. 

The aim of this work is to give an overview of construction mechanisms for random tessellations. Analytical results are sketched briefly. Additionally, aspects of simulation and fitting of random tessellation models are discussed.

Complementary chapters on random tessellations can be found in the stochastic geometry books by Chiu, Stoyan, Kendall and Mecke \cite{skm13} and Schneider and Weil \cite{SchWei08}. Monographs focusing on Voronoi tessellations are \cite{AurenhammerBook} and \cite{OkaBooSugChi00}. A recent monograph on Poisson hyperplane tessellations is \cite{HugSchneider2024}.
A review of asymptotic results for random tessellations, mostly Poisson-Voronoi and Poisson hyperplane tessellations, is given in the book chapters by Calka \cite{Calka2009_NewPerspectives,Calka2013LNM}.

\section{Stochastic geometry background}
\label{sec:RandomTess}
We start by summarizing some concepts from stochastic geometry that will be needed to introduce the tessellation models. For a detailed introduction into the field we refer to \cite{skm13,SchWei08}.

\subsection{Point and hyperplane processes}\label{sec_defs_1}

We start by defining the notion of a point process. Roughly speaking, point processes are models for random collections of points in a given space. In most cases, we consider the Euclidean space $\R^d$. However, we also want to consider the case that "points" are lines or compact sets in $\R^d$. Hence, the following definition is formulated for a general space $E$.

\begin{definition}[Point process]
Let $E$ be a locally compact space with countable basis and Borel $\sigma$-algebra $\bc$. Denote by $\Ns(E)$ the set of all locally finite measures of the form
$$
\phi= \sum_{i \in \N} \delta_{x_i}
$$
where $x_i\in \R^d$ and $\delta_{x_i}$ are Dirac measures. 
Equip $\Ns(E)$ with the $\sigma$-algebra $\Nf$ generated by the mappings $\phi \mapsto \phi(B)$ for $B\in \bc$.
Then a {\it point process} is a random variable $\Phi$ on a probability space $(\Omega, \Ac, \P)$ taking values in the measurable space $(\Ns(E), \Nf)$.
If $E= \R^d \times M$, where $M$ is a locally compact space with countable basis, we call $\Phi$ a marked point process with mark space $M$. In this case, a point $(x,m) \in \Phi$ is interpreted as a point $x\in \R^d$ to which a mark $m \in M$ is attached. 
\end{definition}

The support $\{x_1, x_2, \ldots\}$ of $\phi$ forms a locally finite collection of points of $E$ which we will also call $\phi$. Hence, realizations of point processes can both be interpreted as measures and as point patterns.

Examples of marks could be arrival times of the points, some random weights associated to the points, or statistical information on the shape of objects located at the points of the point process.   
\begin{definition}[Stationarity and isotropy]
\label{def:Stat}
A point process on $\R^d$ is called {\it stationary} if its distribution is invariant under translations and {\it isotropic} if it is invariant under rotations.
\end{definition}

\begin{definition}[Intensity measure]
For any point process $\Phi$ on $E$, its {\it intensity measure} $\Lambda:\bc\to[0,\infty]$ reports the expected number of points in Borel sets $B$, i.e, 
$$
\Lambda(B)= \E\Phi(B), \quad B \in \bc.
$$
For stationary point processes on $\R^d$ we have $\Lambda= \lambda \nu_d$, where $\nu_d$ denotes the $d$-dimensional Lebesgue measure and $\lambda = \E\Phi([0,1]^d)$ is the {\it intensity} of $\Phi$, that is, the expected number of points in the unit cube. 
\end{definition}

\begin{definition}[Poisson process]
A {\it Poisson (point) process} $\Phi$ on $E$ with intensity measure $\Lambda$ is characterised by the following properties:
\begin{enumerate}
\item
The number of points $\Phi(B)$ contained in a Borel set $B \in \bc$ with $\Lambda(B) <\infty$ has a Poisson distribution with parameter $\Lambda(B)$. 
\item
For arbitrary $k\in \N$, the numbers of points in $k$ disjoint Borel sets are independent random variables.
\end{enumerate}
\end{definition}

A stationary Poisson point process on $\R^d$ is obtained when using $\Lambda(B) = \lambda \nu_d(B)$ for some $\lambda>0$. Such a process is also isotropic.

Let $\Lc$ denote the space of $d-1$-dimensional affine subspaces of $\R^d$. The elements of $\Lc$ are hyperplanes in $\R^d$.  Furthermore, let $\cL_0$ denote the space of hyperplanes intersecting the origin. The line orthogonal to $L \in \cL_0$ is denoted by $L^{\bot}$.
Equip $\Lc$ with the hit-and-miss $\sigma$-algebra generated by
$$ \left\{\{ L_i\}_{i} \, :\, \big[ \bigcup_{i} L_i \big]\cap K \neq \emptyset \right\}, K \subset \R^d \mbox{ compact}.
$$

\begin{definition}[Hyperplane process]
A \emph{hyperplane process} in $\R^d$ is a point process in $\Lc$. A \emph{Poisson hyperplane process} is a Poisson point process on the space $E=\Lc$.
\end{definition}

Stationarity and isotropy for hyperplane processes can be defined in analogy to Definition~\ref{def:Stat}.
Let $X$ be a stationary hyperplane process in $\R^d$ with intensity measure $\Lambda\neq 0$. Then there are a number $\lambda >0$ and a probability measure $\Theta$ on $\cL_0$ with
\begin{equation}
\label{thm:eqLambda}
\int_{\cL} f(E) \Lambda(dE) = \lambda \int_{\cL_0} \int_{L^{\bot}} f(L+x) \nu_{L^{\bot}}(dx) \Theta(dL)   
\end{equation}
for all nonnegative measurable functions $f$ on $\cL$. $\lambda$ and $\Theta$ are uniquely determined by $\Lambda$. 
The number $\lambda$ is the \emph{intensity} of the hyperplane process. It can be interpreted as the mean total $(d-1)$-content of the hyperplanes per unit volume, see \cite[Theorem 4.4.3]{SchWei08}. 
The distribution $\Theta$ is the \emph{directional distribution} of $X$. By considering the normal direction of a hyperplane $H$, $\Theta$ induces a distribution $\Rc$ on $S^{d-1}$. $\Rc$ is called the \emph{rose of normal directions}. 

If $X$ is isotropic, then $\Theta$ is rotation invariant. In this case, $\Theta$ is the Haar measure $\nu$ on $\cL_0$ and the normal directions of the hyperplanes follow a uniform distribution on the unit sphere $S^{d-1}$.

\subsection{Tessellations and their properties}

\begin{definition}[Tessellation]
\label{def:tessellation}
A {\it tessellation} of $\R^d$ is a locally finite collection $T = \{C_i\, :\, i\in \N\}$ of compact sets $C_i$ with interior points such that $\topint(C_i) \cap \topint(C_j) = \emptyset$ for $i \neq j$ and
$\bigcup\limits_{i \in \N} C_i =\R^d$.
Locally finite means that $\#\{C \in T \,:\, C \cap B \neq \emptyset \}< \infty$ for all bounded $B \subset \R^d$. The sets $C_i \in T$ are the {\it cells} of the tessellation $T$.
\end{definition}

For the rest of the section, we will assume that $T$ is a tessellation with convex cells. In this case, the cells are $d$-dimensional polytopes \cite[Lemma 10.1.1]{SchWei08}.

\begin{definition}[k-faces]
The {\it faces} of a convex polytope $P$ are the intersections of $P$ with its supporting hyperplanes \cite[Section 2.4]{Schn93}. Let $P$ be a $d$-dimensional polytope and $k\in \{0,\ldots,d-1\}$. A $k$-dimensional face of $P$ is called a {\it $k$-face}. Then the $0$-faces of $P$ are the {\it vertices}, the $1$-faces the {\it edges}, and the $(d-1)$-faces the {\it facets}. 
For convenience, the polytope $P$ is considered as a $d$-face.  
Write $\Sc_k(P)$ for the set of $k$-faces of a polytope $P$ and 
$\Sc_k(T)=\bigcup_{C\in T} \Sc_k(C)$ for the set of $k$-faces of all cells $C$ of $T$. Furthermore, let \begin{equation} \label{eq:Faces}
F(y)=\bigcap_{C \in T: y\in C} C,\quad y\in \R^d,\end{equation}
be the intersection of all cells of the tessellation containing the point $y$. Then $F(y)$ is a finite intersection of $d$-polytopes and, since it is nonempty, $F(y)$ is a $k$-dimensional polytope for some $k \in \{0,\ldots,d\}$. Therefore, we may introduce
\[ \Fc_k(T)=\{F(y)\,:\, \dim F(y)=k, y \in \R^d \},\quad k=0, \ldots, d,\]
the set of{ \it $k$-faces of the tessellation} $T$. A $k$-face $H \in \Sc_k(T)$ of a cell $C$ of $T$ is the union of all those $k$-faces $F \in \Fc_k(T)$ of the tessellation contained in $H$.
\end{definition}

\begin{figure}[t]
\begin{center}
\includegraphics[width=3.7cm]{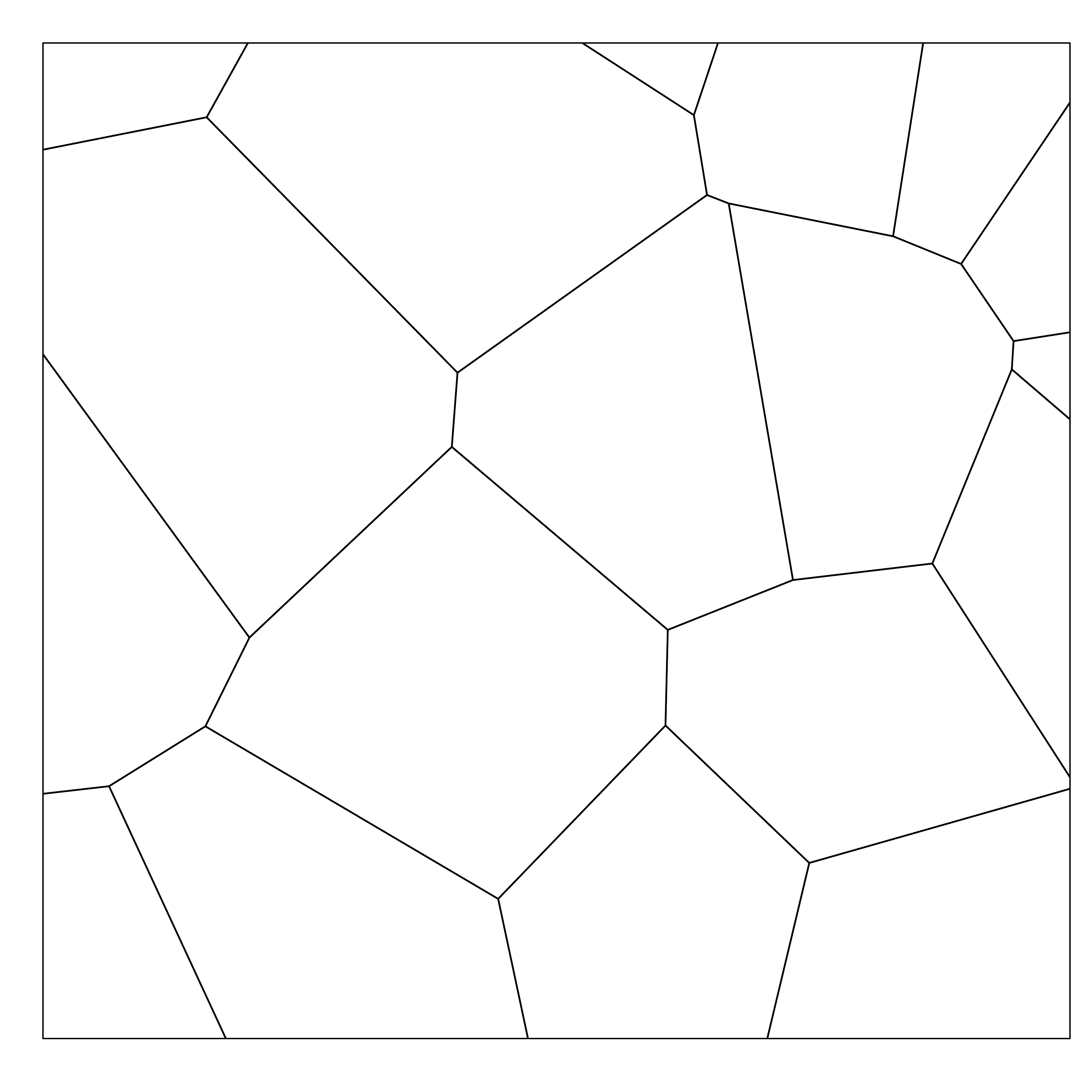}  \includegraphics[width=3.7cm]{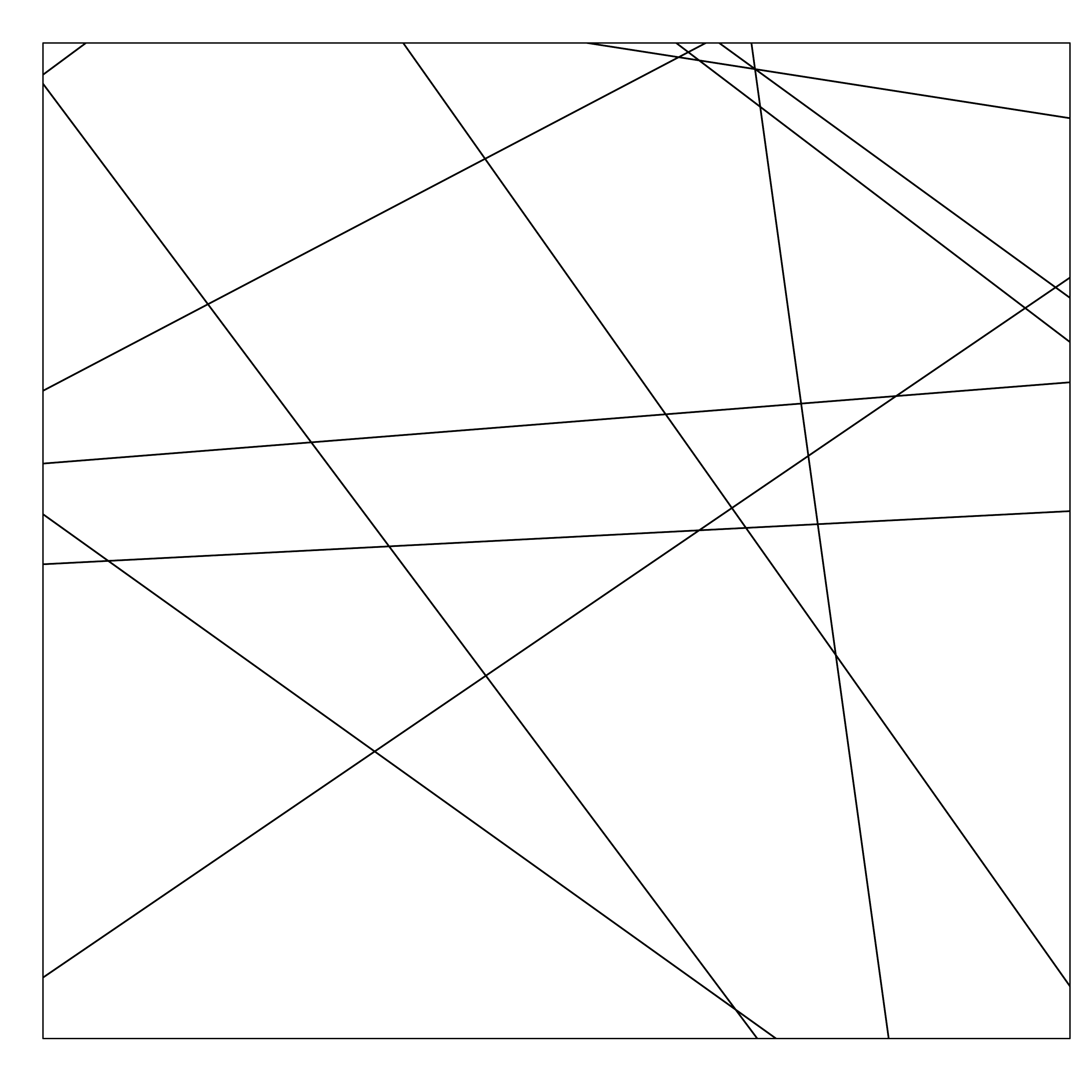}  \includegraphics[width=3.7cm]{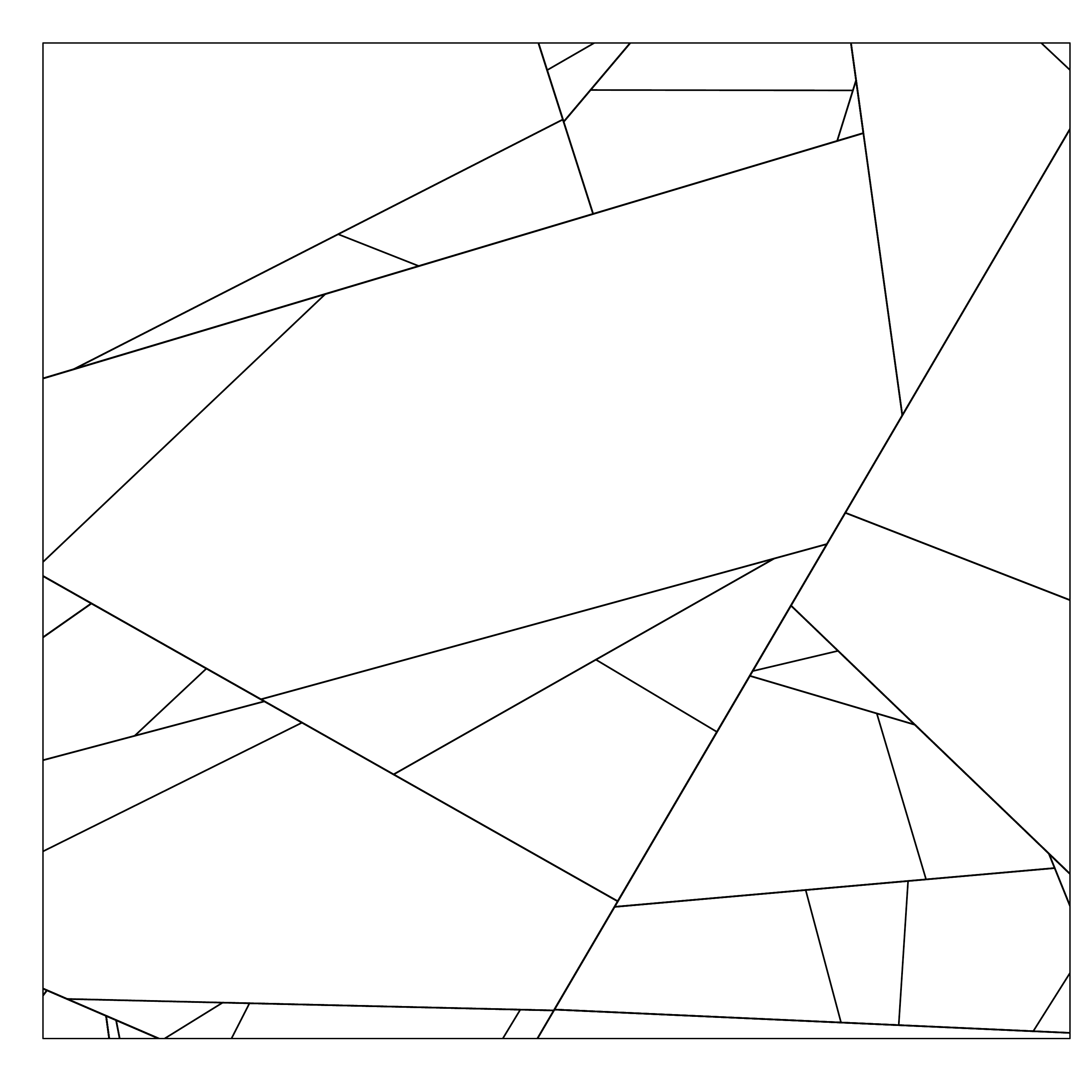}\\
\caption{\label{fig:Properties} Left: A normal tessellation in $\R^2$. Every vertex is contained in exactly three cells. Vertices are Y-shaped. Middle: A face-to-face, but not normal tessellation.  Vertices belong to four cells and are X-shaped. Right: A tessellation that is not face-to-face. Intersections of cells are not necessarily edges of cells. Vertices are T-shaped.}
\end{center}
\end{figure}

\begin{definition}[Face-to-face and normal tessellation]
A tessellation $T$ is called {\it face-to-face} if the faces of the cells and the faces of the tessellation coincide, i.e. if $\Sc_k({T})=\Fc_k({T})$ for all $k=0, \ldots, d$. For $k=0$ and $k=d$ this is always true. In the case of face-to-face tessellations we will unify the notation writing $\Fc_k(C)$ for the set of $k$-faces of a cell $C$ of $T$. 
A tessellation $T$ is called {\it normal} if it is face-to-face and every $k$-face of $T$ is contained in the boundary of exactly $d-k+1$ cells for $k=0, \ldots, d-1$. See Fig.~\ref{fig:Properties} for illustrations of these concepts.
\end{definition}

Face-to-face tessellations are sometimes also called \emph{side-to-side} or \emph{regular}. Here, we use the term 'regular' as a property of cell shapes.

Probably the most striking difference of the tessellation models shown in Fig.~\ref{fig:Properties} is the different structure of the vertices. 

\begin{definition}[X-, Y-, and T-vertices]
A vertex of a planar tessellation is called
\begin{enumerate}
\item 
\emph{X-shaped} or an \emph{X-vertex}, if it lies in the intersection of exactly four edges. 
\item 
\emph{Y-shaped} or a \emph{Y-vertex}, if it lies in the intersection of exactly three edges, none of which are aligned.
\item 
\emph{T-shaped} or a \emph{T-vertex}, if it lies in the intersection of exactly three edges, two of which are aligned.
\end{enumerate}
   
\end{definition}

\subsection{Random tessellations}
\label{Sec:RandomTess}
We write $\Tb$ for the set of all tessellations in $\R^d$. It is equipped with the $\sigma$-algebra $\Tc$ generated by $$ \left\{\{ C_i\}_{i} \, :\, \big[ \bigcup_{i} \partial C_i \big]\cap K \neq \emptyset \right\}, K \subset \R^d \mbox{ compact}.
\label{Eq:sigmaAlgebra}
$$

\begin{definition}[Random tessellation]
A {\it random tessellation} in $\R^d$ is a random variable $X$ on a probability space $(\Omega, \Ac, \P)$ with range $(\Tb, \Tc)$. It is called {\it normal }and {\it face-to-face} if its realisations are almost surely normal and face-to-face, respectively. 

The translation and the rotation of a tessellation $T\in \Tb$ are defined via
\begin{align*}
T+y&=\{C+y\,:\, C\in T \}, \quad y\in \R^d, \mbox{ and}\\
\vartheta T&=\{\vartheta C\, :\, C\in T \}, \quad \vartheta \in SO_d.
\end{align*}
A random tessellation is called {\it stationary} if its distribution is invariant under translations and {\it isotropic} if it is invariant under rotations.   
\end{definition}

Stationary random tessellations contain infinitely many cells. For the definition of geometric tessellation characteristics, a finite subsample of the aggregate of cells is investigated. For this purpose, only cells with centroid in a given reference set are considered.

\begin{definition}[Centroid of a compact set]
Write $\Cc'$ for the system of compact nonempty sets in $\R^d$. Let $c:\Cc' \to \R^d$ be a measurable function such that 
\begin{equation}
\label{centroid}
c(C+y)=y+c(C), \quad y \in \R^d, \quad C\in\Cc'.\end{equation} The point $c(C)$ is called the {\it centroid} of the set $C\in\Cc'$. 
\end{definition}

Typical choices of centroids are the centre of gravity of the set $C$, the centre of its surrounding ball, or the ``extreme" point of $C$ with respect to a given direction. 

\begin{figure}[b]
\begin{center}
\includegraphics[width=3.7cm]{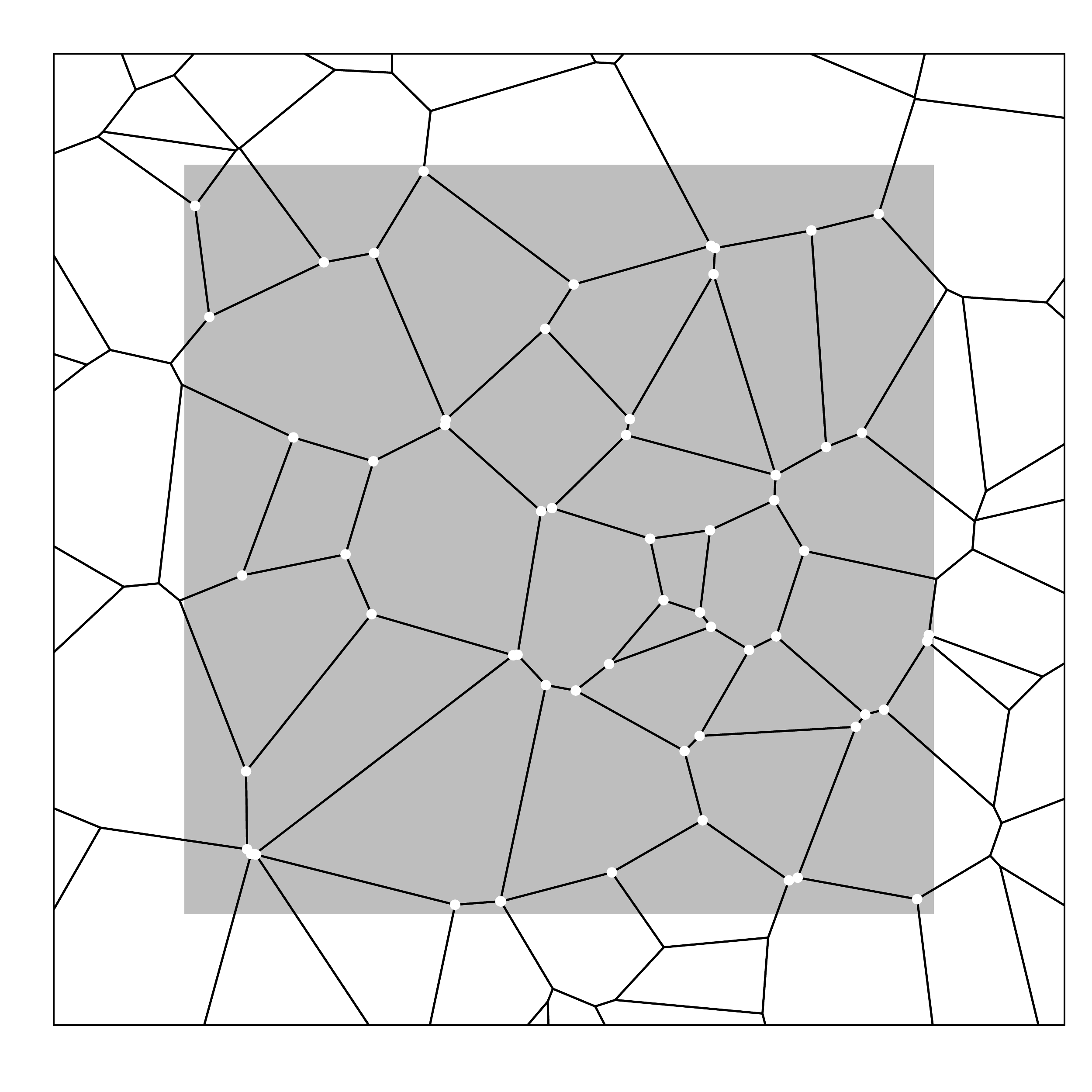} \includegraphics[width=3.7cm]{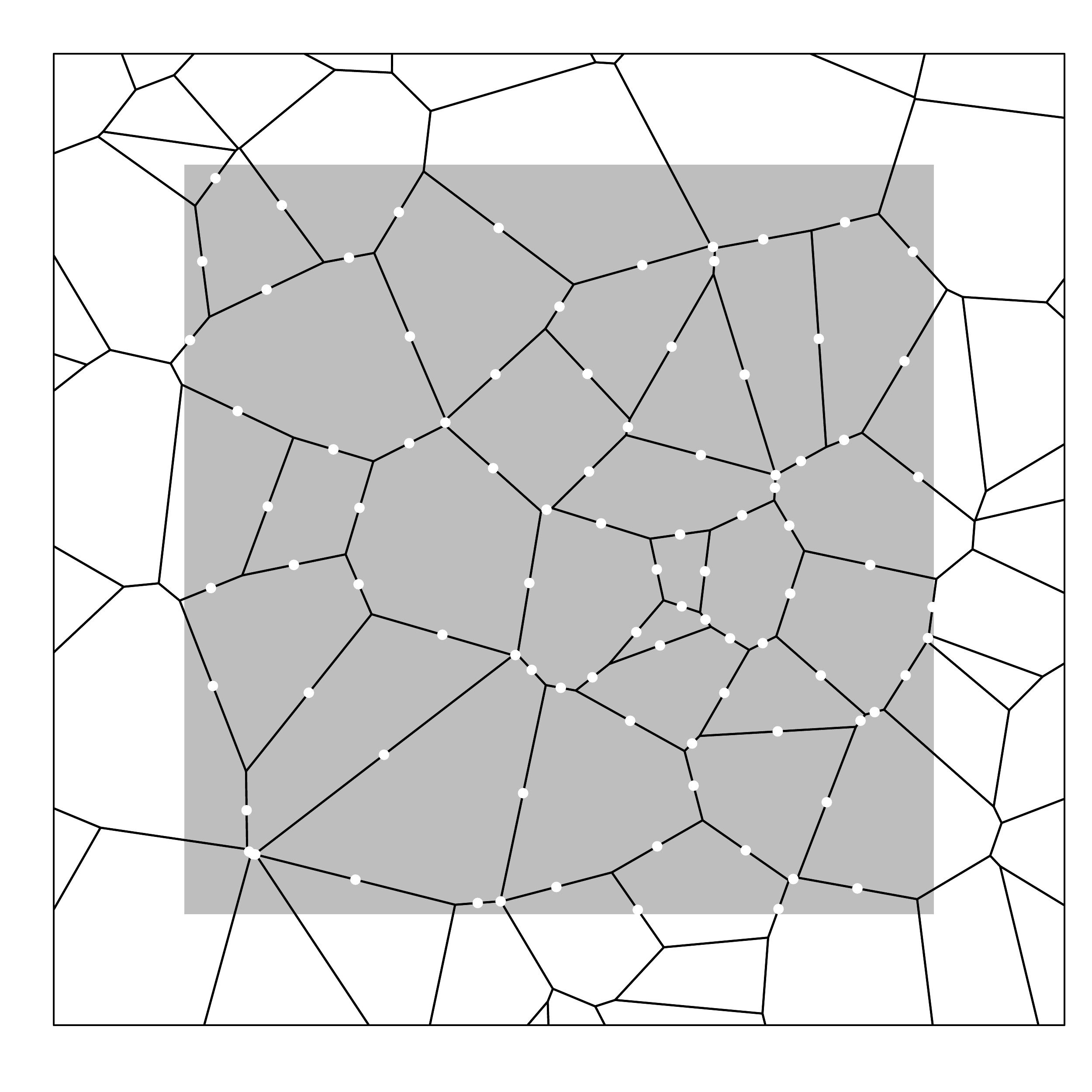} \includegraphics[width=3.7cm]{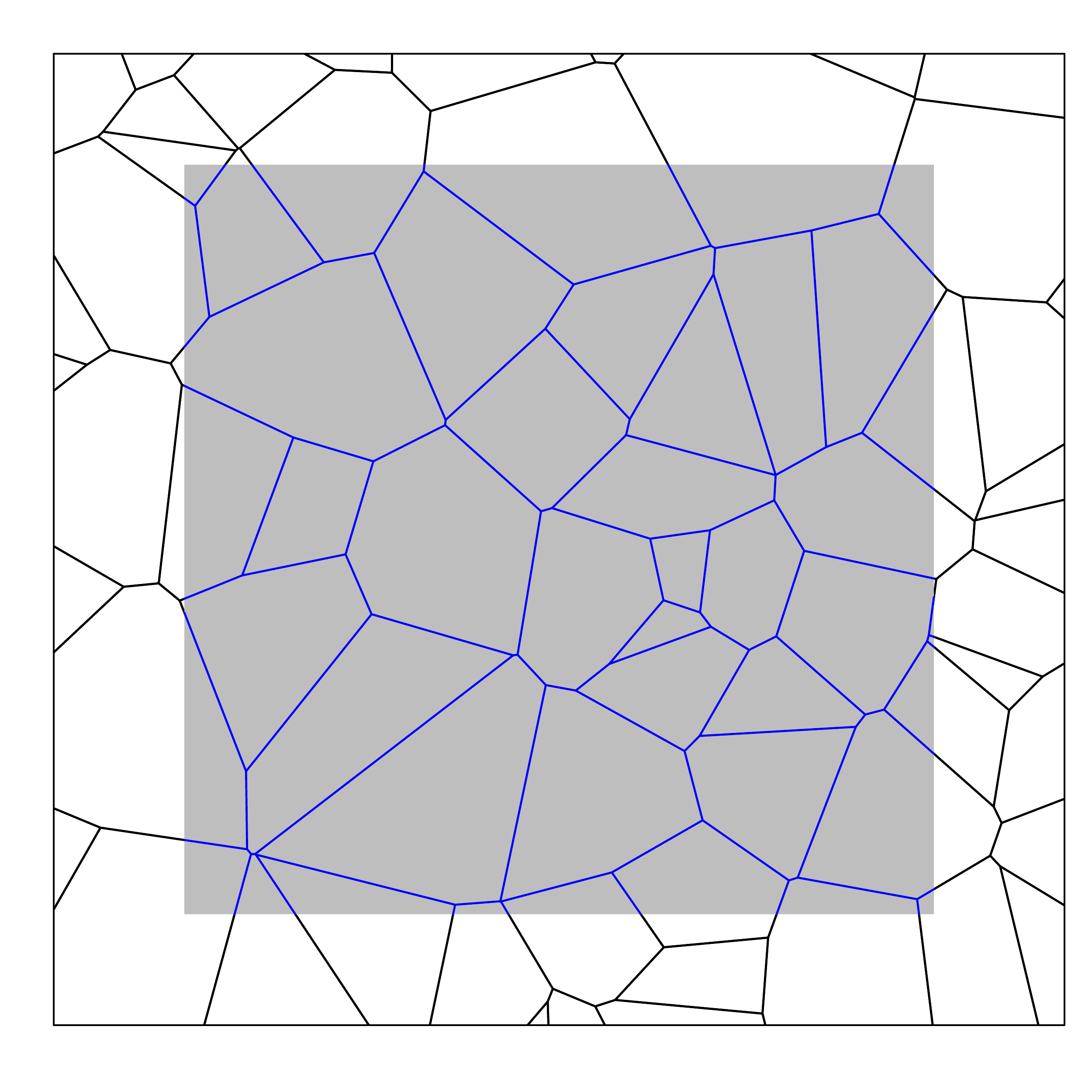}\\
\caption{ \label{fig:Intensities}Left: Point process $\Phi_0$ of vertices. Middle: Point process $\Phi_1$ of edge centers. The intensities $\gamma_0$ and $\gamma_1$ are the expected numbers of points of $\Phi_0$ and $\Phi_1$, respectively, in the unit square (grey). Right: The measure $M_1$ measures the total edge length in a given set $B$. The intensity $\mu_1$ is the expected total edge length (blue) in the unit square (grey).}
\end{center}
\end{figure}

Let $c_k$ denote a centroid function acting on the set of $k$-faces of a random tessellation $X$. Then we can define the point process $\Phi_k$ of centres of the $k$-faces of $X$ as 
\begin{equation*}
\Phi_k(X)=\sum\limits_{F\in \Fc_k(X)} \delta_{c_k(F)}.
\end{equation*} 
See Fig.~\ref{fig:Intensities} for an illustration of $\Phi_0$ and $\Phi_1$.

Stationarity of the random tessellation $X$ and property \eqref{centroid} imply stationarity of the point processes $\Phi_k(X)$.
The intensity $\gamma_k$ of $\Phi_k$ is given by the formula
\begin{equation*}
\gamma_k=\E \left[\sum_{F \in \Fc_k(X)} \ind_{[0,1]^d}(c_k(F)) \right], \quad k=0,\ldots,d,
\end{equation*} 
and can be interpreted as the mean number of $k$-faces per unit volume. The value of $\gamma_k$ does not depend on the choice of the centroid function $c_k$ \cite{Mol89,SchWei08}.

Further random measures induced by a random tessellation are the measures $M_k:\bc\to[0,\infty]$ given by
\begin{equation*}
M_k(B):=  \sum_{F \in \Fc_k(X)}\cH^k(F\cap B), \quad k=0,\ldots,d, B\in \bc,
\end{equation*}
where $\cH^k$ is the $k$-dimensional Hausdorff measure. E.g. $M_1(B)$ measures the mean total edge length in the set $B$.
In the stationary case, their intensities 
\begin{equation*}
\mu_k=  \E \left[\sum_{C \in \Fc_k(X)}\cH^k(C\cap [0,1]^d )\right], \quad k=0,\ldots,d,
\end{equation*}
can be interpreted as the mean total $k$-content of the $k$-faces of the tessellation per unit volume. For $k=0$, we have $\mu_0=\gamma_0$. As a tessellation is space filling, we additionally get $\mu_d=1$.

In the following, we will also use the more established notation from stereology. For a tessellation in $\R^2$, $\mu_1= L_A$ denotes the mean total edge length per unit area. In the 3D case, we have $\mu_1=L_V$, the mean total edge length per unit volume, and $\mu_2=S_V$, the mean total cell surface area per unit volume.

The above characteristics carry some information on the aggregate of the tessellation's cells. In applications it is often interesting to investigate characteristics of the single cells and their faces leading e.g. to the distributions of cell volumes or edge lengths. This is formalised using the typical $k$-face of the stationary tessellation $X$ which is defined by means of Palm theory \cite{SchWei08}. 

\begin{definition}[Typical k-face] 
The \textit{typical $k$-face} $\Cc_k$ of a random stationary tessellation is a $\Cc'$-valued random variable such that
\begin{equation}
\label{eq:typical}
\E f(\Cc_k)=\frac{1}{\gamma_k \nu_d(B)} \E \left( \sum_{ \{C \in \Fc_k \, :\, c_k(C) \in B \} } f(C-c_k(C)) \right).
\end{equation}
\end{definition}
for every bounded, measurable and translation invariant function $f: \Cc' \to \R$ and every Borel set $B \in \bc$ with $0< \nu_d(B) <\infty$.

The typical $k$-face can be interpreted as a $k$-dimensional polytope picked at random from the system of $k$-faces of the tessellation, where each cell is drawn with the same probability. To avoid sampling from an infinite number of $k$-faces, sampling is restricted to a reference set $B$. By choosing $f= \one_A$ with $A \in \B(\Fc')$, equation \eqref{eq:typical} allows to compute the probability that a $k$-face of the tessellation has the property represented by $A$.

\begin{figure}[t]
\begin{center}
\includegraphics[width=3.7cm]{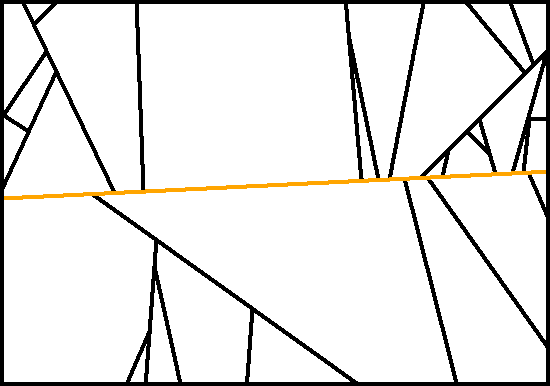} \hspace{0.1cm} \includegraphics[width=3.7cm]{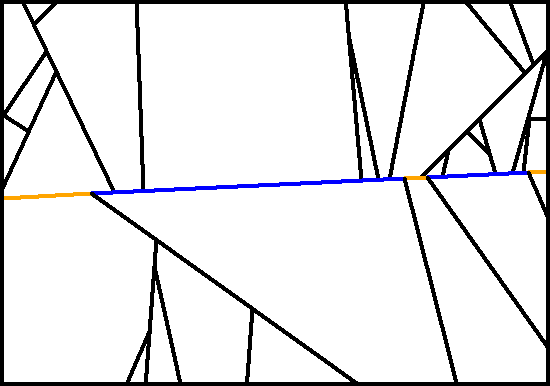}\hspace{0.1cm} 
\includegraphics[width=3.7cm]{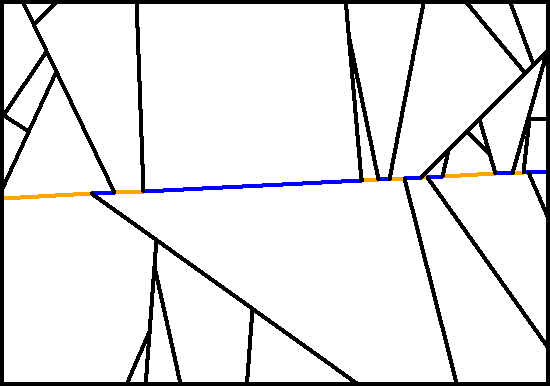}
\caption{\label{fig:STITsegments} Left: An I-segment in a tessellation with T-vertices. Middle: J-segments defined by the cells below the I-segment. Right: K-segments.}
\end{center}
\end{figure}

For non-face-to-face tessellations, notions beyond the $k$-faces of the tessellation cells have been introduced. For edges of planar tessellations with T-vertices, the notions K-, J-, and I- segment are used \cite{mackisack_miles_1996}. A \emph{K-segment} is bounded by two vertices of the tessellation with no vertex in its interior. That is, K-segments are the elements of the set $\Fc_1(X)$. \emph{J-segments} correspond to the elements of the set $\Sc_1(X)$, the edges (1-faces) of the cells. \emph{I-segments} are the maximal unions of collinear line segments. See Fig.~\ref{fig:STITsegments} for an illustration. A T-vertex that is contained in the relative interior of a J-segment is called a \emph{$\pi$-vertex}. A tessellation is face-to-face if it does not possess any $\pi$-vertices. For a generalization to the three-dimensional case, we refer to \cite{TW}.

For a tessellation in $\R^2$, we can define the following further characteristics:
\begin{align*}
N_{kl}&\mbox{ - the expected number of }
       l\mbox{-faces adjacent to the typical }k\mbox{-face, } k,l \in \{0,1,2\}, \\
      & \mbox{e.~g. } N_{02}\mbox{ denotes the expected number of cells neighbouring the typical vertex},\\
L_1& \mbox{ - the expected length of the typical edge},\\
P_2& \mbox{ - the expected perimeter (total edge length) of the typical cell},\\
A_2& \mbox{ - the expected area of the typical cell}
\end{align*}

The above characteristics are closely related. 

\begin{theorem}[\cite{SchWei08}, Theorem 10.1.6]
Let $T$ be a stationary random tessellation in $\R^2$. If $T$ is face-to-face, the  following mean value relations hold:
\begin{equation}
\label{eq:TessChar1}
\gamma_1 = \gamma_0+ \gamma_2
\end{equation}
\begin{equation}
\label{eq:TessChar2}
N_{02}= 2+2 \frac{\gamma_2}{\gamma_0}, \quad N_{20}= 2 +2 \frac{\gamma_0}{\gamma_2}
\end{equation}
\begin{equation}
\label{eq:TessChar3}
A_2= \frac{1}{\gamma_2}, \quad  P_2 = 2 \frac{\gamma_1}{\gamma_2} L_1, \quad  \mu_1=L_1 \gamma_1
\end{equation}
\begin{equation}
\label{eq:TessChar4}
3 \le N_{02}, N_{20} \le 6
\end{equation}
If $T$ is normal, then $N_{02}=3$ and $N_{20}=6$. Then \eqref{eq:TessChar2} and  \eqref{eq:TessChar1} yield
$$
\gamma_0=2\gamma_2, \quad \gamma_1= 3\gamma_2.
$$
\end{theorem}

This theorem shows in particular, that for a planar normal random tessellation in $\mathbb{R}^2$ all mean values can be expressed by the two parameters $\mu_0(=\gamma_0)$ and $\mu_1(=L_A)$. For nonnormal but face-to-face tessellations,  the three parameters $\mu_0(=\gamma_0)$, $\gamma_1$, and $\mu_1(=L_A)$ are required. For non-face-to-face tessellations, a fourth parameter is needed: the proportion $\phi$ of $\pi$-vertices \cite{Weiss:2011:TRS}.

Mean value formulas for stationary rectangular tessellations that are not necessarily face-to-face are proven in \cite{mackisack_miles_1996}. They express all mean values in terms of three parameters: the intensity and the mean length of I-segments, and the mean number of X-vertices on the typical I-segment.

Similar results can be obtained for the following characteristics of random tessellations in $\R^3$.
\begin{align*}
N_{kl}&\mbox{ - the expected number of }
       l\mbox{-faces adjacent to the typical }k\mbox{-face, } k,l \in \{0,1,2,3\}, \\
      & \mbox{e.~g. } N_{13}\mbox{ denotes the expected number of cells neighbouring the typical edge},\\
L_1& \mbox{ - the expected length of the typical edge},\\
P_2& \mbox{ - the expected perimeter of the typical face},\\
A_2& \mbox{ - the expected area of the typical face},\\
B_3& \mbox{ - the expected mean width of the typical cell},\\
L_3& \mbox{ - the expected total edge length of the typical cell},\\
S_3& \mbox{ - the expected surface area of the typical cell},\\
V_3& \mbox{ - the expected volume of the typical cell}.
\end{align*}

 For normal tessellations in $\R^3$, mean value relations can be formulated in terms of $\gamma_3$, $\mu_0 (=\gamma_0)$, $\mu_1(=L_V)$, and $\mu_2(=S_V)$ \cite{Mec84}. For nonnormal, but face-to-face tessellations, three additional parameters are required. A list of the corresponding mean value relations can, e.g., be found in \cite{skm13}. For tessellations that are not face-to-face, the set of parameters has to be extended by another four parameters \cite{Weiss:2011:TRS}. The range of topological parameters that can be realized is given in \cite{cowan_constraints_2015}.

Besides the typical cell also the \emph{zero cell} $\Cc_0$ of a random tessellation can be considered. It is the almost surely unique cell containing the origin (if it exists).
	
\begin{theorem}{\cite[Theorem 10.4.1]{SchWei08}}
	The zero cell $\Cc_0$ and the typical cell $\Cc$ are related by
	$$
	\E f(\mathcal{C}_0) = \lambda_d \E(f(\mathcal{C})V (\mathcal{C}))
	$$
	for any translation invariant, nonnegative, measurable function $f$, where $V$ denotes volume ($d$-dimensional Lebesgue measure). 
\end{theorem}	
Hence, the distribution of the zero cell is the volume-weighted distribution of the typical cell. This result implies that the zero cell is stochastically larger than the typical cell. 
Intuitively, this can be explained as follows. As described above, the typical cell is obtained by sampling at random from all cells of the tessellation having their centroid in a reference set $B$. Each cell is chosen with the same probability irrespective of its size or shape. 
A sampling procedure for the zero cell is to sample a uniform random point from a window $B$ that marks the origin. In the planar case, this can be imagined as throwing a dart at $B$ and calling the hitting point the origin. Clearly, larger cells are more likely to be hit by the dart. Indeed, the hitting probability is proportional to the area of the cell which explains the finding above.


\section{Hyperplane tessellations}
\label{sec:HyperplaneTess}

A tessellation in the plane can be induced by cutting the plane along a set of straight lines. In higher dimensions, sections along hyperplanes induce tessellations in a similar manner.

\begin{definition}[Hyperplane tessellation]
\label{def:hyperplane}
Let $\cH$ be a locally finite set of hyperplanes in $\R^d$.  
Then the connected components of $\R^d\setminus \bigcup_{H\in \cH} H$ form a system of open subsets of $\R^d$. Their closures build the \textit{hyperplane tessellation} induced by $\cH$. 
\end{definition}

Definition \ref{def:hyperplane} only yields a tessellation in the sense of Definition \ref{def:tessellation} if the resulting cells are bounded. For instance, systems of hyperplanes that are all parallel to each other do not induce a tessellation of $\R^d$.

One option to ensure that the hyperplane tessellation is well defined is to choose $\cH$ as a realization of a Poisson hyperplane process such that $\Rc$ is not concentrated on a great sphere in $S^{d-1}$ (a single direction in case of $d=2$). We call such a process \emph{nondegenerate}.
Hyperplane tessellations generated by nondegenerate Poisson hyperplane processes are called \textit{Poisson hyperplane tessellations}. 

\begin{figure}[t]
\begin{center}
\includegraphics[width=3.9cm]{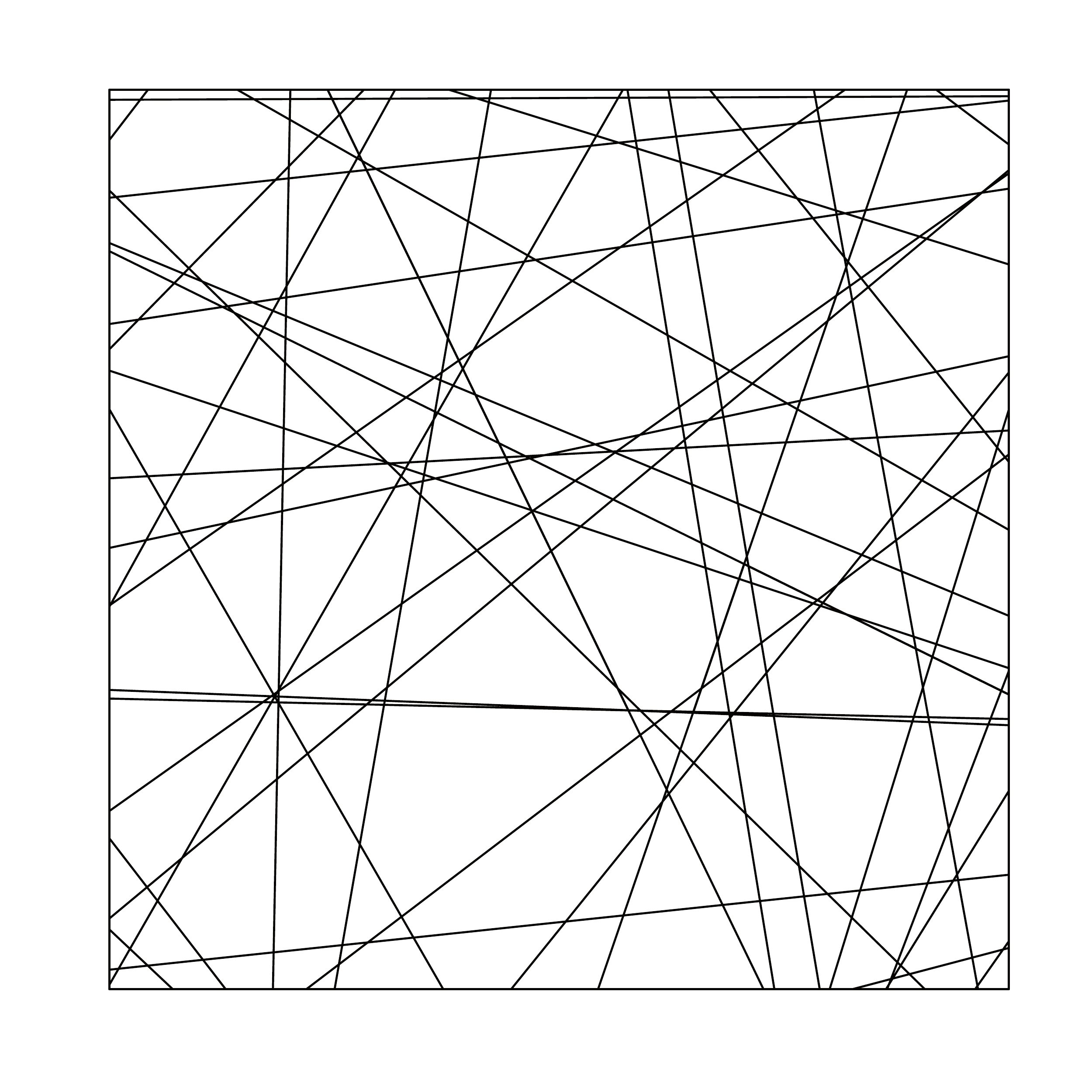}\includegraphics[width=3.9cm]{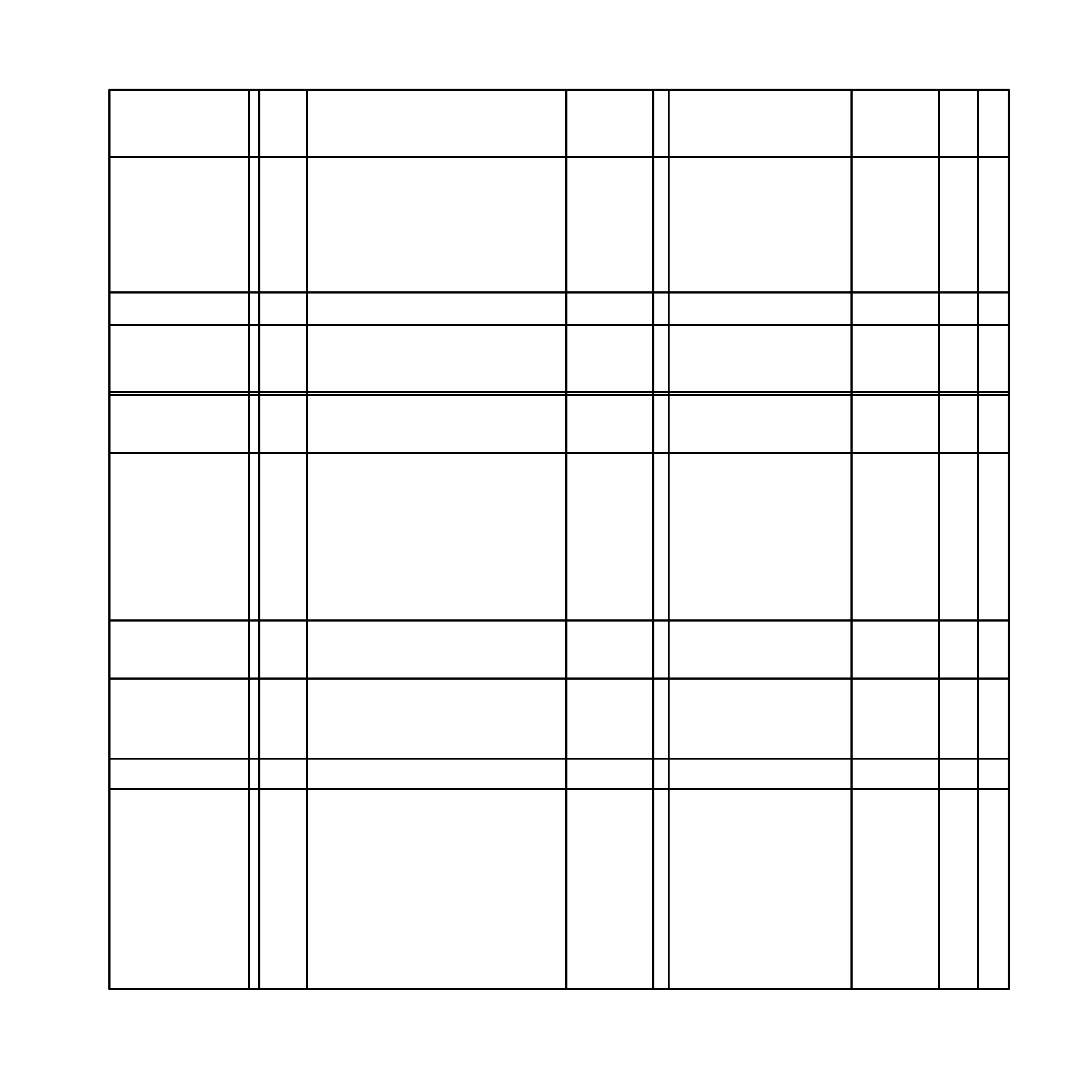}\includegraphics[width=3.7cm]{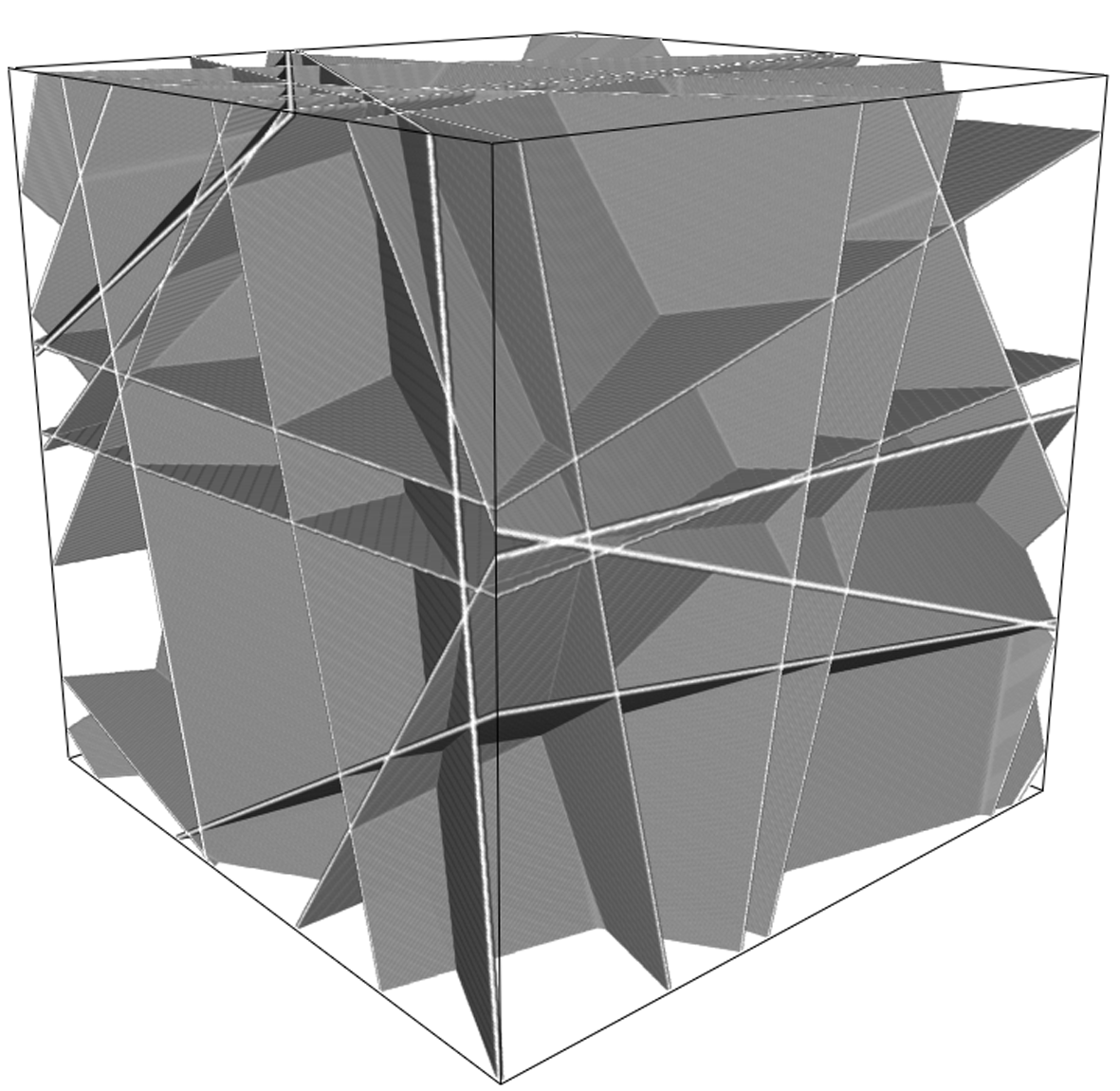}
\caption{\label{fig:vis_Hyperplane} Isotropic Poisson line tessellation in $\R^2$ (left) and Poisson hyperplane tessellation in $\R^3$ (right). Poisson line tessellation with direction distribution concentrated on the coordinate directions (middle).}
\end{center}
\end{figure}

In the stationary case, this model is parametrized by the intensity of the Poisson hyperplane process (which corresponds to $\mu_{d-1}$) and the distribution $\Rc$ of normal directions of the hyperplanes. Realizations of Poisson hyperplane tessellations are shown in Fig.~\ref{fig:vis_Hyperplane}.

Note that Poisson hyperplane tessellations are face-to-face but not normal. In the planar case, all vertices are X-vertices. This is due to the fact that almost surely not more than two lines of a Poisson line process intersect in a single point. 

First analytic results for Poisson line tessellations were obtained by Miles \cite{Miles1964,Miles1964b,MILES1973}.
An overview of mean value formulas as well as second order and distributional results for Poisson hyperplane tessellations and their typical cell can be found in \cite{matheron1974random}, \cite[Section 10.3]{SchWei08}, \cite[Section 9.5]{skm13}, and \cite{HugSchneider2024}.

Here, we only state some mean value relations for the planar case.

\begin{theorem}
Consider an isotropic Poisson line tesselation in $\R^2$.
Let $\rho= \frac{2 L_A}{\pi}$ the mean number of lines intersected by a test line segment of length 1. 
Then, 
$$L_1= \frac{1}{\rho}, \quad A_2 =\frac{4}{\pi \rho^2}, \quad N_{02} = 4,$$  
$$\gamma_0 = \frac{\pi \rho^2}{4}, \quad \gamma_1 = \frac{\pi \rho^2}{2}, \mbox{ and } \gamma_2 = \frac{\pi \rho^2}{4}.$$
\end{theorem}

Mean value formulas for anisotropic Poisson line tessellations in $\R^2$ as well as mean value formulas for spatial Poisson hyperplane tessellations can be found in \cite[p. 373-375]{skm13}.

\section{Voronoi tessellations and their generalizations}
\label{sec:VoronoiTess}

The construction of hyperplane tessellations is based on cutting space by lines/planes. In Voronoi tessellations, the initial input is a set of points. The cell of a point $x$ is then given by all points $y$ in space that have $x$ as their nearest neighbour in the given set. Depending on the metric that is chosen to measure distance and the arrangement of the generator points, very different types of cell systems can be generated. An overview of results for Voronoi tessellations is given in the textbooks 
\cite{AurenhammerBook}, \cite{mol94} and \cite{OkaBooSugChi00}. 
Voronoi tessellations can be defined for deterministic generator sets $\phi$. To obtain random tessellations, $\phi$ is chosen as a realization of a suitable point process $\Phi$.

\subsection{Voronoi tessellations}

\begin{definition}[Voronoi tessellation]
Let $\phi$ be a locally finite subset of $\R^d$. The \textit{Voronoi cell} $C(x,\phi)$ of $x \in \phi$ is defined as 
\begin{equation}
\label{eq:Voronoi}    
C(x,\phi) = \{ y \in \R^d \, :\, ||y-x|| \le ||y-x' || \quad \text{for all } x' \in \phi\}, \end{equation}
where $||\cdot||$ denotes the Euclidean norm.
The {\it Voronoi tessellation} of $\phi$ is the set $V(\phi) = \{ C(x,\phi) \, :\, x \in \phi \}.$
\end{definition}

The Voronoi tessellation is also known as the \emph{Dirichlet} or the \emph{Thiessen tessellation}. The cells of a Voronoi tessellation are convex polyhedra. Cell facets are located on hyperplanes that orthogonally intersect the line segment connecting two generating points in its midpoint. Voronoi tessellations are face-to-face. They are normal, if the points of $\phi$ are in \emph{general position}, i.e., for any $k= 2, \ldots, d$ no $k+1$ points are contained in a $(k-1)$ dimensional affine subspace of $\R^d$ and 
no $d+2$ points are contained in the surface of a sphere.

\begin{figure}[t]
\begin{center}
\includegraphics[width=3.7cm]{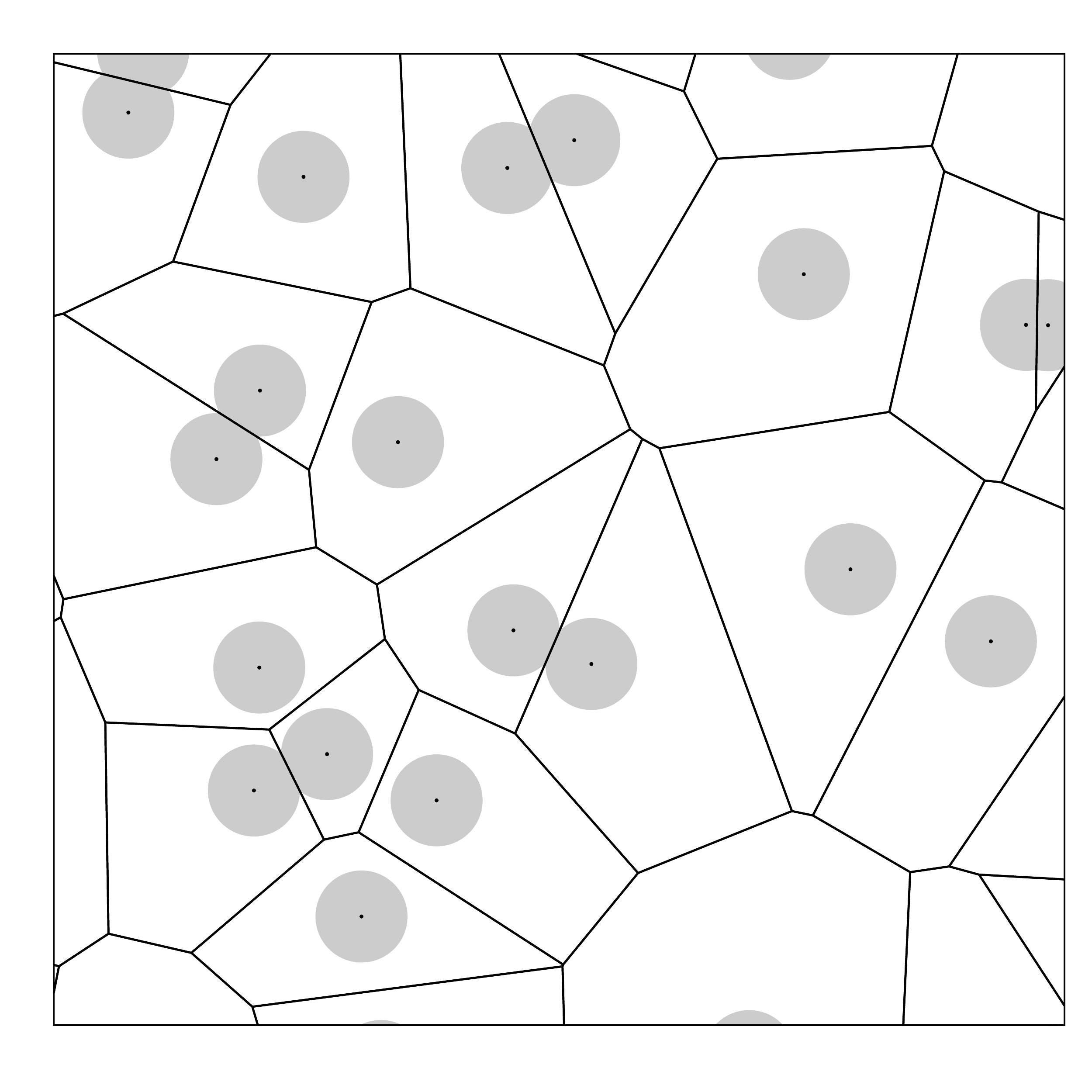}  \includegraphics[width=3.7cm]{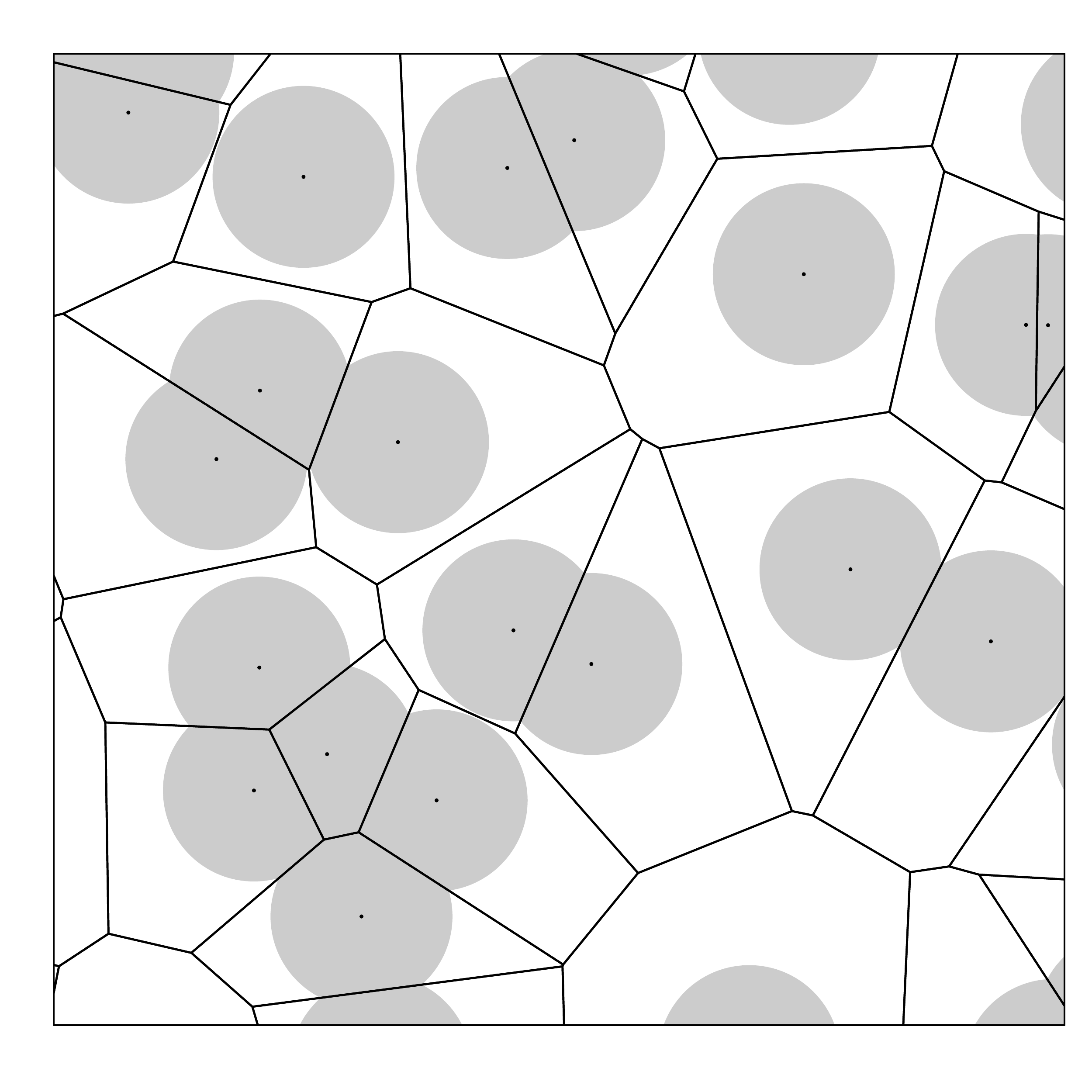} \includegraphics[width=3.7cm]{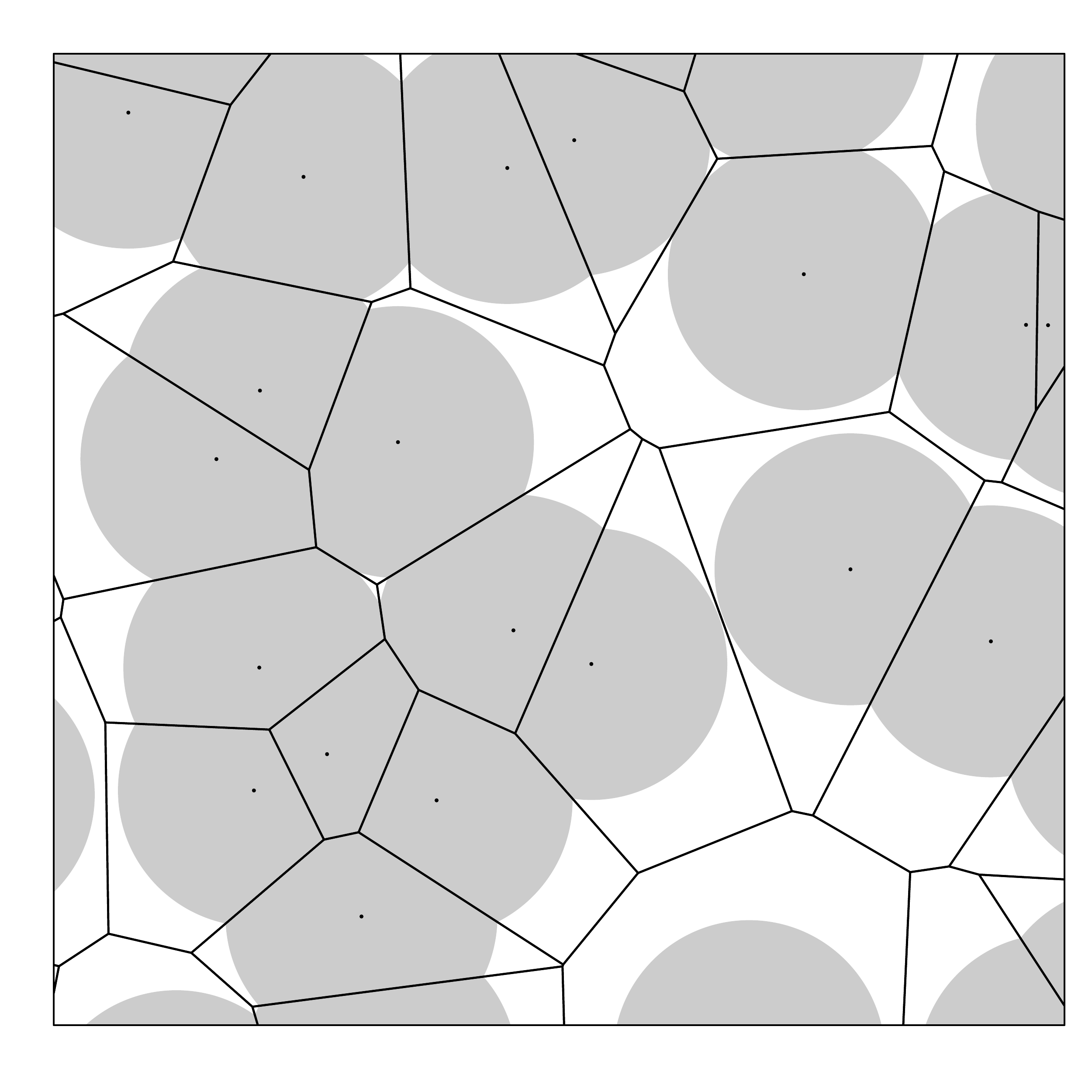}
\caption{Illustration of the Voronoi growth process at three equidistant time points.\label{fig:vis_VoronoiGrowth}}
\end{center}
\end{figure}

Voronoi tessellations can also be interpreted as the result of a growth or crystallization process. At time 0, growth starts with constant and equal speed in all generator points resulting in systems of balls of equal radius centered in the generator points. Each point $y$ in space is assigned to the generator $x$ that reached $y$ first. See Fig.~\ref{fig:vis_VoronoiGrowth} for an illustration. 

\begin{figure}[b]
\begin{center}
\includegraphics[width=3.7cm]{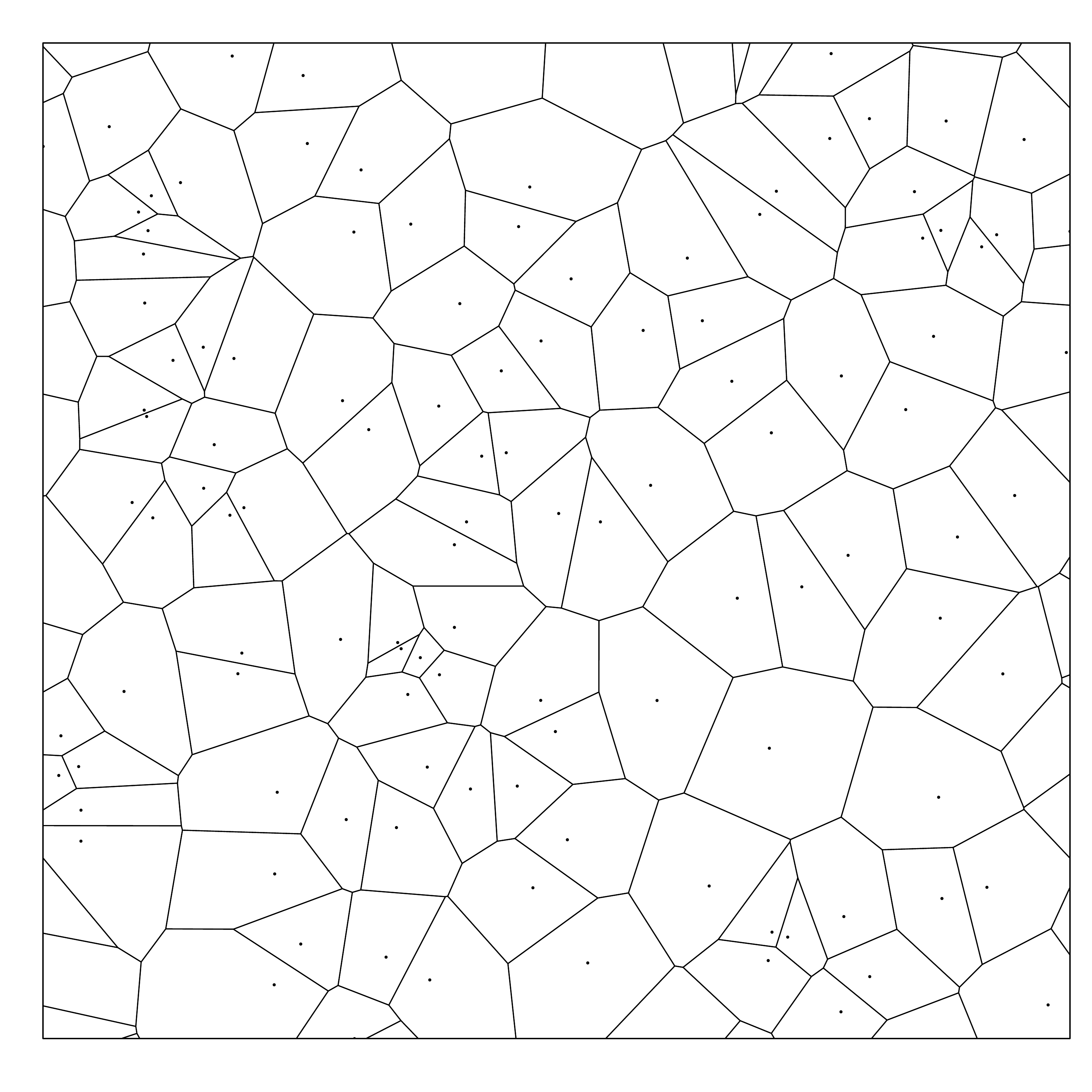}  \includegraphics[width=3.7cm]{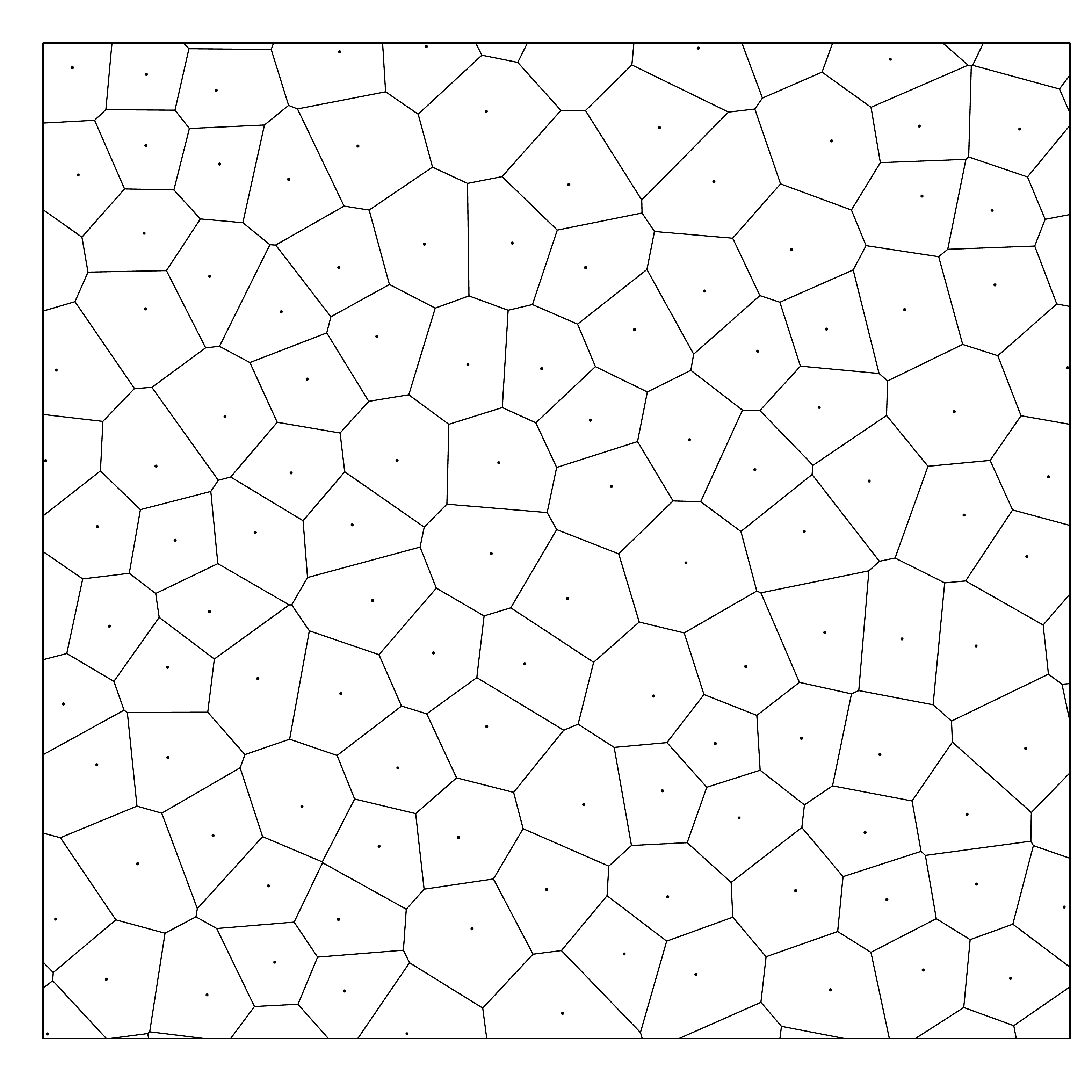} \includegraphics[width=3.7cm]{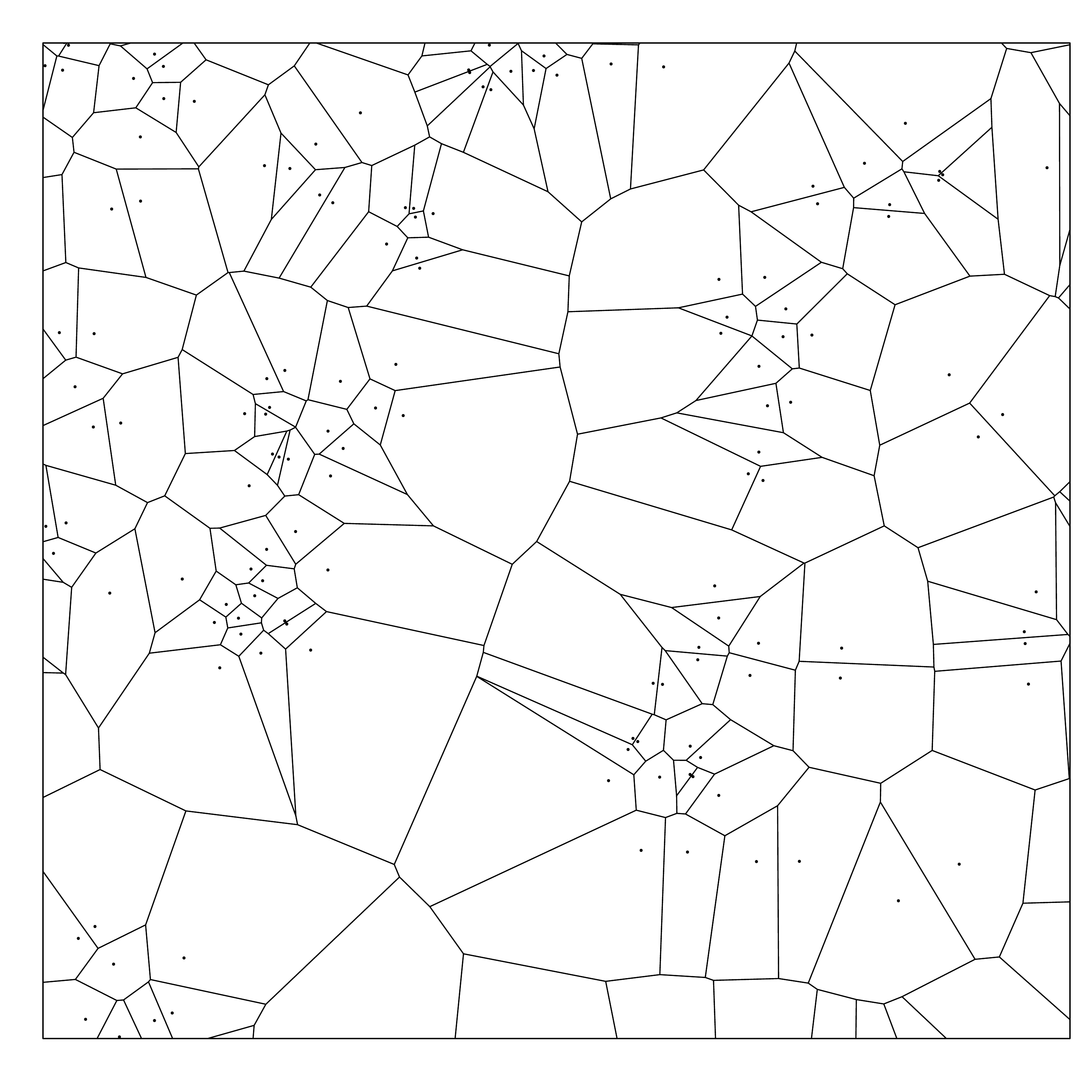}\\
\includegraphics[width=3.7cm]{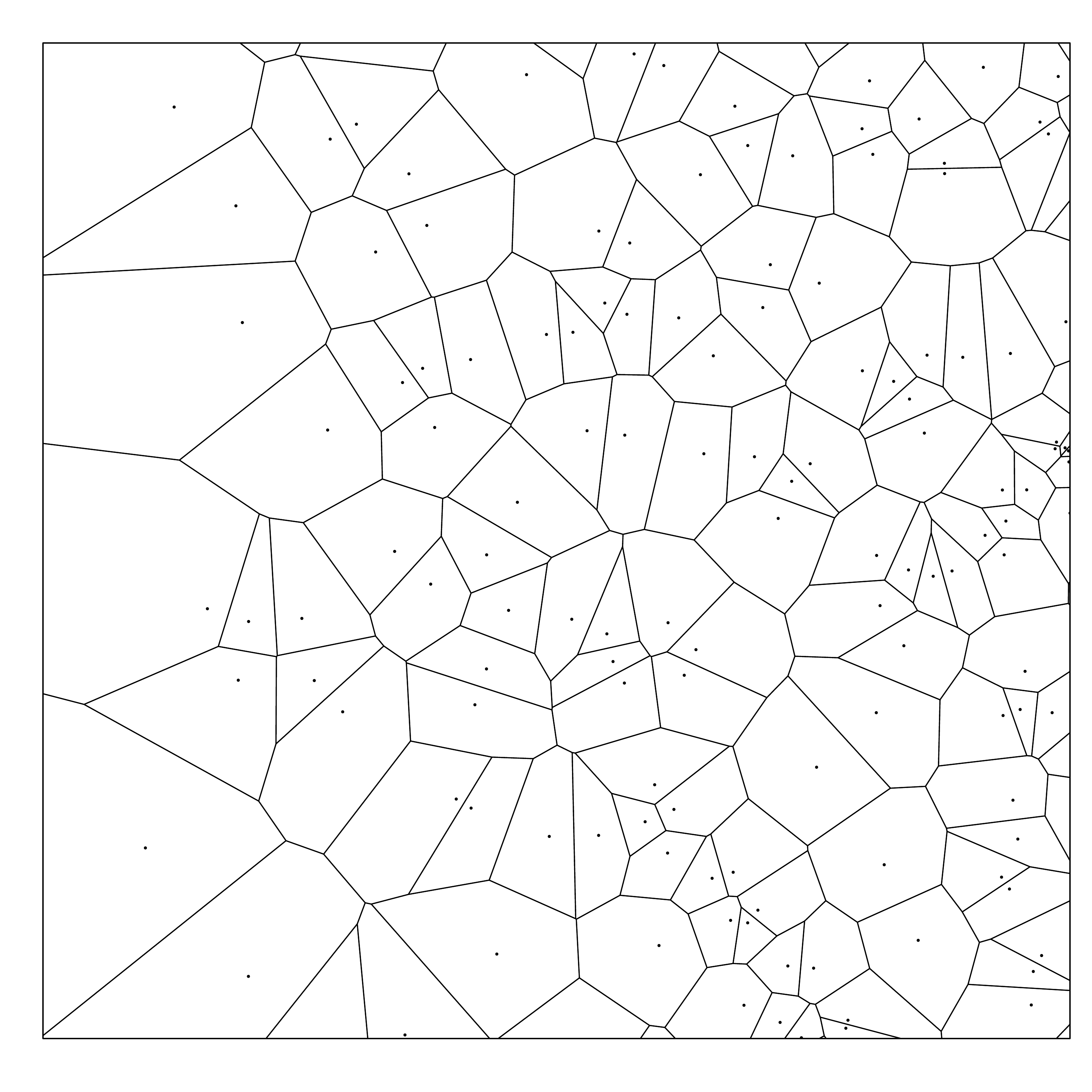}  \includegraphics[width=3.7cm]{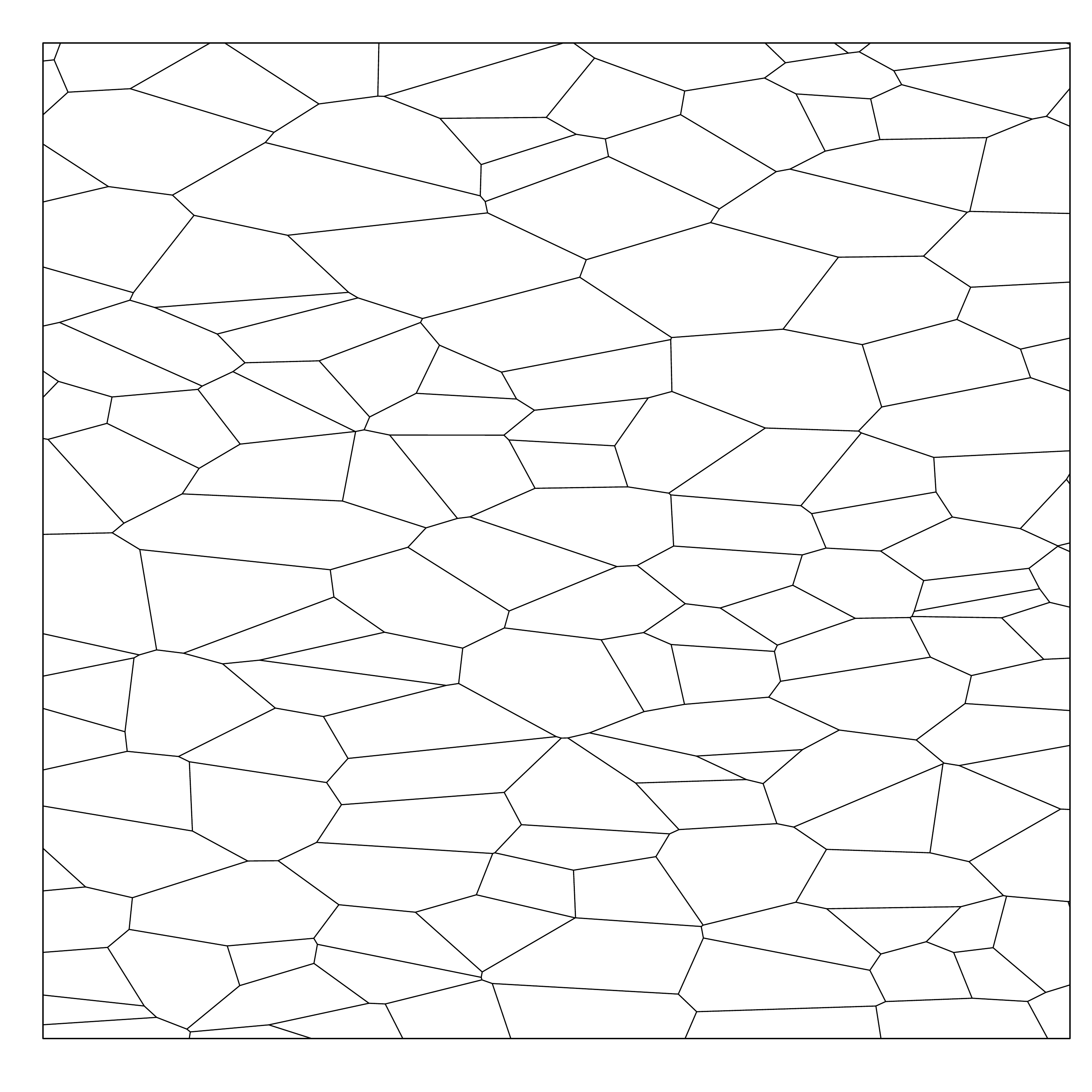} \includegraphics[width=3.7cm]{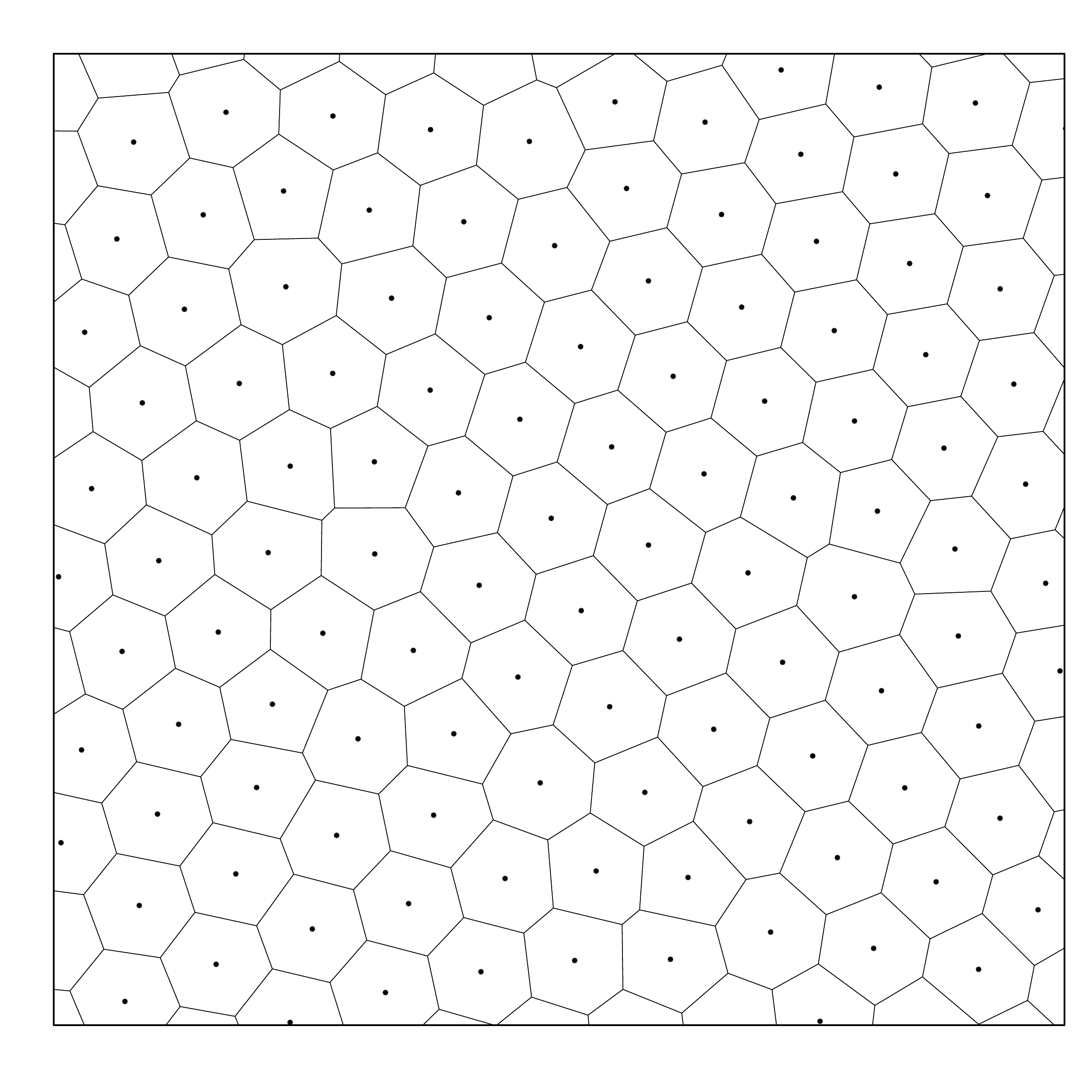}
\caption{\label{fig:vis_Voronoi} Voronoi tessellations in $\R^2$. Top from left to right: Stationary tessellation with generators drawn from a stationary Poisson process, an SSI hardcore process and a Matérn cluster process, see \cite{Spatstat} for a description of the models. Bottom from left to right: Nonstationary Poisson-Voronoi tessellation, anisotropic Voronoi tessellation by global scaling, centroidal Voronoi tessellation.}
\end{center}
\end{figure}

The cell intensity $\gamma_d$ of a stationary Voronoi tessellation is inherited from the intensity of the generating point process. The distributions of cell shapes and sizes are highly influenced by the interaction between the points. Additionally, stationarity and isotropy of the tessellation are inherited from the properties of the generating point process. Some examples are shown in Fig.~\ref{fig:vis_Voronoi}. 

The range of cell patterns that can be generated by Voronoi tessellations can be extended by transformations of the cell system. For instance, anisotropy can be introduced by scaling the tessellation with different factors in the coordinate directions. Tessellations with very regular cells are obtained from the \emph{centroidal Voronoi tessellation} \cite{Centroidal}. This notion refers to a Voronoi tessellation where the generating points coincide with the centroids (centers of mass) of the cells. An example is also shown in Fig.~\ref{fig:vis_Voronoi}.

Tessellations generated by Poisson processes allow for analytical results. In the stationary case, all mean values of the tessellation are completely determined by the intensity $\lambda$ of the generating Poisson process.

\begin{theorem}{\cite[Theorem 7.2]{Mol89}}
Let $X$ be a stationary Poisson-Voronoi tessellation with intensity $\lambda$ in $\R^d$.
The parameters $\mu_k$ are given by
$$\mu_k= \lambda^{\frac{d-k}{d}} \frac{2^{d-k+1} \pi^{\frac{d-k}{2}}}{d(d-k+1)!} \frac{\Gamma\left(d-k+\frac{k}{d}\right) \Gamma\left(\frac{d^2-dk+k+1}{2}\right)  \Gamma\left( 1+ \frac{d}{2}\right)^{d-k+\frac{k}{d}}  }{ \Gamma \left( \frac{k+1}{2} \right)  \Gamma \left( \frac{d^2-dk+k}{2} \right) \Gamma \left(\frac{d+1}{2} \right)^{d-k}} .$$
Together with $\gamma_d=\lambda$ this allows for a complete mean value characterization in $\R^2$ and $\R^3$.

For $d=2$, we have
\begin{eqnarray*}
\gamma_0= 2 \lambda, \quad \gamma_1 = 3 \lambda, \quad\gamma_2= \lambda  \\
N_{02}=3, \quad  N_{20}=6\\
L_1= \frac{2}{3 \sqrt{\lambda}}, \quad P_2= \frac{4}{\sqrt{\lambda}}, \quad A_2= \frac{1}{\lambda}.
\end{eqnarray*}

For $d=3$, the mean values are
\begin{eqnarray*}
\gamma_0&=& \frac{24 \pi^2}{35} \lambda, \quad  \gamma_1 = \frac{48 \pi^2}{35} \lambda, \quad \gamma_2 = \left( \frac{24 \pi^2}{35} +1 \right) \lambda, \quad \gamma_3 = \lambda\\
\mu_1&=& \frac{16}{15} \left(\frac{3}{4}\right)^{1/3} \pi^{5/3} \Gamma (4/3) \lambda^{2/3} \approx 5.832 \lambda^{2/3}\\
\mu_2&=& 4 \left(\frac{\pi}{6}\right)^{1/3}  \Gamma (5/3) \lambda^{1/3} \approx 2.910 \lambda^{1/3}\\
N_{21}&=& \frac{144\pi^2}{24\pi^2+35} \approx 5.23, \quad  N_{30}= \frac{96 \pi^2}{35} \approx 27.07\\
N_{31} &=& \frac{144 \pi^2}{35} \approx 40.61, \quad N_{32}= \frac{48 \pi^2}{35}+2 \approx 15.54\\
L_1&=& \frac{7 \Gamma(1/3)}{9(36 \pi)^{1/3}} \lambda^{-1/3}, \quad  A_2= \frac{35 ^\cdot 2^{8/3}\Gamma(2/3) \pi^{1/3}}{(24 \pi^2+35) 9^{2/3}} \lambda^{-2/3},\\P_2&=& \frac{7 \cdot 2^{10/3}\Gamma(1/3) \pi^{5/3}}{(24 \pi^2+35)9^{1/3}} \lambda^{-1/3}\\
A_3&=& \lambda^{-1} \quad  S_3= \left(\frac{256 \pi}{3}\right)^{1/3}  \Gamma (5/3) \lambda^{-2/3} \approx 5.821 \lambda^{-2/3}\\B_3&=&\frac{1}{5} \left(\frac{16 \pi^5}{243}  \right)^{1/3}  \Gamma (1/3)\lambda^{-1/3} \approx 1.458 \lambda^{-1/3}
\end{eqnarray*}
\end{theorem}

A tabulated overview can be found in \cite{OkaBooSugChi00}.

\begin{figure}[t]
\begin{center}
\includegraphics[width=3.7cm]{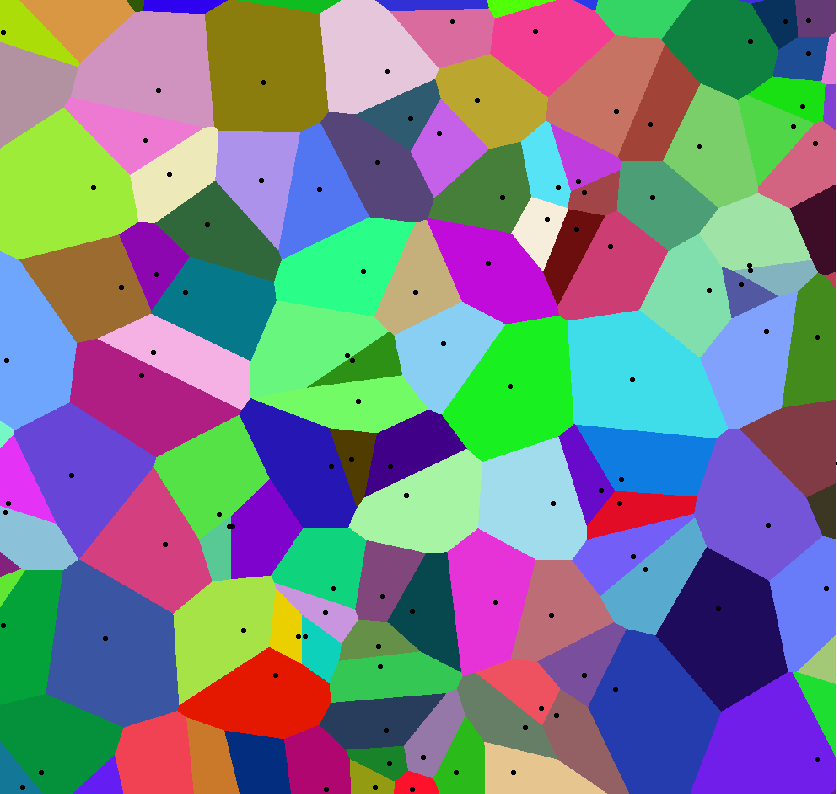}  \includegraphics[width=3.7cm]{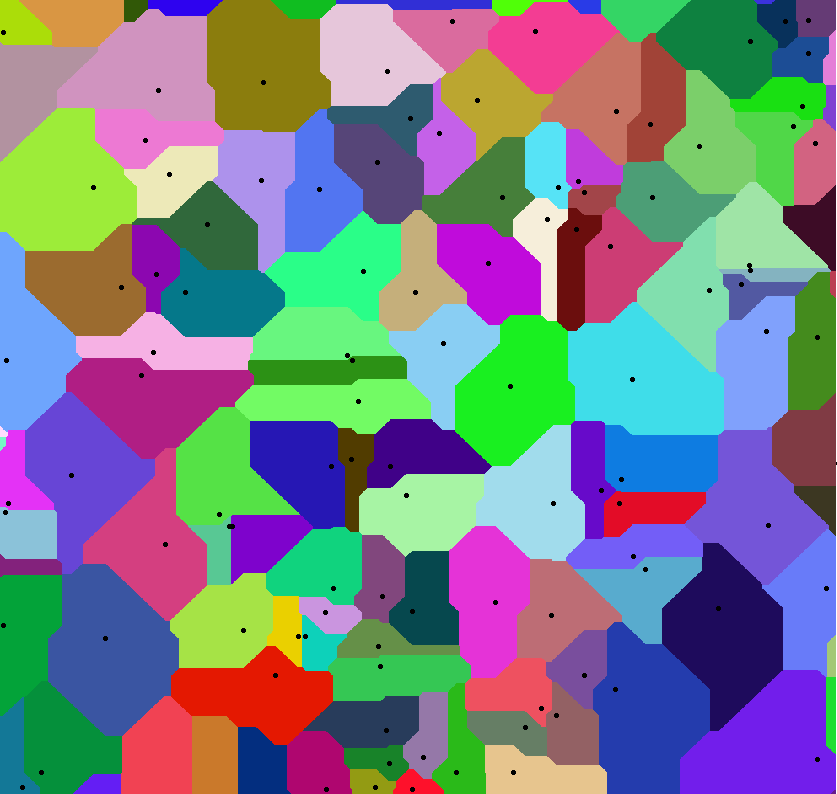} \includegraphics[width=3.7cm]{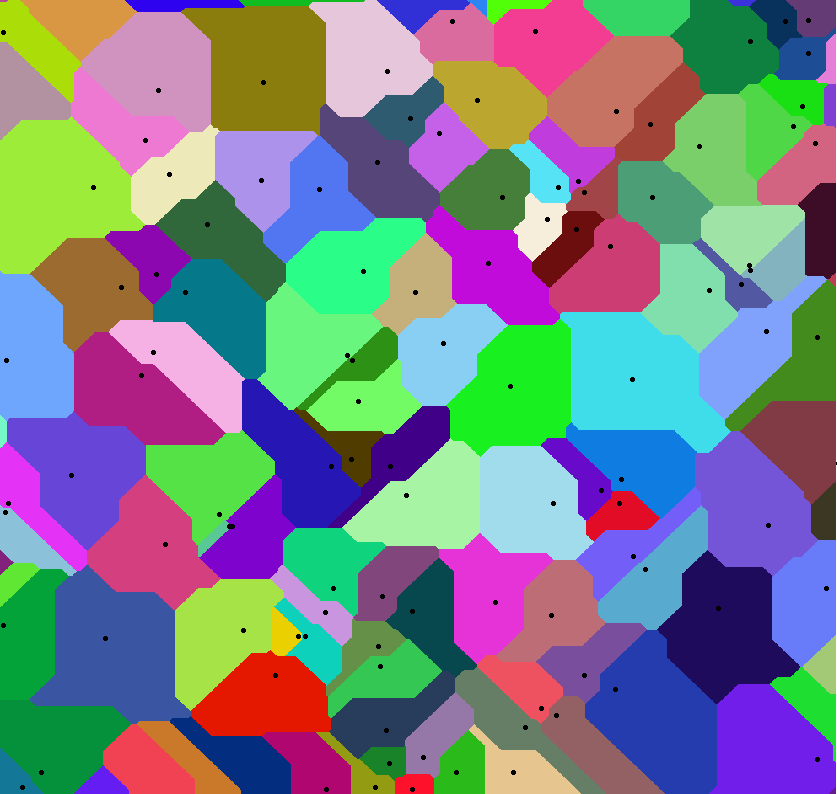}
\caption{\label{fig:VoronoiOtherMetrics} Voronoi tessellations in $\R^2$ when using different metrics: Euclidean metric $d(x,y)=((x_1-y_1)^2+(x_2-y_2)^2)^{1/2}$ based on 2-norm (left), Manhattan metric $d(x,y)=|x_1-y_1|+|x_2-y_2|$ based on 1-norm (middle), metric $d(x,y)=max(|x_1-y_1|,|x_2-y_2|)$ based on max-norm (right). The same realization of a Poisson process is used to generate all three images.}
\end{center}
\end{figure}

It is most common to consider the Euclidean metric in \eqref{eq:Voronoi} when defining Voronoi cells. However, other metrics can also be used. Fig.~\ref{fig:VoronoiOtherMetrics} shows planar Voronoi tessellations with respect to the Euclidean metric, the Manhattan metric and the maximum metric.

Another modification of the Voronoi construction is to consider diagrams of higher order \cite{edelsbrunner_poissondelaunay_2019}. In the original Voronoi tessellation (order 1), cells are formed by all points having the same nearest neighbour among the generators. In the Voronoi tessellation of order $k$, cells are defined by all points having the same $k$ nearest neighbours in the generator set.

\subsection{Weighted Voronoi tessellations}

In the Voronoi tessellation, all generator points have an equal weight. The cells can be interpreted as the result of a growth process where each point appears at the same time and all points grow with the same constant speed. 

In practice, it is often desirable to assign weights to the points which make the points appear at different times or let some cells grow faster than others. 

In Voronoi tessellations, the variation in cell structures can only be controlled by the choice of the generating point process, see Fig.~\ref{fig:vis_Voronoi}. By addition of the weight, much more flexibility in the range of cell systems is obtained. 

The way in which the weights are incorporated in the distance metric influences the shape of the resulting cells. For some models, cells are no longer convex. In addition, the speed of growth can depend on directions such that the cell systems are no longer isotropic.

\begin{figure}[b]
\begin{center}
\includegraphics[width=2.8cm]{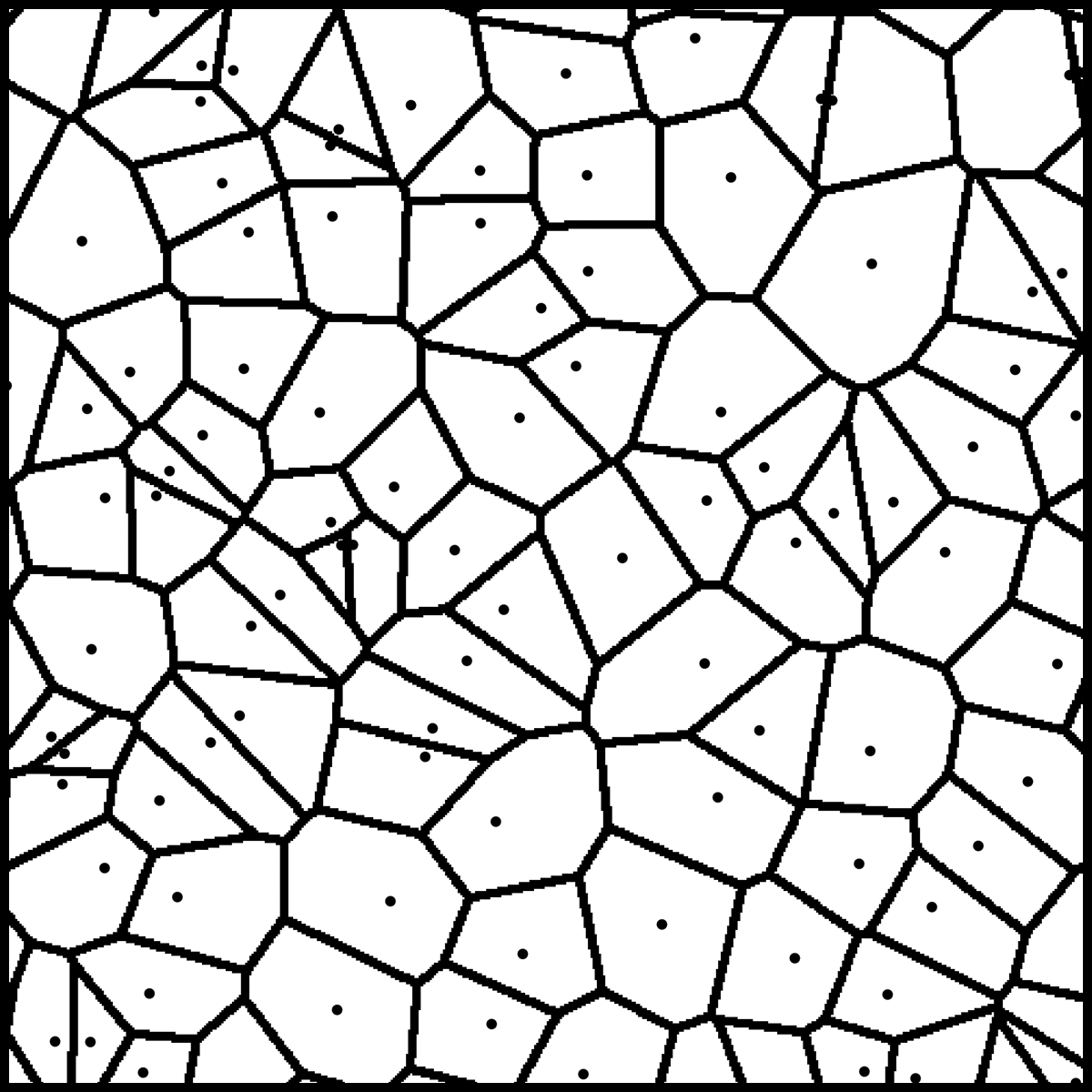}  \includegraphics[width=2.8cm]{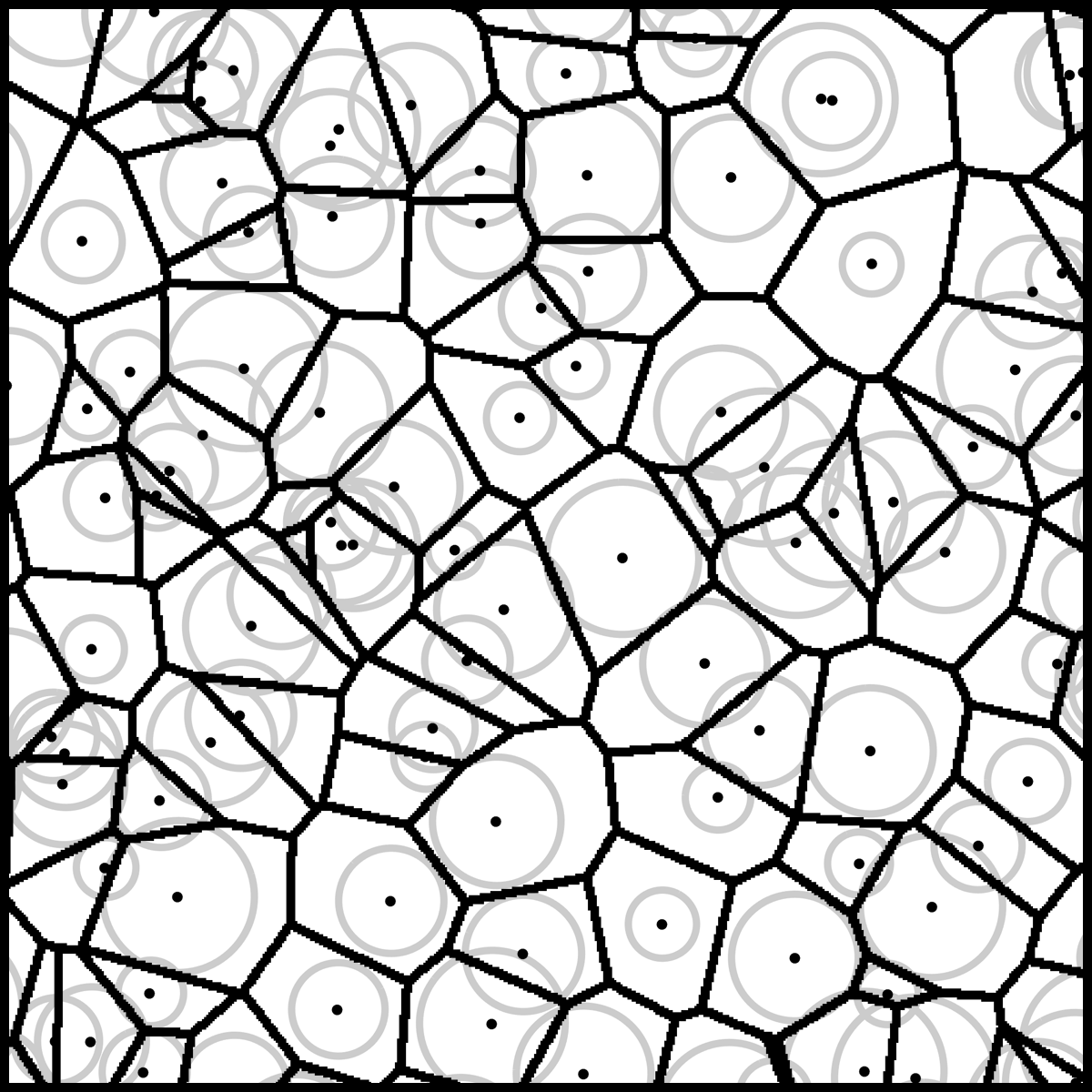} \includegraphics[width=2.8cm]{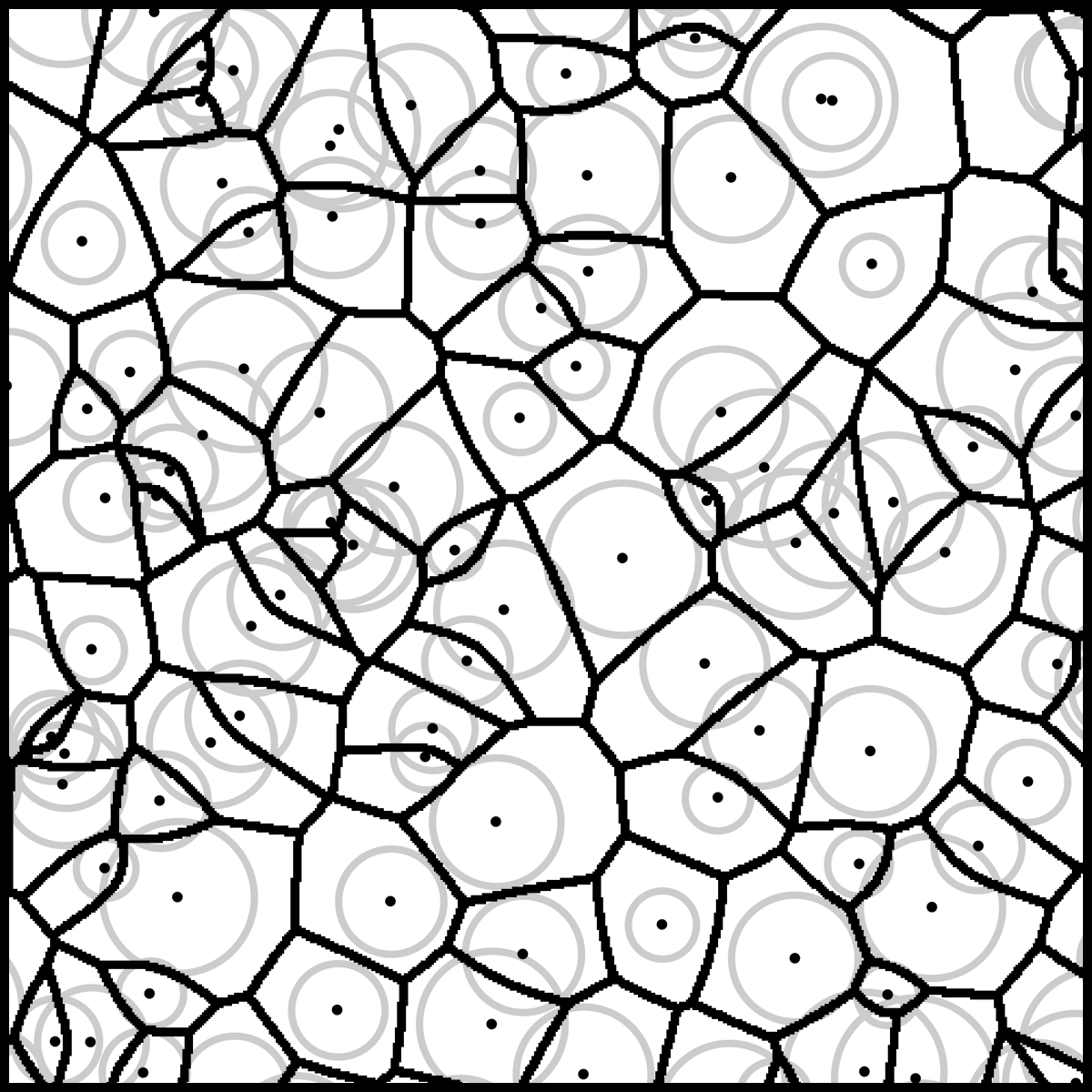}
\includegraphics[width=2.8cm]{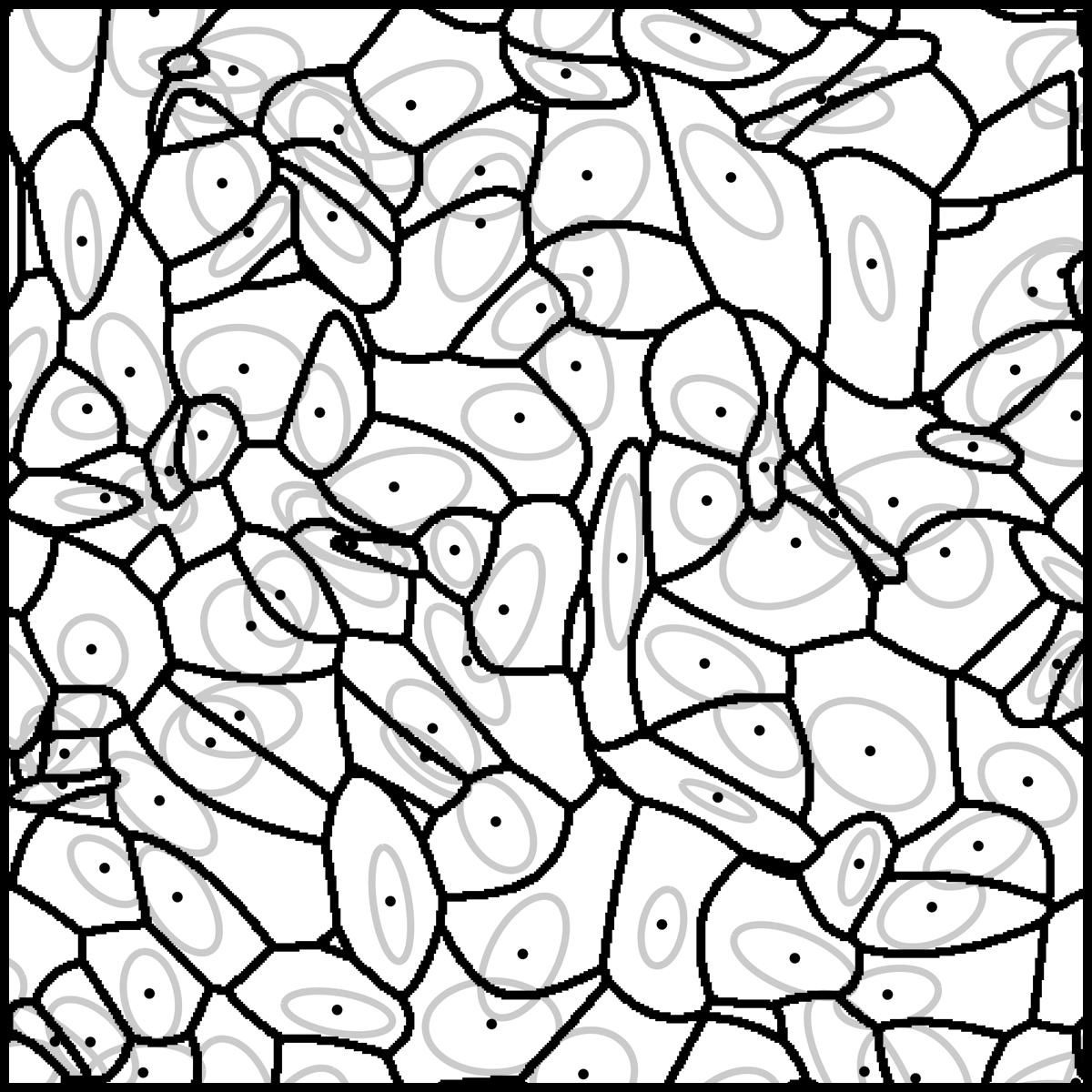} 
\caption{\label{fig:vis_WeightedVoronoi} Several weighted Voronoi tessellation models in $\R^2$. The generating point process is a Poisson process of intensity $\lambda=100$. The window is the unit square. From left to right: Voronoi tessellation, Laguerre tessellation with radii from a uniform distribution on the interval [0.025,0.075], Johnson-Mehl tessellation of the same point pattern, GBPD generated from ellipses with centers as above and semi-major axes lengths uniformly distributed on [0.4,0.7], semi-minor axes lengths on [0.1,0.4] and rotations on [0,$\pi$]. The weights of the GBPD are the same as the radii of the Laguerre tessellation.}
\end{center}
\end{figure}

An overview of weighted distance metrics that have been suggested in the literature is given in Table~\ref{Tab:Metrics}. Unfortunately, the nomenclature for the resulting models is not unified such that most models are known under several names. Visualizations of realizations of the models in $\R^2$ are shown in Fig.~\ref{fig:vis_WeightedVoronoi}. A general concept of Voronoi tessellations with respect to local metrics is presented in \cite{JEULIN2014139,Jeulin2021}.

In the following, the most prominent models that have also been investigated analytically are introduced in more detail.

\begin{table}[!t]
\caption{Overview of metrics for generalized Voronoi tessellations. Let $x,y \in \R^d$ with $x$ a generator point and $y$ arbitrary, $w, r >0$ weights where $r$ is interpreted as the radius of a ball centered in $x$, $K$ a compact set in $\R^d$ and $M$ a positive definite $d \times d$ matrix.}
\label{Tab:Metrics}      

\begin{tabular}{p{4cm}p{5.2cm}p{2cm}}
\hline\noalign{\smallskip}
Names & Metric & References\\
\hline\noalign{\smallskip}
\multirow{3}{*}{\minitab[l]{Voronoi  \\
Dirichlet \\ Thiessen }}& \\
\multirow{3}{*} & $d(x,y) = ||x-y||$ &\cite{Mol89,OkaBooSugChi00}\\
\multirow{3}{*} & &\\
\hline
\multirow{3}{*}{\minitab[l]{Laguerre  \\ power  \\radical }} & $d((x,w),y) = ||x-y||^2- w$ or &\\
\multirow{3}{*} & $d((x,r),y) = ||x-y||^2- r^2$ &\cite{lautensack06:_random_laguer,redenbach08random}\\ 
\multirow{3}{*} & \\
\hline
\multirow{3}{*}{\minitab[l]{Johnson-Mehl \\ Apollonius \\Additively weighted}} & $d((x,r),y) = ||x-y||- r$ &\cite{Moller1995_johnsonMehl}\\
\multirow{3}{*} & \\
\multirow{3}{*} & \\ \hline
Multiplicatively weighted & $d((x,w),y) = w ||x-y||$&\cite{aurenhammer-mwv, held2020efficient} \\ \hline
\multirow{2}{*}{\minitab[l]{Set Voronoi \\ Voronoi-S}}   &  $d(y,K) = \min \{ ||y-x||: x \in K \}$ & \cite{schaller}\\ 
\multirow{2}{*} & \\ \hline
\multirow{2}{*}{\minitab[l]{Ellipsoid \\ Anisotropic Voronoi}}  & $d((x,M),y) = (x-y)^TM(x-y)$ &\cite{altendorf,LSorig}\\ 
\multirow{2}{*} & \\ \hline
Generalized balanced power & $d((x,M,w),y) = (x-y)^TM(x-y)-w$ &\cite{linproGBPD,alpers_turning_2022,alpers_dynamic_2023,jung2023analytical}
\\
\noalign{\smallskip}\hline\noalign{\smallskip}
\end{tabular}
\end{table}

\subsubsection{Laguerre tessellations}

In Laguerre tessellations, each generator point $x$ is assigned a (usually) nonnegative weight $r$. An increase of $r$ is supposed to increase the size of the cell generated by $x$.

\begin{definition}[Laguerre tessellation]
Let $\phi$ be a locally finite subset of $\R^d \times \R_+$. The \textit{Laguerre cell} of $(x,r) \in \phi$ is defined as 
\[C((x,r),\phi) = \{ y \in \R^d \, :\, ||y-x||^2-r^2 \le ||y-x'||^2-r'^2 \quad \text{for all } (x', r') \in \phi.\} \]
The {\it Laguerre tessellation} of $\phi$ is the set $L(\phi) = \{ C((x,r),\phi) \, :\, (x,r) \in \phi \}.$ The 'distance' $ \pow((x,r),y) =||y-x||^2-r^2$ is called the {\it power} of $y$ w.r.t. $(x,r)$. 
\end{definition}

A generator $(x,r) \in \phi$ can be interpreted as a ball with centre $x$ and radius $r$. For points $y$ outside the  ball, the power distance $\pow((x,r),y)$ measures the squared length of the tangent line from $y$ to the ball. If all radii are equal, the special case of a Voronoi tessellation is obtained. 

While each cell of a Voronoi tessellation contains its generating point, this is no longer true for Laguerre tessellations (see Fig.~\ref{fig:vis_WeightedVoronoi}, left). In fact, there may be points which do not generate a cell at all. However, if the system of generators consists of nonoverlapping balls, then each ball is completely contained in its cell. Under a condition  on the set $\phi$ that generalizes the general position assumption for Voronoi tessellations, $L(\phi)$ is a normal random tessellation. An important reason to consider this model is the fact that for $d\ge 3$ each normal tessellation of $\R^d$ is a Laguerre tessellation \cite{Aurenhammer87a, redenbach08random}. 

In \cite{redenbach08random}, integral formulae for the intensities $\mu_k$ of the Laguerre tessellation generated by a stationary Poisson process $\Phi$ with intensity $\lambda$ are given. Due to the lack of symmetries in the Laguerre tessellation these formulae are less explicit than the ones for the Poisson-Voronoi tessellation:
For $m \in \N$ and $x_0,\ldots, x_m \in \R^m $ let $\Delta_m(x_0,
\ldots, x_m)$ be the $m$-dimensional volume of the convex hull of
$x_0, \ldots, x_m$ in $\R^m$. For $w_0,\ldots,w_m\ge 0$ define
\begin{displaymath}
  V_{m,k}(w_0,\ldots,w_m)= (m!)^{k+1}\int\limits_{S^{m-1}} \ldots
\int\limits_{S^{m-1}} \Delta_m^{k+1}(w_0 u_0, \ldots, w_mu_m )\,
\sigma(du_0)\ldots \sigma(du_m),
\end{displaymath}
where $\sigma$ is the surface measure on $S^{d-1}$.
Furthermore, let 
\begin{displaymath}
  p(t)= \exp \Big( - \lambda \kappa_d \int\limits_0^{\infty}
  \big([t+r^2]_+\big) ^{\frac{d}{2}}Q(dr) \Big),
\end{displaymath}
where $t_+= \max \{t,0\}$ and $\kappa_d$ is the volume of the $d$-dimensional unit ball. Then $p(t)$ is the probability that the power from the
origin to each point of $\Phi$ exceeds $t$.

\begin{theorem}{\cite[Theorem 4.3]{redenbach08random}}
Let $\Phi$ be a stationary, independently marked (marks are independent and independent of the point locations) Poisson process in $\R^d$ with intensity $\lambda$ and mark distribution $Q$ with finite $d$-th moment.
The intensities $\mu_k, 0< k<d$, of the Laguerre tessellation generated by $\Phi$ are
  given by the formula
  \begin{align}
    \begin{split}
      \label{redenbach_meancontent}
      \mu_k&=\frac{\lambda^{m+1}}{4(m+1)!} c_{dm} \sigma_{k}
      \int\limits_0^{\infty} \ldots \int\limits_0^{\infty}
      \int\limits_{-\min\limits_i r_i^2}^{\infty}
      \prod\limits_{i=0}^m (t+r_i^2)^{\frac{m-2}{2}} V_{m,k}\Big((t+r_0^2)^{\frac{1}{2}},\ldots,
      (t+r_m^2)^{\frac{1}{2}}\Big) \\&\quad
       \phantom{\frac{\lambda^{m+1}}{4(m+1)!} c_{dm} \sigma_{k}
      \int\limits_0^{\infty} \ldots \int\limits_0^{\infty}} \times \int\limits_0^{\infty} p(s+t) s^{\frac{k-2}{2}} \,ds 
       \, dt\,Q(dr_0)\ldots Q(dr_m), 
    \end{split}
  \end{align}
  where $m=d-k$, $\sigma_k$ is the surface area of $S^{k-1}$, and $c_{dm}=\frac{\sigma_{d-m+1}\ldots
    \sigma_{d}}{\sigma_{1}\ldots \sigma_{m}}$. For $k=d$ we have
  $\mu_d=1$, and for $k=0$,
  \begin{align*}
    \mu_0&= \frac{\lambda^{d+1}}{2(d+1)!} \int\limits_0^{\infty}
    \ldots \int\limits_0^{\infty}
    \int\limits_{-\min\limits_i r_i^2}^{\infty} \prod\limits_{i=0}^d
    (t+r_i^2)^{\frac{d-2}{2}}  V_{d,0}\Big((t+r_0^2)^{\frac{1}{2}},\ldots,
    (t+r_d^2)^{\frac{1}{2}}\Big) p(t) \, dt\\&
    \phantom{\frac{\lambda^{d+1}}{2(d+1)!} \int\limits_0^{\infty}
    \ldots \int\limits_0^{\infty}
    \int\limits_{-\min\limits_i r_i^2}^{\infty}} \times Q(dr_0)\ldots Q(dr_d).
  \end{align*}
\end{theorem}

For the Poisson-Voronoi tessellation, obviously $\gamma_d =\lambda$. For Poisson-Laguerre tessellations no explicit formula for $\gamma_d$ is known for general dimension. Nevertheless, for $d=2$, the cell intensity can be computed via $\gamma_2= \mu_0/2$. For $d=3$, an explicit formula for $\gamma_3$ would be of particular interest as $\gamma_3$ is one of the four parameters determining the mean value characteristics of the tessellation, see Section~\ref{Sec:RandomTess}.

\subsubsection{Johnson-Mehl tessellations}

Johnson-Mehl tessellations, also known as additively weighted Voronoi tessellations or Apollonius diagrams, are based on the same type of generator systems as Laguerre tessellations. However, the distance from a point $y$ to a ball $(x,r)$ is measured by the Euclidean distance to the ball surface rather than the tangential length. Johnson-Mehl tessellations are one example of a tessellation model with nonconvex cells.

\begin{definition}[Johnson-Mehl tessellation]
Let $\phi$ be a locally finite subset of $\R^d \times \R_+$. The \textit{Johnson-Mehl cell} of $(x,r) \in \phi$ is defined as 
\[C((x,r),\phi) = \{ y \in \R^d \, :\, ||y-x||-r \le ||y-x'||-r' \quad \text{ for all } (x', r') \in \phi.\} \]
\end{definition}

Analytical results for random Johnson-Mehl tessellations can be found in \cite{Moller1995_johnsonMehl}. As an example, we state explicit formulas for the densities $\mu_k$ in Johnson-Mehl tessellations generated by Poisson point processes.

\begin{theorem}{\cite[Theorem 4.1]{Moller1995_johnsonMehl}}
Let $\Phi$ be a stationary, independently marked Poisson process with intensity $\lambda$ and mark distribution $Q$.
For $0<k<d$ and $m=d-k$ we have 
$$
\mu_k=\lambda^{md+k} c_{dk} \int\limits_0^{\infty} \int\limits_0^t \left( (t-s)^{d-1} Q(ds)  \right)^{m+1}  \exp \left(-\lambda^d \kappa_d \int\limits_0^t (t-s)^d Q(ds) \right) dt
$$   
with
$$ c_{dk} = \frac{2^{m+1} \pi^{\frac{(m+1)d}{2}} \Gamma\left( \frac{dm+k+1}{2}\right)}{(m+1)! \Gamma\left( \frac{dm+k}{2}\right) \Gamma\left( \frac{d+1}{2}\right)^m \Gamma\left( \frac{k+1}{2}\right)} .
$$
\end{theorem}

\subsubsection{Generalized balanced power diagrams}

In addition to the weighting, the generalized balanced power diagram allows to choose individual directions of elongation of the cells. To this end, the weights of the points are extended by positive definite $d\times d$ matrices. Let $\mathcal{M}$ denote the set of all such matrices.

\begin{definition}[Generalized balanced power diagram]
Let $\phi$ be a locally finite subset of $\R^d \times \mathcal{M} \times \R_+$. The cells of the \textit{generalized balanced power diagram (GBPD)} generated by $\mathbf{x}:= (x,M,r) \in \phi$ are  defined as 
\[C(\mathbf{x},\phi) = \{ y \in \R^d  :(y-x)^TM(y-x)-r \le (y-x')^TM'(y-x')-r' \,\text{ for all } \mathbf{x'} \in \phi\}. \]
\end{definition}

As Johnson-Mehl tessellations, GBPDs have nonconvex cells. However, the cells in a GBPD need not be connected. An analytic representation of the GBPD in $\R^2$ has been derived in \cite{jung2023analytical}.

The generators of a GBPD can be interpreted as ellipsoids. This is seen by decomposing the positive definite matrix $M$ into
\begin{equation}\label{decompo}
    M= U\Lambda^{-1} U^\top
\end{equation}
where $U=(u_1,\ldots,u_d)$ is an orthogonal matrix and $\Lambda= \text{diag}(a_1,\ldots,a_d)$ is a diagonal matrix. Consider the equation
\begin{equation}\label{drawellipse}
(y-x)^\top M (y-x) - w = 1.
\end{equation}
For $z=(z_1,\ldots,z_d)^T\coloneqq U^T(y-x)$, Equation \eqref{drawellipse} reads
\begin{equation}\label{drawellipse1}
\frac{z_1^2}{a_1}+\ldots+\frac{z_d^2}{a_d} - w = 1.
\end{equation}
For $w=0$, this equation defines an ellipsoid. Because $z$ is obtained by rotation or reflection of $(y-x)$, the pair $(x,M)$ can be interpreted as an ellipsoid centered at $x$ with semi-axes $u_1,\ldots,u_d$ and semi-axis lengths $\sqrt{a_1},\ldots,\sqrt{a_d}$.

When choosing $M_i=I_d$ and $w_i=r_i^2$, the Laguerre tessellation is obtained. For diagonal matrices $M_i = c \cdot I_d$ and all $w_i=0$, the GBPD equals the multiplicatively weighted Voronoi diagram.

\section{Delaunay tessellations}
\label{sec:DelaunayTess}

\begin{figure}[b]
\begin{center}
\includegraphics[width=3.7cm]{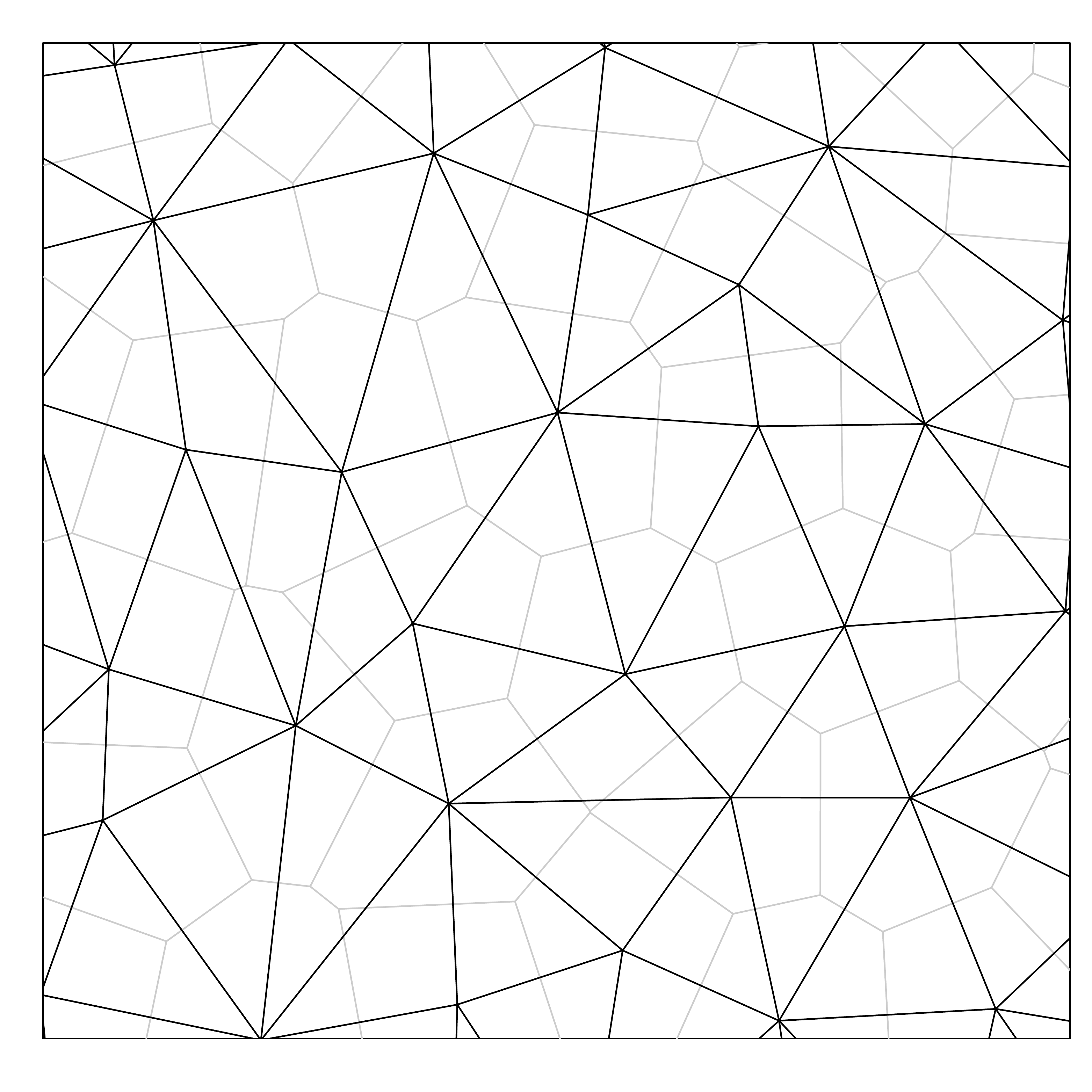} 
\caption{\label{fig:vis_DelaunayConstruction} Voronoi tessellation (grey) and its dual Delaunay tessellation (black).}
\end{center}
\end{figure}

The Delaunay tessellation is the dual of the Voronoi tessellation. Its vertices are given by the point process of generators of the Voronoi tessellation. Each facet in the Voronoi tessellation gives rise to an edge in the Delaunay tessellation which links the two generators of Voronoi cells adjacent to this facet. The cells of a Delaunay tessellation are simplices whose vertices are at equal distance to a vertex of the Voronoi tessellation. This Voronoi vertex can be used as centroid of the Delaunay cells. The construction is illustrated in Fig.~\ref{fig:vis_DelaunayConstruction}.

The Delaunay tessellation is face-to-face, but not normal. For Delaunay tessellations that are generated by a stationary Poisson process, the complete distribution of the typical cell is known. In particular, mean values of cell characteristics can be derived as functions of the intensity $\lambda$ of the generating point process. 

\begin{theorem}{\cite[Theorems 10.2.9, 10.4.4]{SchWei08}}
Consider a stationary Poisson-Delaunay tessellation with point intensity $\lambda$ in $\R^d$. 

	For $d=2$, we have
		\begin{align*}
		\gamma_0 = \lambda, \quad \gamma_1= 3 \lambda, \quad  \gamma_2= 2 \lambda \\
		L_1= \frac{32}{9 \pi \sqrt{\lambda}}, \quad A= \frac{1}{2 \lambda}, \quad P_2= \frac{32}{3 \pi \sqrt{\lambda}}
		\end{align*}

		If the Voronoi vertices are chosen as centroids of the cells, then the distribution $\mathbb{Q}$ of the typical cell is given by
		\begin{align*}
		\mathbb{Q}(A)&= a_d \lambda^d \int_0^{\infty} \int_{S^{d-1}} \ldots \int_{S^{d-1}} \one_A(\conv\{ru_0, \ldots, r u_d\} ) e^{- \gamma \kappa_d r^d} r^{d^2-1} \\&\quad \Delta_d (u_0, \ldots, u_d) \sigma(du_0) \ldots \sigma(du_d) dr 
		\end{align*}
		for $A$ a Borel set consisting of simplices with center of the surrounding balls in the origin and 
		$$
		a_d= \frac{d^2}{2^{d+1} \pi^{\frac{d-1}{2}}} \frac{\Gamma(\frac{d^2}{2})}{\Gamma (\frac{d^2+1}{2})}\left[\frac{\Gamma(\frac{d+1}{2})}{\Gamma (1+\frac{d}{2})}  \right]^d.
		$$	
\end{theorem}
Hence, the distribution $\mathbb{Q}$ of the typical cell can be simulated as follows. First, $d+1$ unit vectors are drawn with density proportional to $\Delta_d (u_0, \ldots, u_d)$. Independently, a positive number with (unnormalized) density $e^{- \gamma \kappa_d r^d} r^{d^2-1}$ is sampled. The convex hull of the rescaled vectors $\{ru_0, \ldots, r u_d\}$ then yields a sample of the typical Delaunay cell.

\begin{figure}[t]
\begin{center}
\includegraphics[width=3.7cm]{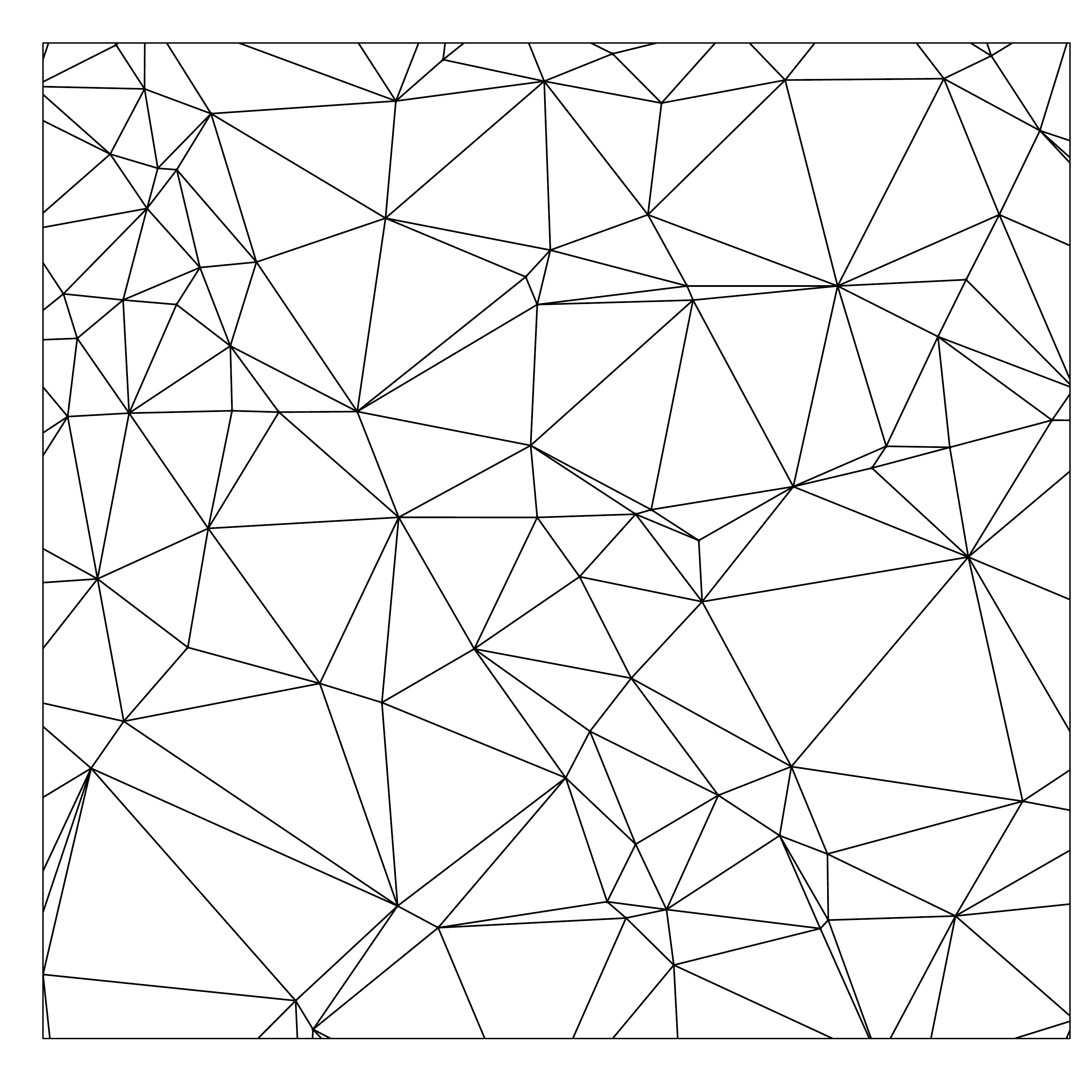} \hspace{0.1cm} \includegraphics[width=3.7cm]{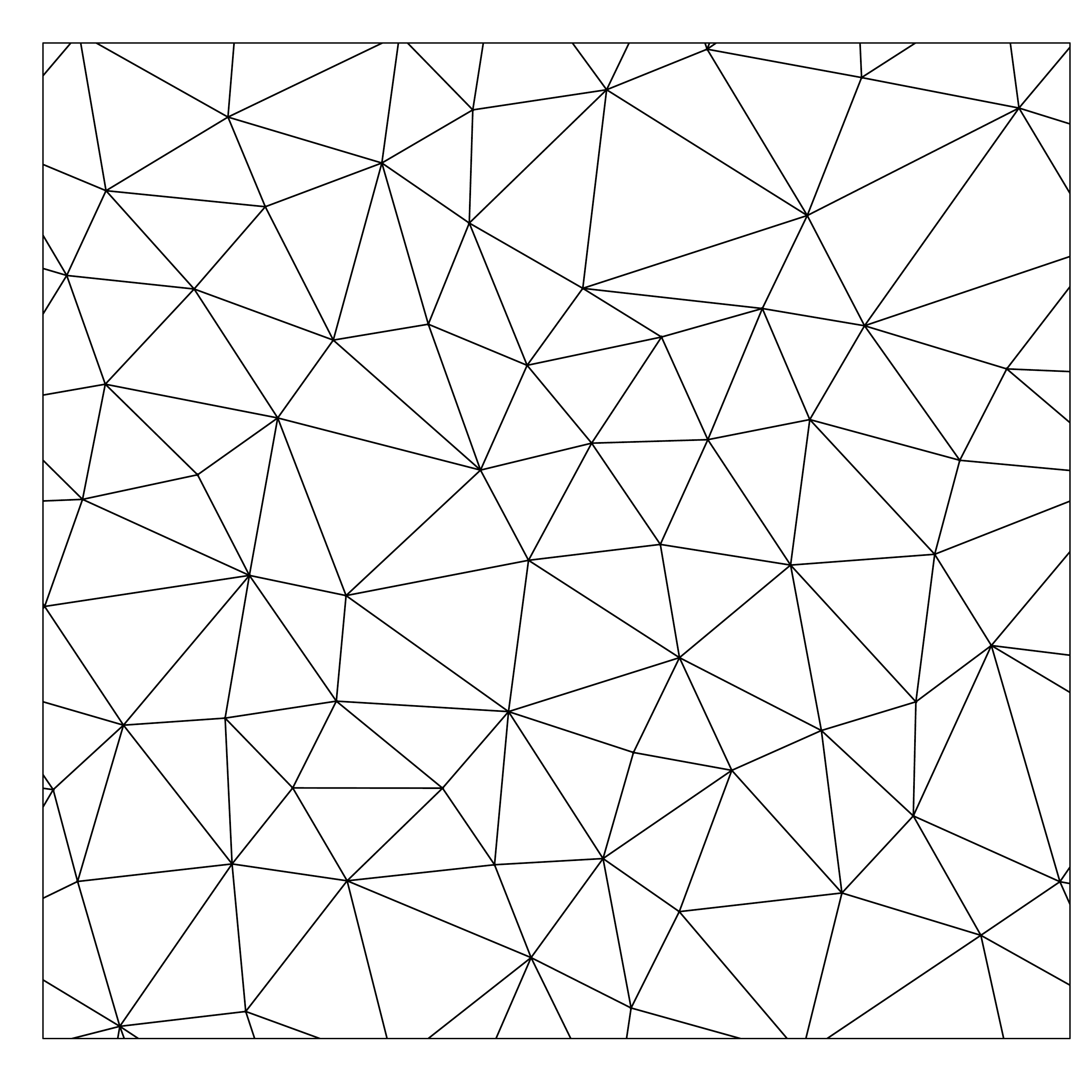} \hspace{0.1cm} \includegraphics[width=3.7cm]{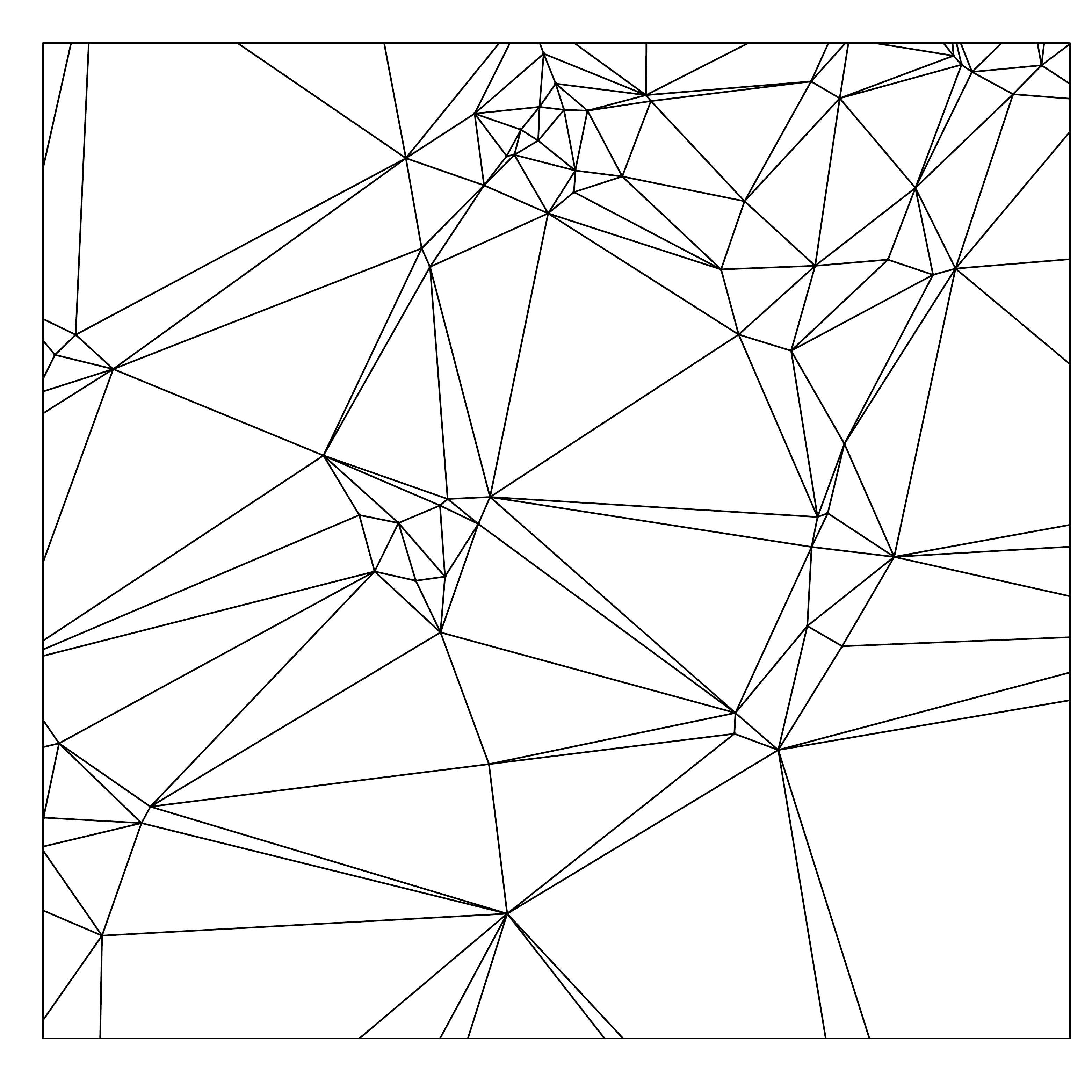}
\caption{\label{fig:vis_Delaunay} Delaunay tessellations in $\R^2$. The same models as in Fig.~\ref{fig:vis_Voronoi} are used: Generators drawn from a stationary Poisson process, an SSI hardcore process and a Matern cluster process (from left to right).}
\end{center}
\end{figure}

Realizations of Delaunay tessellations generated by a selection of different point process models are shown in Fig.~\ref{fig:vis_Delaunay}.

Generalized Delaunay tessellations can also be introduced as duals of weighted Voronoi tessellations, in particular of the Laguerre tessellation. 

 One such example is the $\beta$-Delaunay tessellation that has recently been introduced and studied analytically in a series of papers by Gusakova, Kabluchko, and Th\"ale \cite{gusakova_beta-delaunay_2021-1, gusakova_beta-delaunay_2021,  gusakova_beta-delaunay_2021-2, gusakova_-delaunay_2022}. 

The $\beta$-Delaunay tessellation is a tessellation in $\R^{d-1}$ defined by a Poisson process in $\R^{d-1} \times [0, \infty)$
whose intensity measure has density
$$
(x,h) \mapsto \gamma c_{d, \beta} h^{\beta} \quad \text{ with } c_{d, \gamma}= \frac{\Gamma\left(\frac{d}{2}+\beta+1 \right)}{\pi^{\frac{d}{2}} \Gamma(\beta+1)}.
$$
Here, $\gamma>0$ is an intensity parameter and $\beta>-1$ is a shape parameter.
The tessellation is then defined as the Delaunay tessellation with respect to the Laguerre tessellation in $\R^{d-1}$ with generator points $x$ and weights $w=-h$. 

\begin{figure}[t]
\begin{center}
\includegraphics[width=3.7cm]{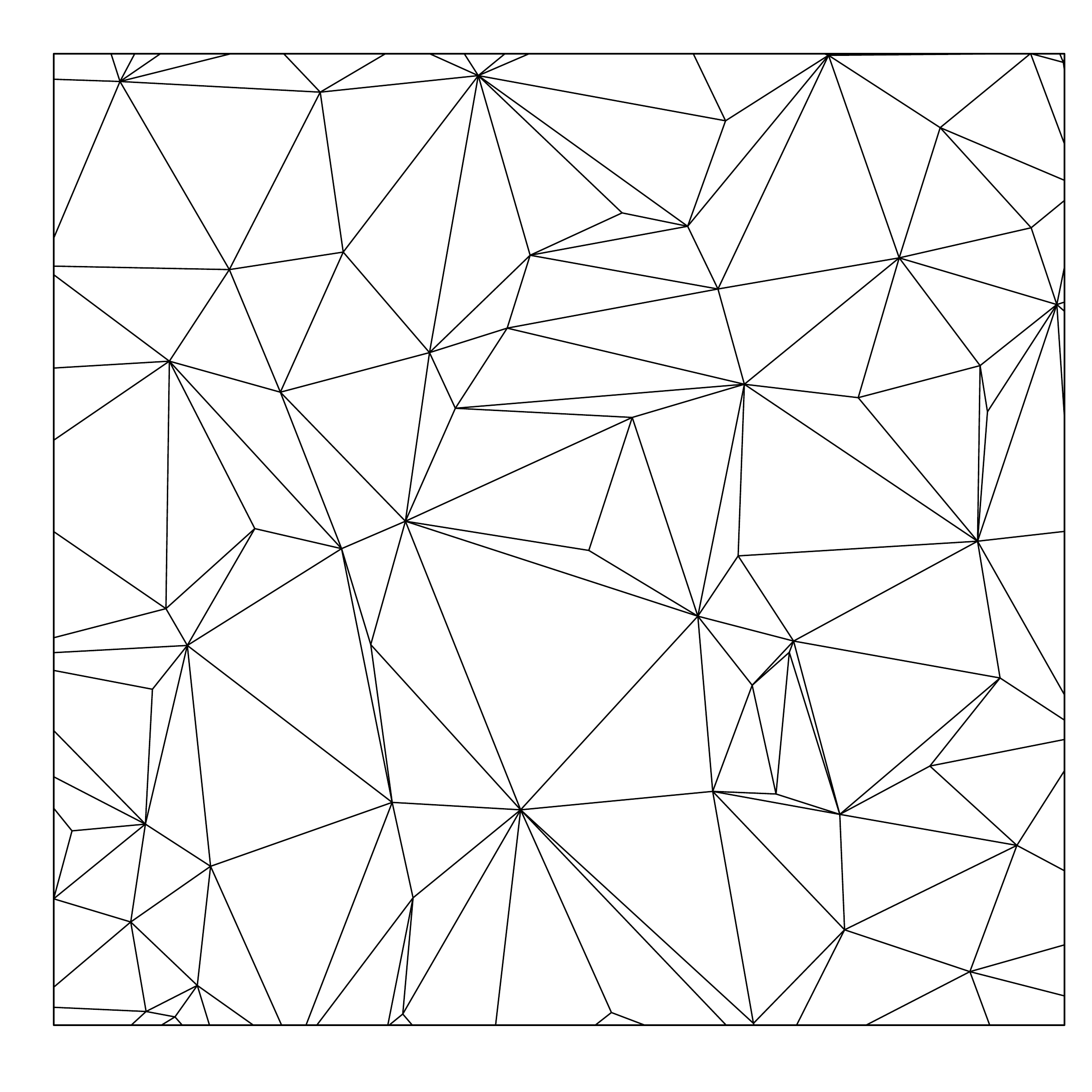} \hspace{0.1cm} \includegraphics[width=3.7cm]{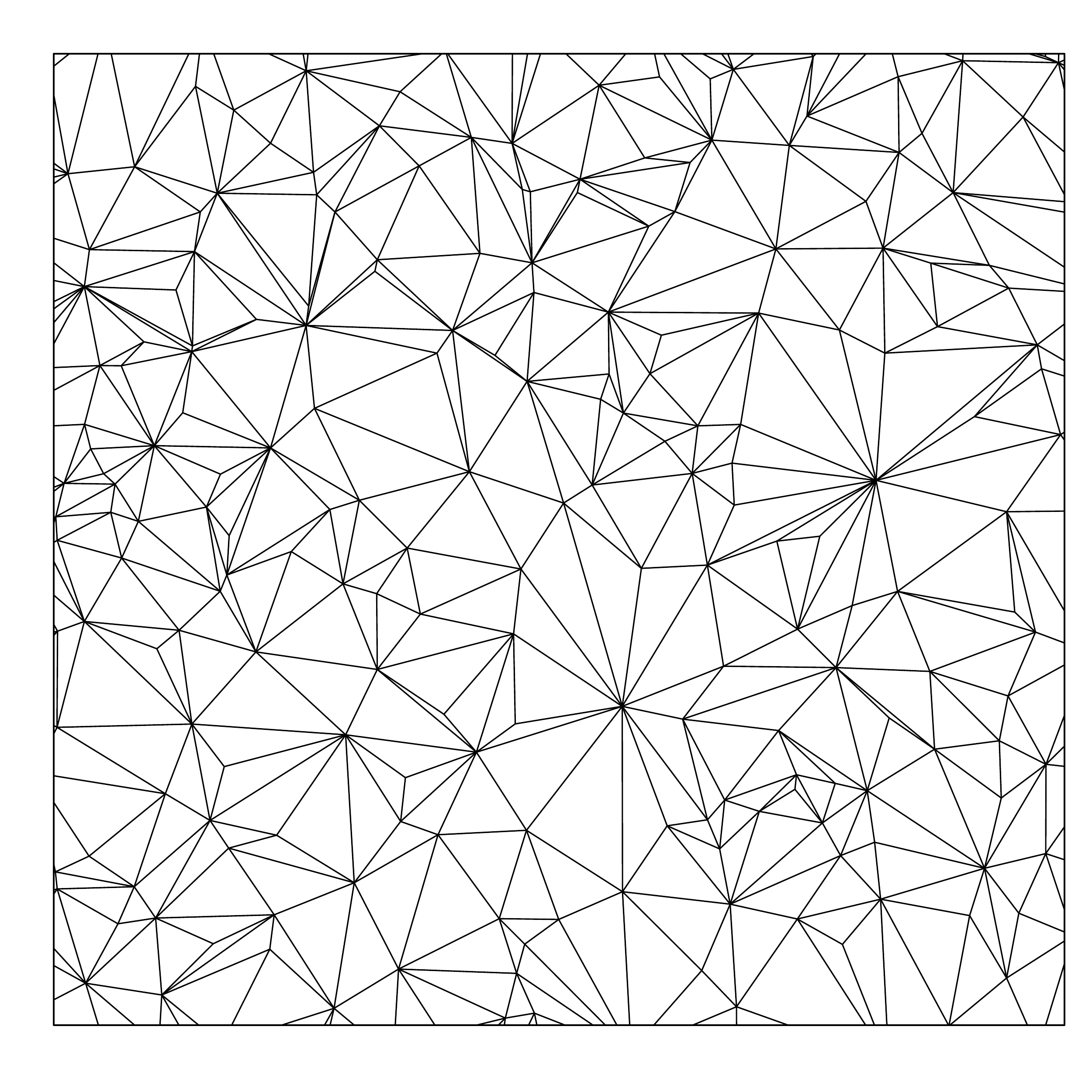} \hspace{0.1cm} \includegraphics[width=3.7cm]{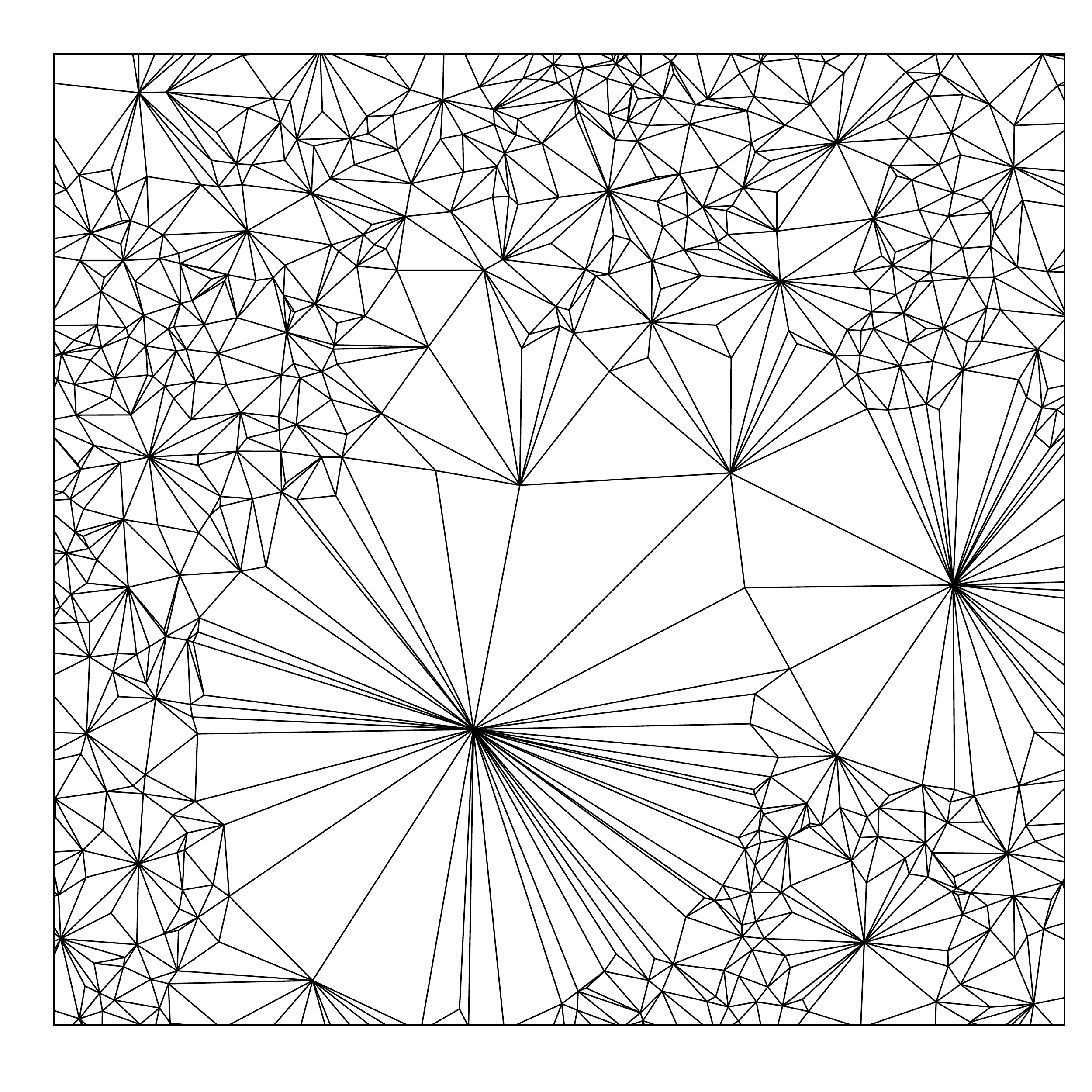}
\caption{\label{fig:vis_betaDelaunay} $\beta$- and $\beta'$-Delaunay tessellations in $[0,5]^2$ with unit intensity $\gamma=1$ and $\beta=5$ (left), $\beta=15$ (middle) and $\beta'=2.5$ (right). }
\end{center}
\end{figure}

The $\beta'$-tessellation is generated in a similar manner by a Poisson process in $\R^{d-1} \times (-\infty,0]$
whose intensity measure has density
$$
(x,h) \mapsto \gamma c'_{d, \beta} (-h)^{-\beta} \quad \text{ with } c'_{d, \gamma}= \frac{\Gamma(\beta)}{\pi^{\frac{d}{2}} \Gamma(\beta-\frac{d}{2})}.
$$
Some realizations of the model are shown in Fig.~\ref{fig:vis_betaDelaunay}.

\section{Gibbs tessellations}
Weighted Voronoi tessellations were introduced to increase the variability of cell patterns that can be generated. Additionally, the weight distribution can introduce some control on distributions of geometric characteristics of the cells, e.g. the cell sizes. An alternative way is the use of Gibbs models. These are models whose distribution is absolutely continuous with respect to the distribution of a suitable reference model, e.g., a Poisson-Voronoi or a Poisson-Delaunay tessellation. The density is used to favour or penalize certain geometric properties of the cells.

Gibbs models for point processes are well-established. A large amount of theoretical results and statistical methods for these models is available, see \cite{moel:waag04,Spatstat}.
The use of Voronoi tessellations in the density was already proposed by Ord in the discussion of \cite{Kfunction}. For the planar case, he suggested to use a density of the form
$$
f(x) = \alpha \beta^n \prod_{x_i \in \Phi} g(V(C(x_i,X))),
$$
see also \cite{BadMol}. This model is now known as Ord's process.

Gibbs-Delaunay tessellation models were studied in a series of papers by Dereudre and coauthors.  
In \cite{dereudre_gibbs_2008}, the existence of infinite Gibbs tessellation models is proven. In these models, the reference distribution is a Poisson-Delaunay tessellation. The density of a tessellation $X$ is defined locally, i.e., on a bounded window $W$ using an energy term (Hamiltonian)
\begin{align*}
E_W(X)=   \sum_{C\in X, C\text{ in } W} \phi_1(C) &+ \sum_{C_1, C_2\in X: C_1 \text{ or } C_2\text{ in } W \text{ and }  C_1 \sim C_2} \phi_2(C_1, C_2)\\
&+ \ldots + \varphi_n(C_1,\ldots,C_n)_{\big|C_1,\ldots,C_n\in X,C_1,\ldots,C_n \textrm{ in } W}\, ,
\end{align*}
where $\phi_1, \phi_2,\ldots,\varphi_n$ are the energies (potentials) of the cell $C$, the pair of cells $C_1,C_2$, and collections of up to $n$ cells $C_1,\ldots,C_n$, respectively. The density is then given by
$$
f(X)= \frac{1}{Z_W} e^{-E_W(X)},
$$
where $Z_W$ is the normalizing constant.

For the higher order interactions, only pairs of cells that are neighbours are taken into account. Writing $C_1 \sim C_2$ typically means that the two cells share a common facet ($(d-1)$-face).

In the literature, several definitions of a cell 'being in' the window $W$ are considered. Typical options are that the cell is completely contained in the window ($C\subset W$), the cell intersects the window ($C \cap W \neq \emptyset$) or that the centroid/generator of the cell is in the window. Some of these definitions as well as the higher order interaction require that knowledge about the tessellation outside $W$ is available. 

\begin{figure}[b]
\begin{center}
\includegraphics[width=2.8cm]{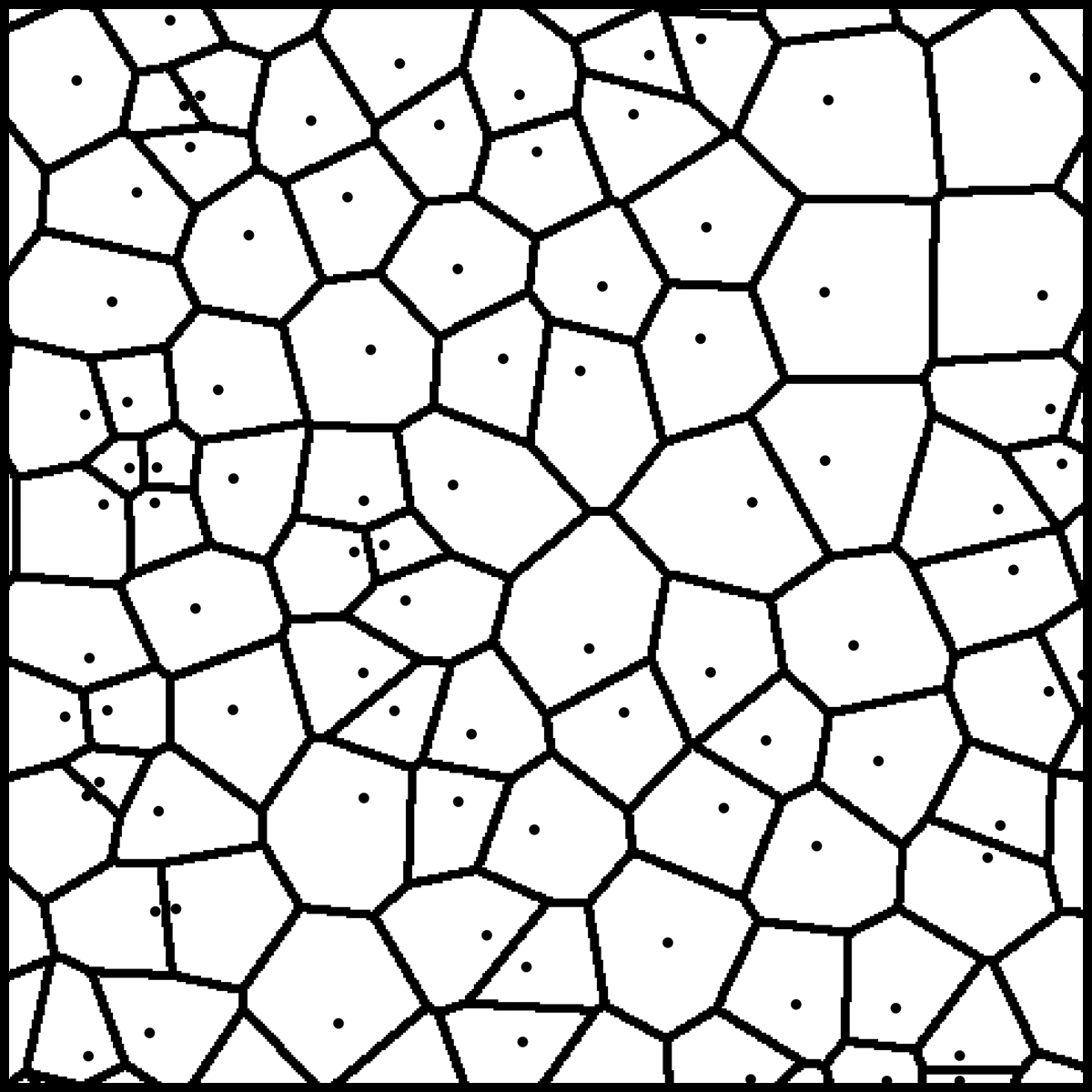} 
\includegraphics[width=2.8cm]{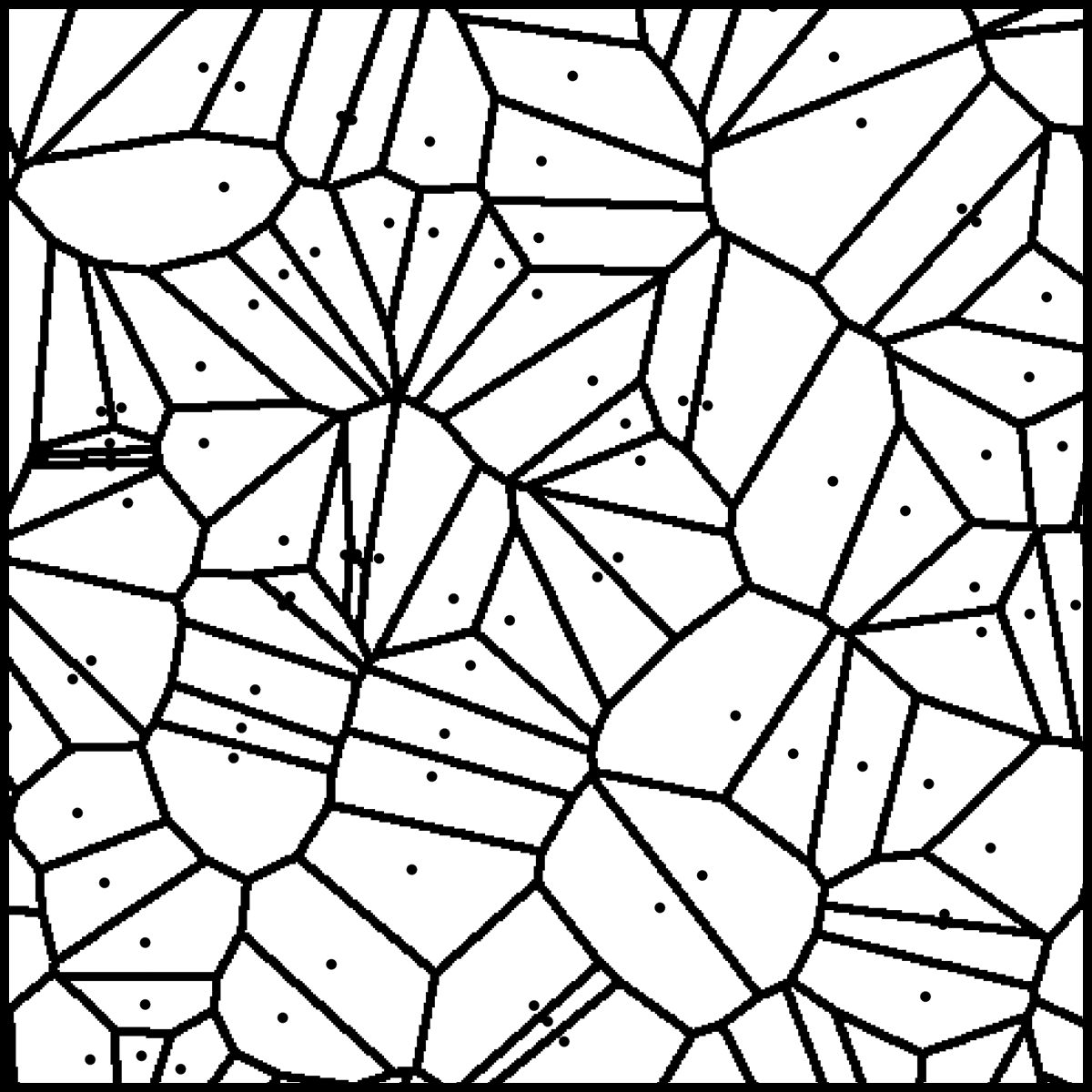}\quad
\includegraphics[width=2.8cm]{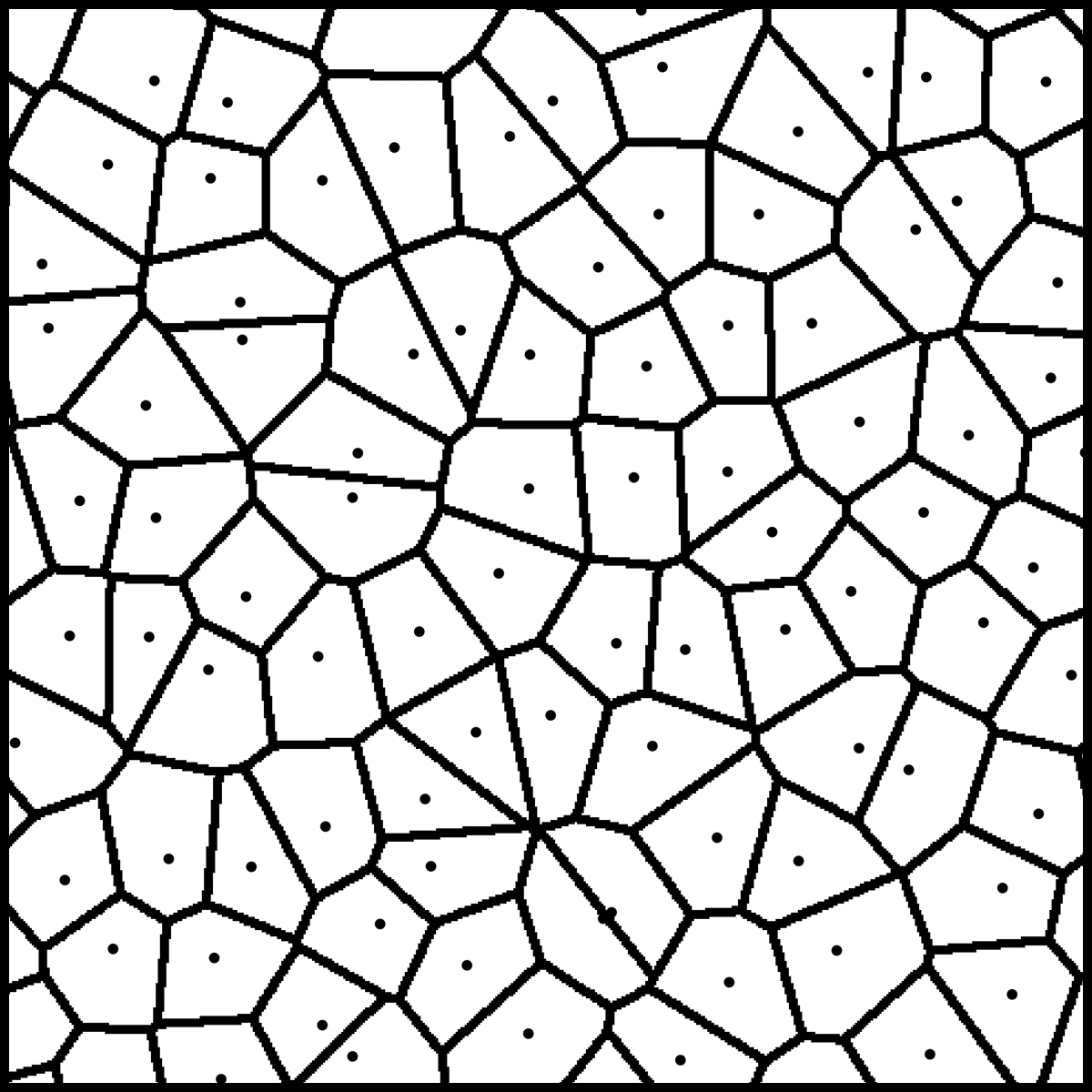} 
\includegraphics[width=2.8cm]{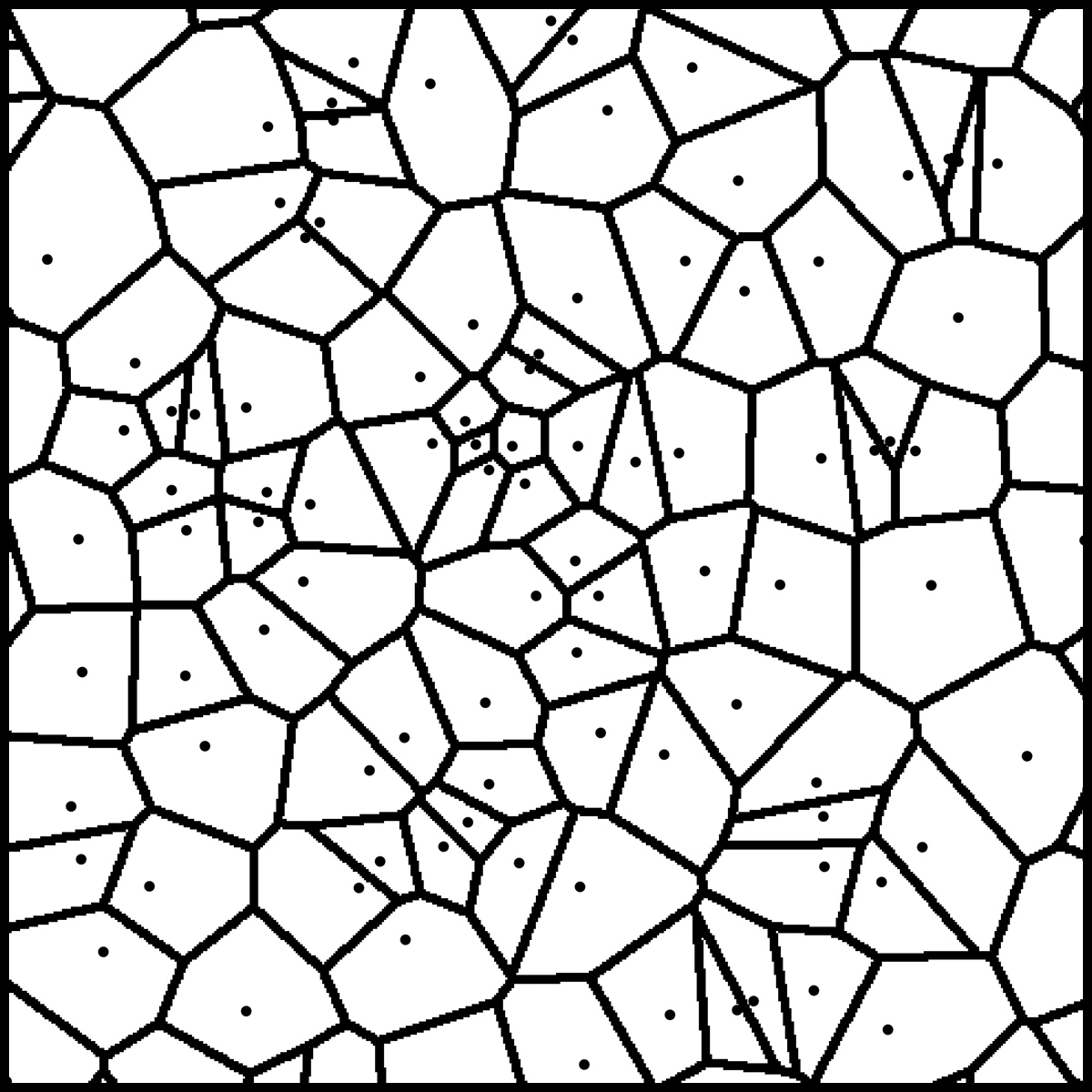} 
\caption{\label{fig:vis_Gibbs} Examples of Gibbs-Voronoi tessellations. Left and left middle: Realizations for ${E_W(X)=\varphi_1(C)=\theta\,\textrm{ar}(C)}$ where $\textrm{ar}(C)\in(0,1]$ denotes the aspect ratio of $C$. Aspect ratios close to 1 are favoured (left) or penalized (left middle). Right middle and right: Realizations for ${E_W(X)=\varphi_n(C_1,\ldots,C_n) = \theta\,\textrm{std}(A(C_1),\ldots,A(C_n))}$ where std$(A(C_1),\ldots,A(C_n))$ denotes the standard deviation of the cells' areas. A small standard deviation is favoured (right middle) or penalized (right).}
\end{center}
\end{figure}

Typical examples for $\phi_1$ are $\phi_1(C)= \theta \size(C)$, where $\size$ can, e.g., be measured by the volume, the surface area,  the total edge length or the mean width of a cell. The sign of the parameter $\theta$ implies that large cells are favoured (negative) or penalized (positive). The value of $\theta$ controls the strength of this effect. Instead of the size, also shape characteristics such as the isoperimetric shape factor or aspect ratios can be used.   
Some examples for such models are shown in Fig.~\ref{fig:vis_Gibbs}.

As an additional option, Dereudre \cite{dereudre_gibbs_2008} studies models with hardcore interactions. In such models, the density is set to zero (or the energy to $+\infty$) if a tessellation fulfils certain geometric properties, e.g., if it contains cells whose size is above or below a prescribed threshold. 

Second order potentials can be used to enforce that neighbouring cells are of similar or dissimilar size or shape.

Theoretical results for Gibbs point processes are often built on the assumption that the density is \emph{hereditary}. That means that any subconfiguration $\mathbf{y}$ of a point pattern $\mathbf{x}$  with nonzero density also has nonzero density. Formally, 
$$
f(\mathbf{x}) > 0 \Rightarrow f(\mathbf{y}) >0 \mbox{ for any } \mathbf{y} \subset \mathbf{x}. 
$$
This assumption is not necessarily met in case of Gibbs models for random tessellations. E.g., assume that a tessellation $X$ does not contain any cell with volume larger than a threshold $t$. Upon removal of a generator point of $X$, the region covered by this point's cell will be assigned to neighbouring cells. As a consequence their volume may then exceed $t$. A theory for nonhereditary densities is developed in \cite{DereudreLavancier2009}. General existence results that can be applied to energies formulated based on Delaunay or Voronoi tessellations are proven in \cite{dereudre_existence_2012}.
Simulation and parameter estimation for the models is discussed in \cite{dereudre_practical_2011}. An approach for model validation based on the notion of residuals for Gibbs processes \cite{BaddeleyResiduals2005} is also presented.

A corresponding theory for Gibbs-Laguerre tessellations has also been developed. Existence of the models is established in \cite{jahn_existence_2020}. Practical issues including simulation by a Birth-Death-Move Metropolis-Hastings algorithm and parameter estimation by pseudo-likelihood are discussed in \cite{seitl_exploration_2021}.


\section{Models for T-tessellations}
\label{Sec:Temporal}
In this section we will summarize several approaches for the construction of tessellations that are not face-to-face. The models are characterized by their T-shaped vertices but differ in the distributions of sizes and shapes of their cells. In most cases, the constructions are restricted to $\R^2$.

Various constructions for T-tessellation models have been introduced. Some of them turn out to generalize or be special cases of others. One of the main characteristics of such models is the direction distribution of the edges. In the literature, mostly two choices are considered: isotropic models, where directions are drawn from a uniform distribution on the unit circle, and rectangular tessellations, where only directions parallel to the coordinate axes are considered. Often, models have been formulated for one particular choice but can be straightforwardly generalized to the other case or even arbitrary direction distributions (as long as they are not concentrated on a single direction).
Below, some prominent models are listed in historic order. 

\subsection{The Gilbert tessellation}
The construction of the Gilbert tessellation  \cite{Gilbert} starts by sampling a realisation of a stationary Poisson process in the plane. Each point is independently assigned a direction from the uniform distribution on the circle. Then, cracks start to grow from the Poisson points at uniform speed in both the chosen direction and its opposite. When a crack hits another one, the growth in this direction stops (but may continue in the other one). This construction gives rise to a non-face-to-face tessellation with T-vertices.

\begin{figure}[t]
\begin{center}
\includegraphics[width=3.7cm]{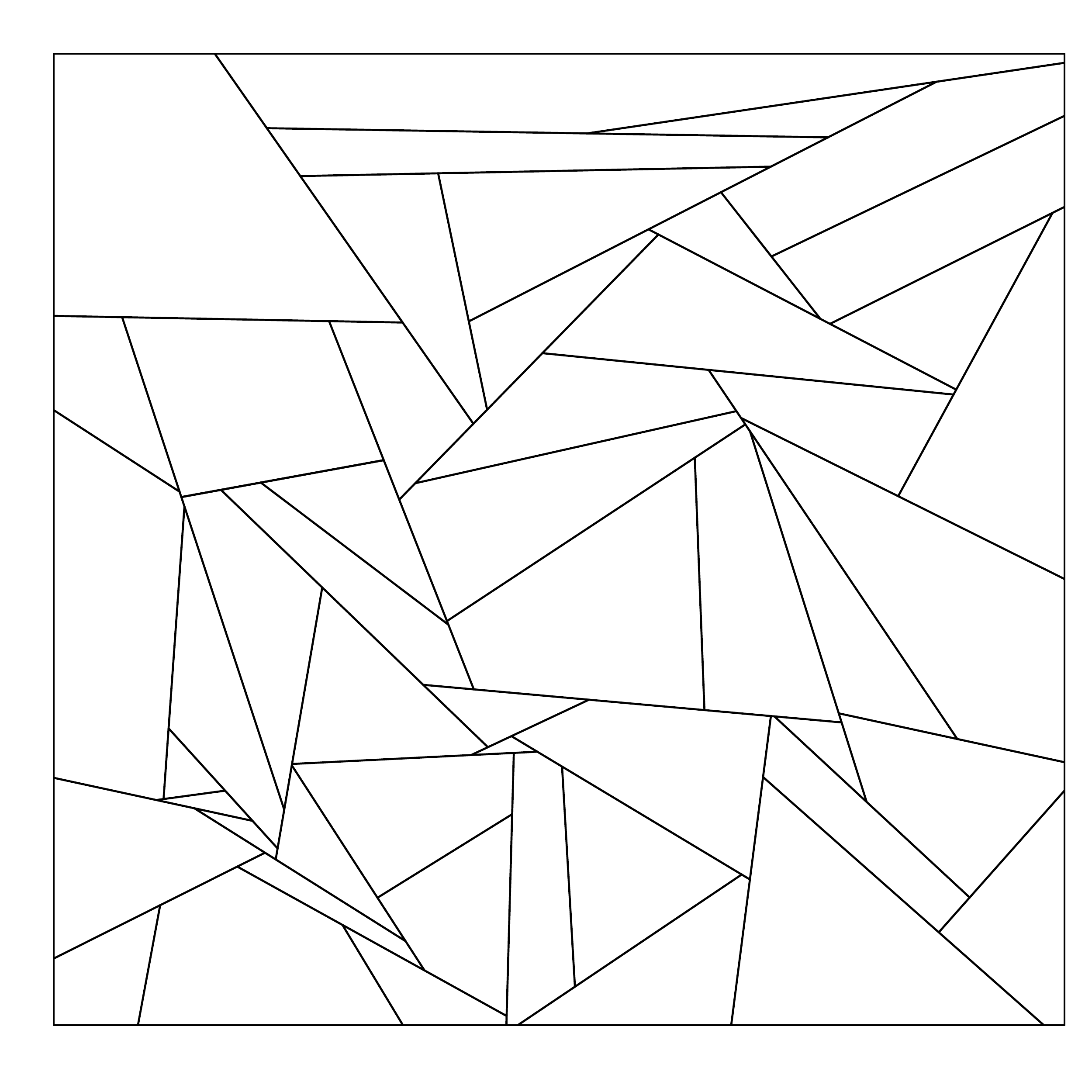} \hspace{0.1cm} \includegraphics[width=3.7cm]{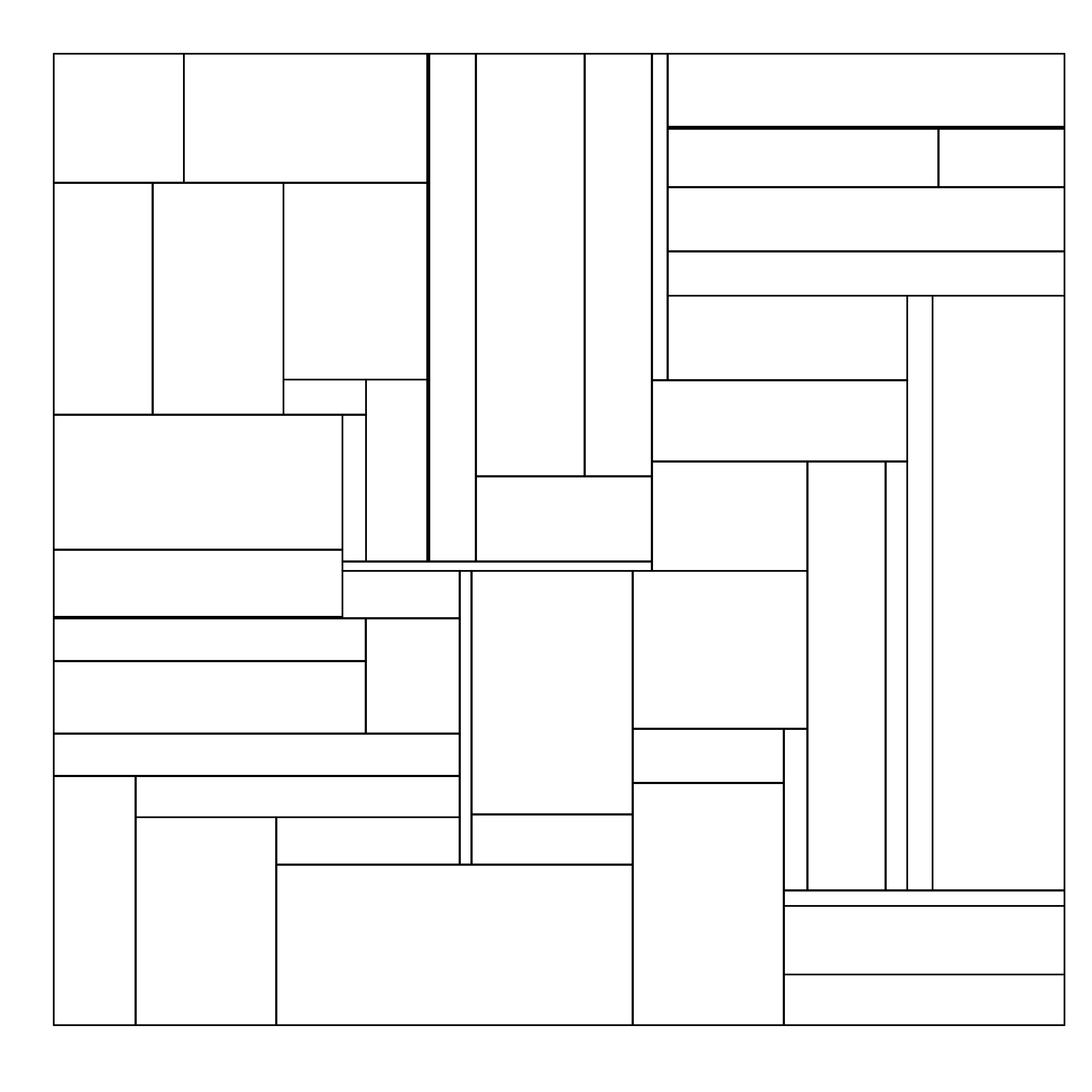} \hspace{0.1cm} \includegraphics[trim={2cm 2cm 1.5cm 1.5cm},clip,width=3.7cm]{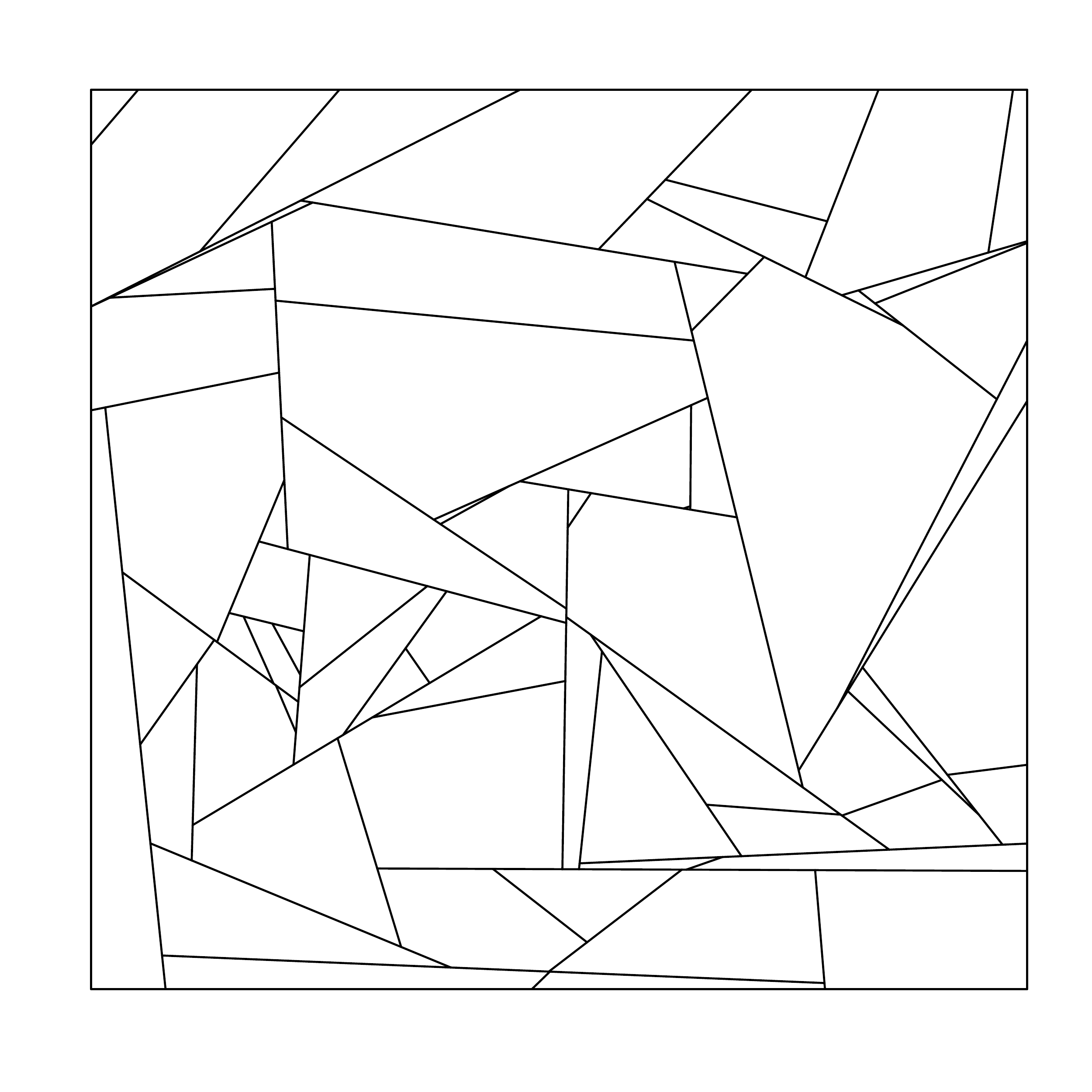}
\caption{\label{fig:Gilbert} Left: Gilbert tessellation. Middle: Gilbert tessellation with rectangular cells. Right: Arak-Clifford-Surgailis tessellation.}
\end{center}
\end{figure}

A variant of this model with rectangular cells is obtained when restricting the growth directions to the coordinate axis directions. This model is studied by Mackisack and Miles in \cite{mackisack_miles_1996} where the authors come to the conclusion that it is not tractable to significant theoretical analysis.
Simulations of the isotropic and the rectangular Gilbert model are shown in Fig.~\ref{fig:Gilbert}.

Burridge et al. \cite{burridge_full-_2013} further adapt the rectangular Gilbert model to a half-Gilbert tessellation. It is obtained when only pairs of east and south or north and west growing segments block each other in the growth process. They also provide a series expansion for the edge length distribution.

 An iterated ($k$th order) version of the Gilbert model is introduced in \cite{BACCELLI2022752}. In this model, any particle will survive $k$ collisions with other particles before it dies.

\subsection{The Arak-Clifford-Surgailis tessellation}
\label{sec:Arak}
This tessellation model in $\R^2$ is constructed as a Markov random graph \cite{arak_point-based_1993,arak_markov_1989}. The tessellation is generated by the trajectories of particles moving in a compact and convex observation window $W$.

As in Definition \ref{def:hyperplane}, let $\cL$ be the space of lines in $\R^2$. For any Borel set $B$ let $[B]$ denote the elements of $\cL$ hitting $B$. Let $\Lambda$ be a nonnegative Borel measure on $\cL$ with $\Lambda([B])< \infty $ for any bounded Borel set $B$ and $\Lambda([\{x\}])=0$ for any $x\in \R^2$. Denote by $\cH$ a Poisson line process with intensity measure $\Lambda$ and by $\cH_W$ its restriction to $W$. 
We will interpret lines as particle trajectories by letting the $x$-coordinate of any point on the line represent time $t:=x$ and the $y$-coordinate $y=y(t)$ the position at time t. This way, particles move from left to right. Any particle moving along a line $L$ in $[W]$ will enter the window at an entry point $in(L;W)$. 
Let $\Gamma$ be the intensity measure of the point process of entry points $in(L;W)$ of the lines of $\cH$ (as a measure on $\partial W$).

The particles generating the tessellation move independently according to the following rules, see \cite{thale_arak-clifford-surgailis_2011} for more details. A visualization is shown in Fig.~\ref{fig:Gilbert}.

 \begin{programcode}{Algorithm Arak-Clifford-Surgailis tessellation}    
\begin{enumerate}
\item 
Particles are born on the boundary $\partial W$ of the window $W$ according to a Poisson point process with intensity measure $\Gamma$.

Additionally, each point $x$ is assigned a mark giving the initial direction of the particle's movement. The marks follow the same directional distribution as the lines $L$ of $\cH_W$ conditioned on $x$ being the entry point of $L$ into $W$.
\item 
The particles move freely and with constant velocity until one of the following events occurs.
\begin{itemize}
\item 
A particle that hits the boundary of the window dies.
\item 
If two particles collide, a fair coin is tossed to determine which of the two particles survives the collision. The surviving particle continues in its initial direction, the other one dies.
\item 
Particle trajectories may also branch. Consider the line $L(e)$ containing a particle trajectory $e$. Branching points are simulated by a Poisson process with intensity measure $B \mapsto\frac{1}{2} \Lambda([B])$, where $B$ is a Borel set on the line $L(e)$. As soon as the particle reaches a branching point, a second particle is born and starts its movement in a random direction which is chosen according to the direction distribution of $\cH$. The original particle continues moving along $L(e)$.
\end{itemize}
\end{enumerate}
 \end{programcode}

In 1995 and 1997, Miles and Mackisack presented a T-tessellation model with rectangular cells on two conferences. It turned out, that their construction can be generalized to arbitrary orientation distribution of the edges, and corresponds to the model by Arak, Clifford, and Surgailis \cite{MilesMackisack2002}. Some analytical results for this model are presented in \cite{MilesMackisack2002}. The most prominent finding was that the distribution of the typical cell of the ACS model corresponds to that of a Poisson line tessellation. 

In \cite{thale_arak-clifford-surgailis_2011}, explicit formulas for second order characteristics are given for the isotropic case. 

Note that we consider only a special case of the polygonal Markov field construction described in \cite{arak_point-based_1993,arak_markov_1989}. In the more general version of the model, particles may also be born inside the window and may change direction during their movement. This way, tessellations with nonconvex cells and V-shaped vertices are obtained. Additionally, one may consider the case that both particles survive a collision resulting in X-shaped vertices.  

\subsection{The STIT tessellation}
\label{sec:STITTess}

STIT tessellations were introduced by Nagel and Weiss in \cite{Nagel:2003:LSS} motivated by non-face-to-face tessellations as they appear in crack patterns on ceramic surfaces or in dry soil. Initially, the tessellation was defined by an iterative cell division process on a bounded window $W$. To define this process, we set up some notation.

 Let $\Rc$ denote an even probability distribution on the unit sphere $S^{d-1}$. We define a translation invariant measure $\Lambda$ on $\cH$ via
$$
\int_{\cH} f(H) \Lambda (dH) = \int_{S^{d-1}} \int_{\R} f(H(u, r)) dr \, \Rc(du)  
$$
for any nonnegative measurable function $f: \cH \to [0, \infty)$. 
In this setting $\Rc$ determines the distribution of normal directions of the hyperplanes. We assume that $\Rc$ is not concentrated on a great circle.

For a given polytope $P$, we define a probability distribution 
\begin{equation}
   \Lambda_{[P]} (\cdot)= \frac{\Lambda(\cdot \cap[P])}{ \Lambda([P])} 
\end{equation}
with $[P]$ as in Section \ref{sec:Arak}.

\begin{definition}[STIT tessellation]
\label{def:STIT}
A random \textit{STIT tessellation} in any bounded Borel set $W\subset{\R}^d$ with $0<\Lambda([W])<\infty$ is obtained by the following spatio-temporal construction. 

 \begin{programcode}{Algorithm STIT}
Let $(\tau_i, \ell_i), i \ge 1$ be an independent and identically distributed (i.i.d.) sequence of pairs where $\tau_i$ follows an exponential distribution with parameter $\Lambda([W])$ and $\ell_i$ is a random hyperplane with distribution $\Lambda_{[W]}$.
Choose a stopping time $a>0$.
Start at time $t=0$ with an empty window $W$.
\begin{itemize}
\item
If $\tau_1>a$, return $W$ as tessellation.
\item
If $\tau_1\le a$, compute the cells $W_+$ and $W_-$ obtained by intersecting $W$ with $\ell_1$.
\item
Iterate the procedure for $W_+$ and $W_-$ independently. We describe the procedure for $W_+$. \\Update $t$ to $t+\tau_2$. 
\begin{itemize}
    \item 
    If $t>a$, stop the evolution of $W_+$ and keep it as cell of the tessellation.
    \item 
    If $t\le a$, check if $\ell_2$ intersects $W_+$. If yes, split $W_+$ by $\ell_2$ into $W_{++}$ and $W_{+-}$ and continue the process on both cells separately. If no, the process is continued on $W_+$.
\end{itemize}
\end{itemize}
\end{programcode}
\end{definition}

\begin{figure}[t]

\begin{center}
\includegraphics[width=3.7cm]{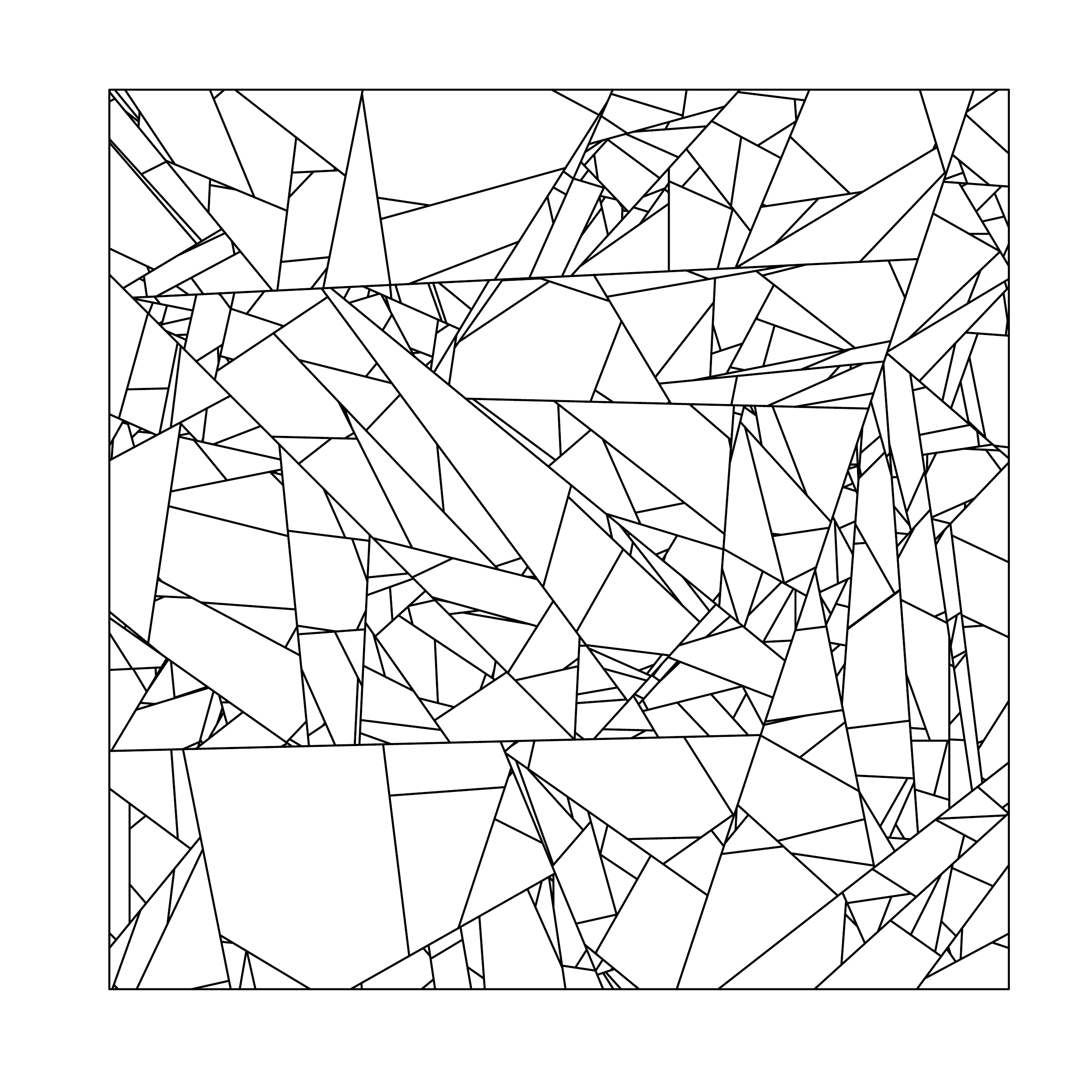} \hspace{0.1cm} \includegraphics[width=3.7cm]{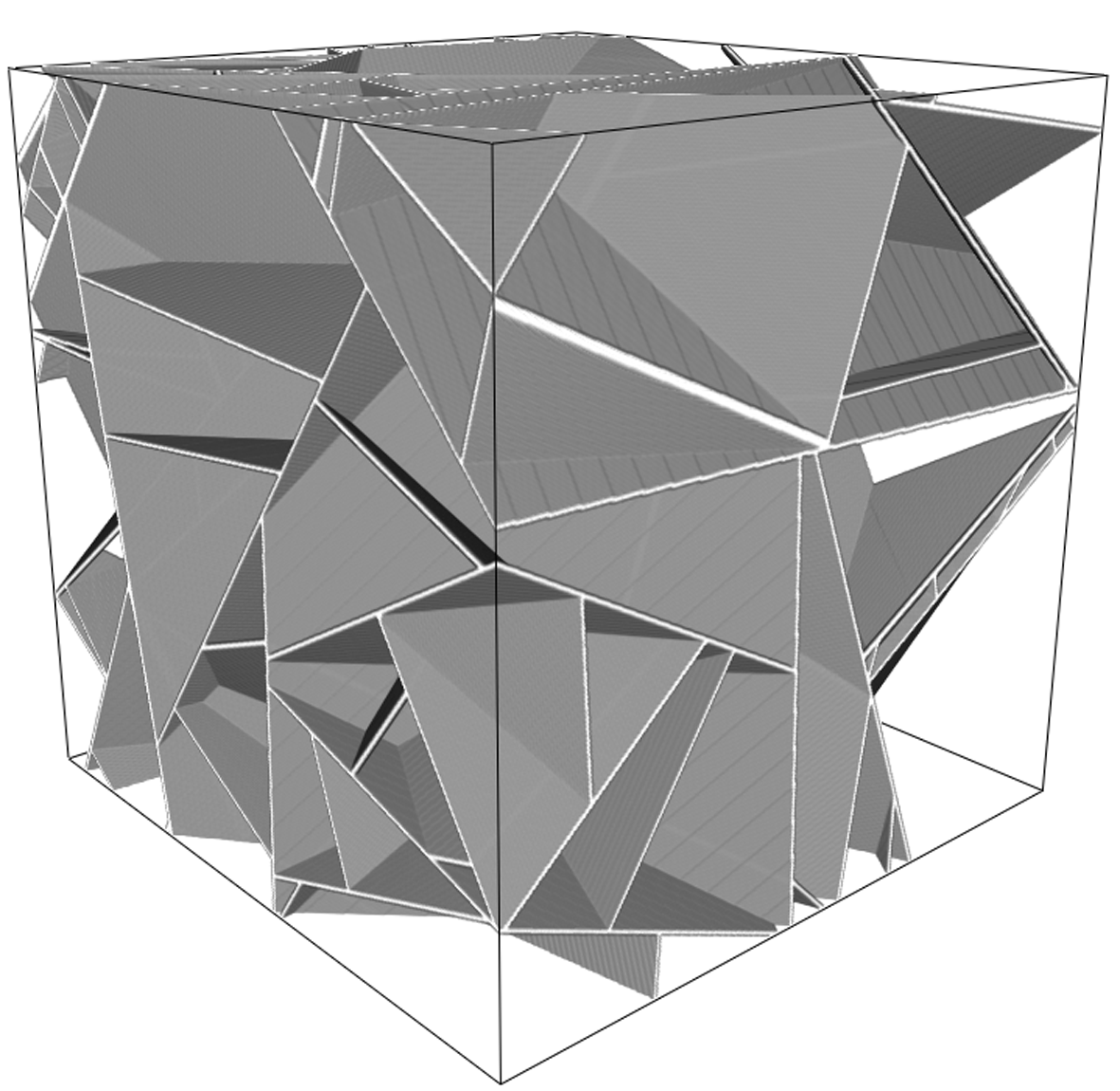} \hspace{0.1cm} \includegraphics[width=3.7cm]{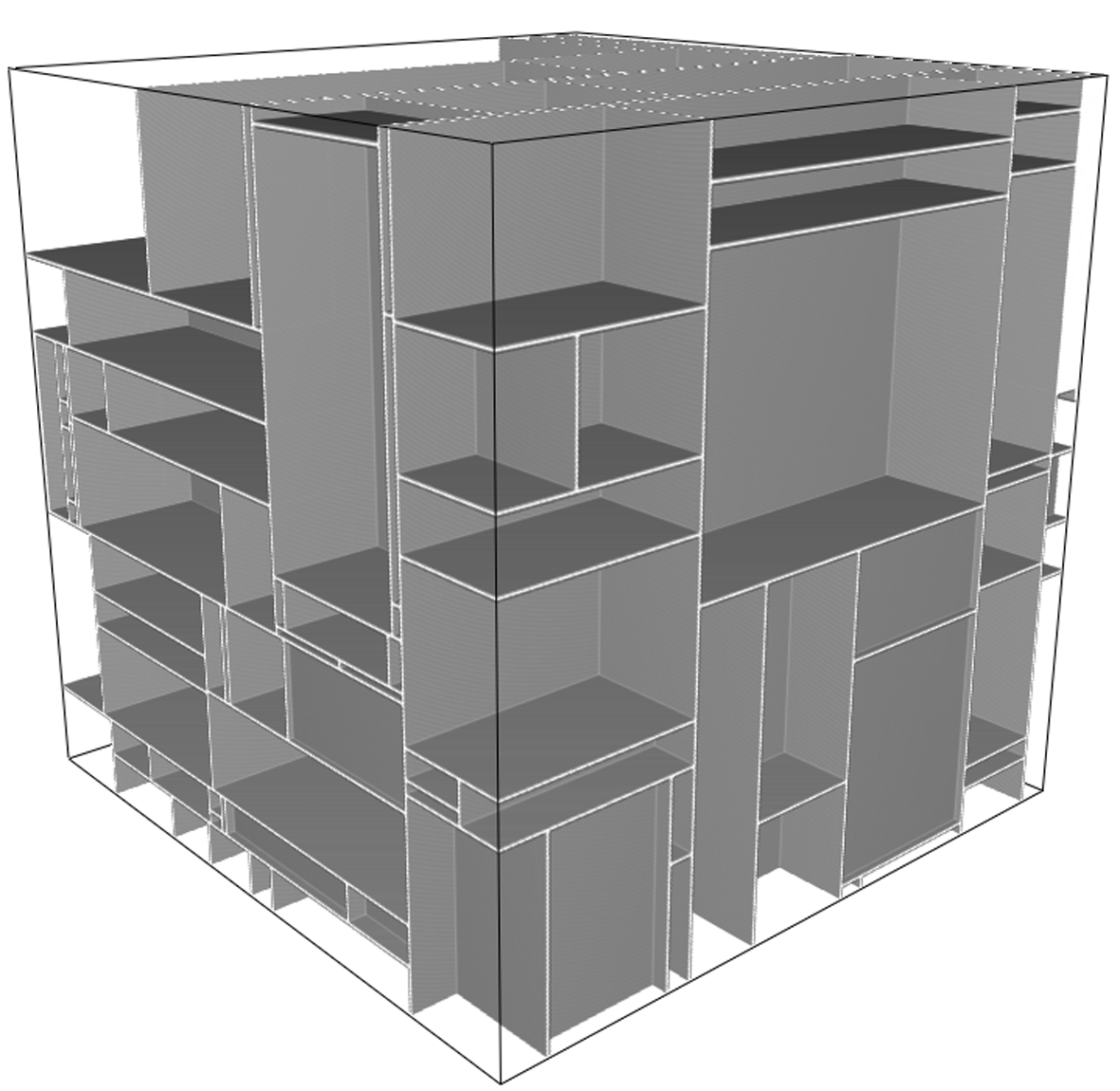}\\
\caption{\label{fig:STIT}Isotropic STIT tessellations in $\R^2$ and $\R^3$ and a STIT tessellation in $\R^3$ where $\Rc$ is concentrated on the coordinate directions (from left to right).}
\end{center}
\end{figure}

A detailed description of the implementation is given in \cite{leon_parameter_2023}. Visualisations of two- and three-dimensional STIT tessellations are shown in Fig.~\ref{fig:STIT}.

Model parameters are the stopping time $a$ (which equals $\mu_{d-1}$) and the distribution $\Rc$ of normal directions of the hyperplanes.

The name STIT stems from the fact that these tessellations are STable with respect to ITeration \cite{nagel_crack_2005}. 

\begin{definition}[Iteration of tessellations]
Let $T_0=\{C_1, C_2, \ldots\}$ be a stationary tessellation and let $\mathcal{T} =(T_1, T_2, \ldots)$ denote a sequence of i.i.d. tessellations that are independent of $T_0$. The \emph{iteration (or nesting)} of $T_0$ and $\mathcal{T}$ is a random tessellation given by
\begin{equation}
\label{eq:iteration}
   I(T_0, \mathcal{T}) = \{C_i \cap C_{ij}  \, : \, i,j \in \N,  C_i \cap C_{ij} \neq \emptyset \},
\end{equation}
the system of (non-empty) cells obtained by intersecting the cells $C_{ij}$ of $T_i$ with cell $C_i$ of $T_0$.
\end{definition}

That is, the part of tessellation $T_i$ that is visible in cell $C_i$ of the initial tessellation $T_0$ is copied into that cell. 

Repeated iteration $I_m$ by i.i.d. sequences $\mathcal{T}_1, \mathcal{T}_2, \ldots$ is defined via
\begin{align}
\label{eq:repeatediteration}
I_2(T_0)&= I(2T_0, 2 \mathcal{T}_1)\\
I_m(T_0)&=I( \frac{m}{m-1} I_{m-1}(T_0), m \mathcal{T}_{m-1}),
\end{align}
where $mT_i$ refers to the tessellation $T_i$ upscaled by the factor $m$.

Iteration according to \eqref{eq:iteration} will result in 'finer' tessellations with increased cell intensity $\gamma_d$ and expected total facet content $\mu_{d-1}$. 
The rescaling in \eqref{eq:repeatediteration} is therefore required to preserve the scale of the tessellations such that $\mu_{d-1}$ is preserved after each iteration step. 

In the following, we will assume that $T_0$  and all tessellations in the sequences $\mathcal{T}_1, \mathcal{T}_2, \ldots$ are i.i.d. 
\begin{definition}[Stable with respect to iteration (STIT)]
A stationary tessellation $T$ is \emph{stable with respect to iteration} if $T$ and $I_m(T)$ are equally distributed for all $m=2,3 \ldots.$ 
\end{definition}

The tessellation constructed as in Definition \ref{def:STIT}, originally referred to as \emph{Crack STIT tessellation}, is stable with respect to iteration which explains the name STIT. The original construction from \cite{Nagel:2003:LSS} is restricted to a bounded window. In \cite{nagel_crack_2005}, the existence of a stationary random tessellation that is stable with respect to iteration and that is locally obtained by the Crack STIT construction was shown. The authors also showed that, for given parameters $a$ and $\mathcal{R}$, the Crack STIT tessellation is the unique stationary tessellation that is STIT. Hence, the model was henceforth simply called \emph{STIT tessellation} in the literature. A global construction of the model is given in \cite{MNW08}. 
STIT tessellations with direction distribution $\mathcal{R}$ concentrated on the coordinate axis directions are also known as \emph{Mondrian tessellations} \cite{betken_second-order_2023} as they visually resemble paintings of Piet Mondrian.

Several first-order properties of this model in the planar and three-dimensional case were derived in \cite{NW06}, \cite{NW08} and \cite{TW}. Second-order theory for STIT tessellations was developed in \cite{ST}. 

It was shown that the interiors of the typical cells of a Poisson hyperplane tessellation and a STIT tessellation with the same parameters have the same distribution \cite{nagel_crack_2005,SCHREIBER2011989}. A comparison  of Poisson-Voronoi, Poisson hyperplane and STIT tessellations with equal total facet content $\mu_{d-1}$ is presented in \cite{RedTha12a}. The authors come to the conclusion that regarding cell sizes, STIT tessellations are closer to the Poisson hyperplane tessellation (due to the equality of the typical cell). In contrast, mean face numbers and topological characteristics are closer to those of a Poisson-Voronoi tessellation. In terms of second order properties, STIT tessellations are positioned in between the other two models. The arrangement of cells in these models has been discussed in \cite{RedTha12}.

As mentioned above, the Arak-Clifford-Surgailis model is another tessellation model sharing the distribution of the typical cell with the Poisson line tessellation. When comparing the pair correlation function of the ACS and STIT tessellations, it turns out that ACS tessellations have a shorter range of dependence than STIT tessellations, see \cite{thale_arak-clifford-surgailis_2011}. 

\subsubsection{STIT extensions: iterative cell division}
\label{Sec:STITExtension}
Attempts to model real patterns by STIT tessellations were not very successful \cite{nagel_tessellation_2008}. This might be explained by differences between physics-driven cell splitting processes and the STIT construction. For instance, in reality large cells might be more likely to split than smaller ones. Additionally, the newly appearing fracture might not come from a uniform distribution but rather split a cell close to its center. Based on these considerations, Cowan \cite{Cowan2010} extended the interative cell division process by introducing several variants of cell selection and cell division rules resulting in the following algorithm. 

\begin{programcode}{Algorithm Iterative cell division} \label{alg:Cowan}
Start at time $t=0$ with a bounded convex window $W_0$ of area 1 centered in the origin\\
   Iterate
		\begin{itemize}
  \item $t \mapsto t+1$
  \item
			Randomly select a cell in $W_{t-1}$ by a given \emph{cell selection rule}
   \item
			Divide $W_{t-1}$ by a random chord drawn according to a given \emph{cell division rule} 
   \item Scale the window by a factor of $\sqrt{t+1}/\sqrt{t}$ to get $W_t$
		\end{itemize}
\end{programcode}

The expansion of the window is chosen such that a mean cell area of 1 is preserved throughout the process. In the planar case, Cowan suggests the cell selection rules \emph{equally likely} (each cell with the same probability), \emph{area, perimeter, and corner weighted} (select cell $i$ with  weight $w_i/\sum_{k=1}^t w_k$ where $w_i$  is the area, perimeter, or number of vertices, respectively).

Cell division rules considered are
\begin{itemize}
\item
choose a chord through a uniform point in the cell,
\item
choose a uniform chord in the space of lines hitting the cell (uniform random),
\item
choose a line from a uniform point on the boundary to an 'opposite' corner
\end{itemize}
The isotropic STIT tessellation is the special case of this model when choosing the perimeter-weighted selection rule and the uniform division rule.

\begin{figure}
\begin{center}
\includegraphics[width=3.7cm]{STIT_ISO.pdf}  
\includegraphics[width=3.7cm]{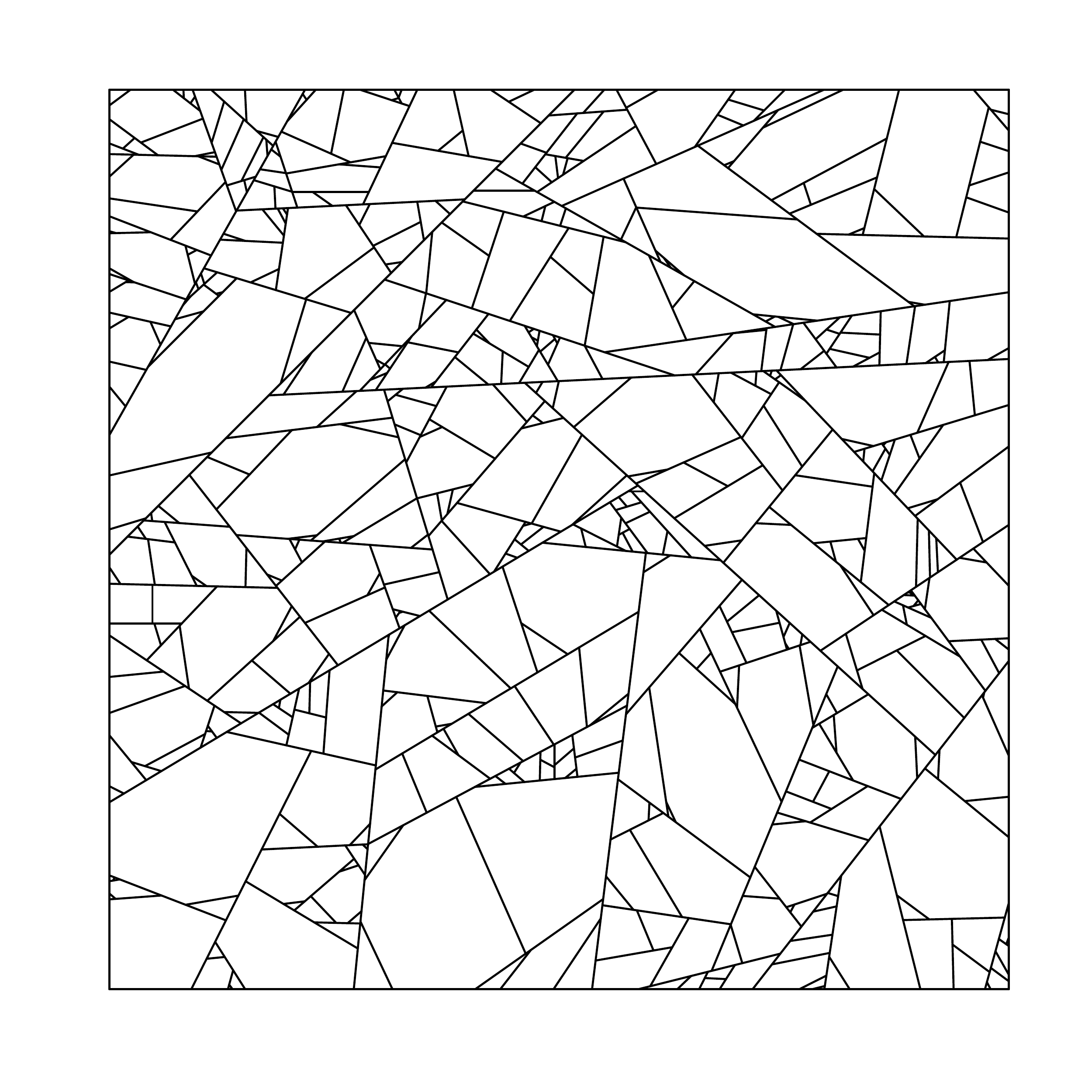} 
\includegraphics[width=3.7cm]{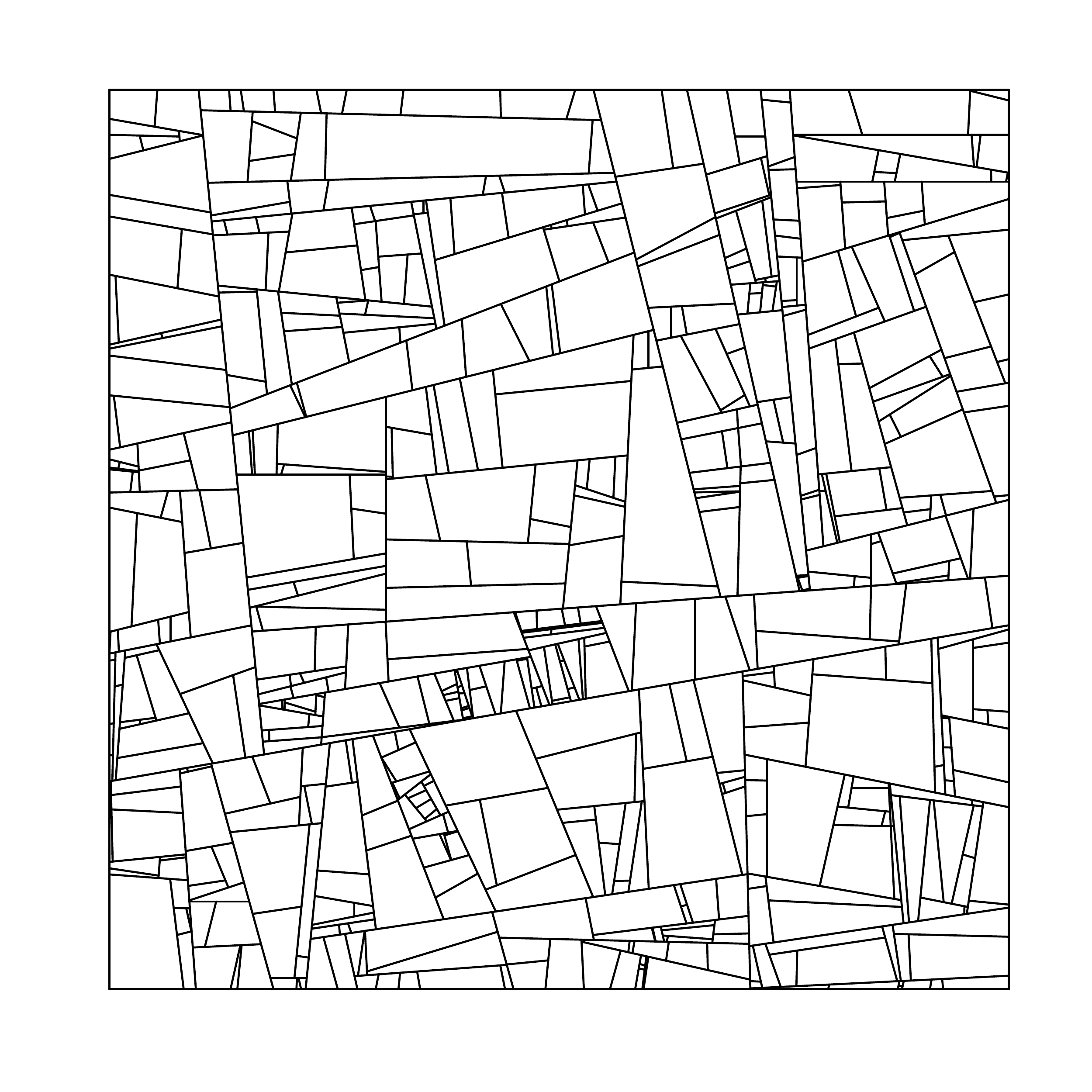} \\
\includegraphics[width=3.7cm]{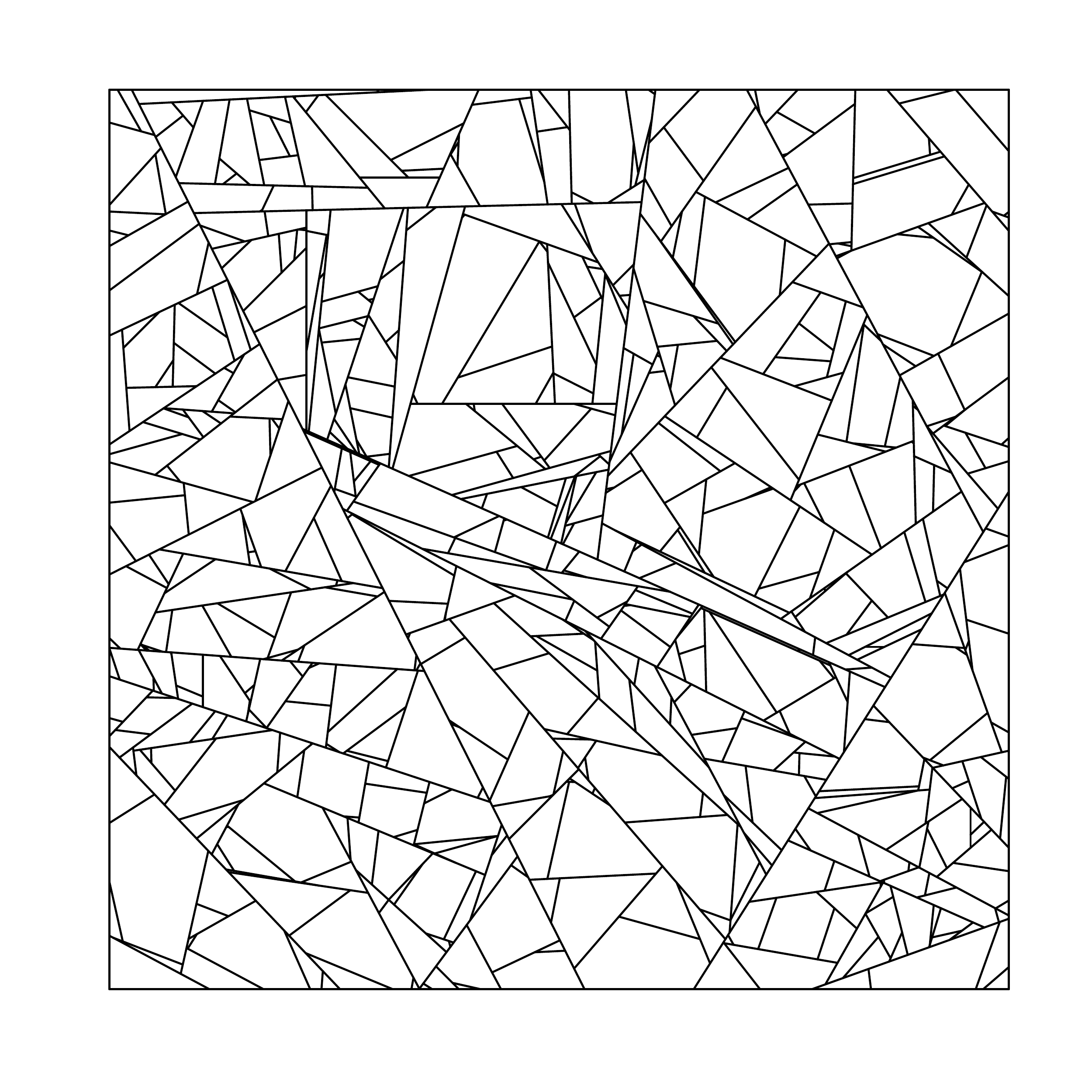} 
\includegraphics[width=3.7cm]{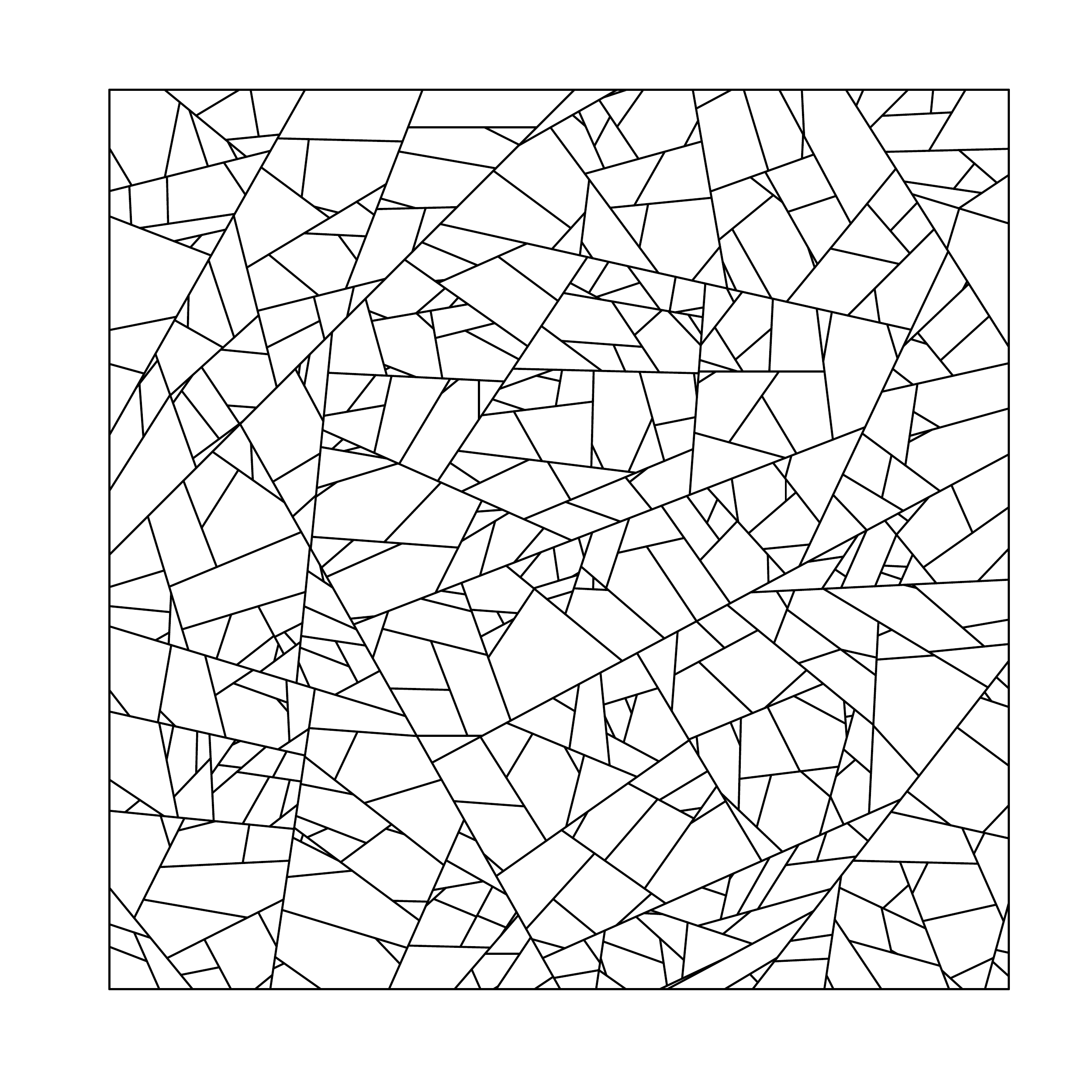} \includegraphics[width=3.7cm]{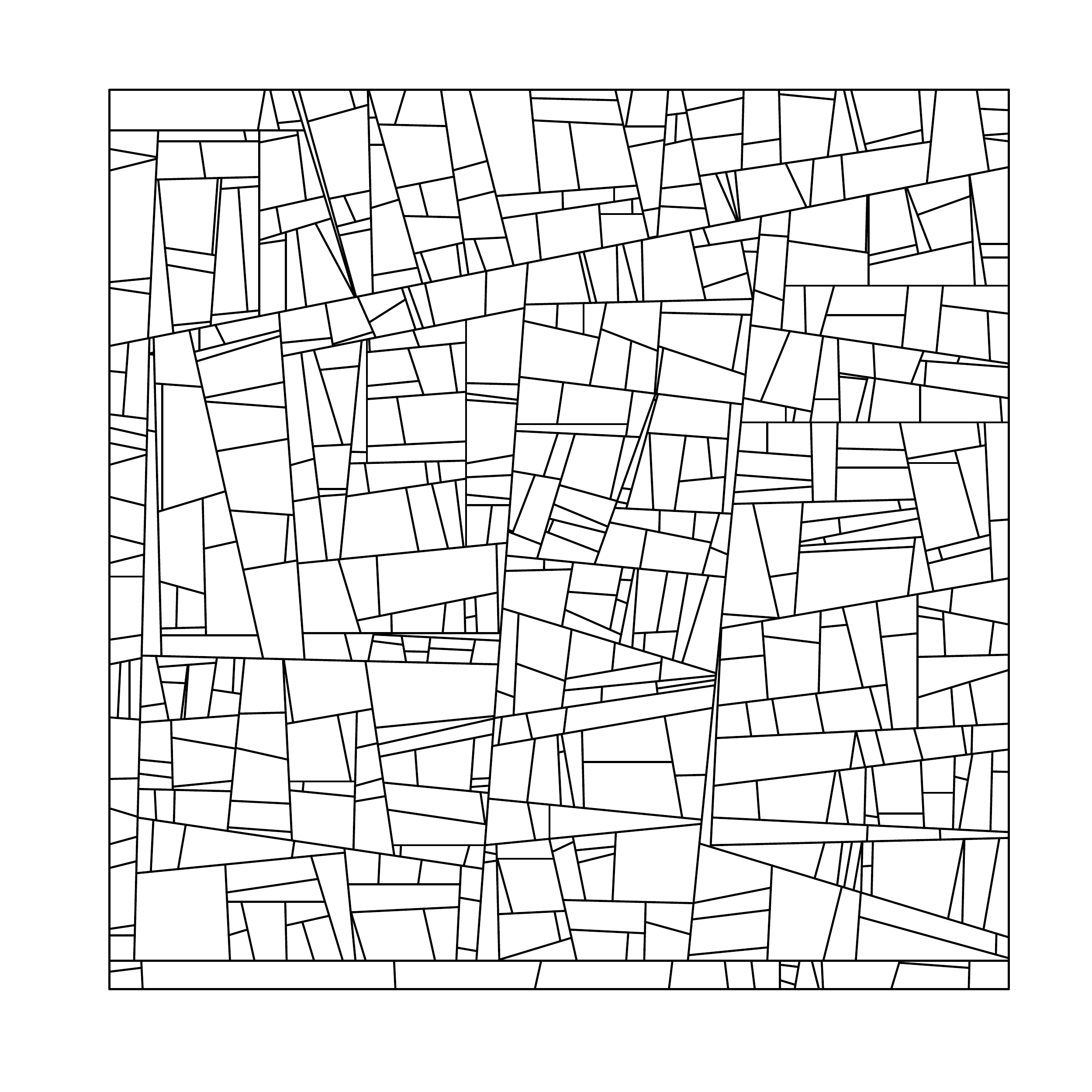}
\caption{\label{fig:VariantsSTIT} Some variants of the STIT model as introduced in \cite{leon_parameter_2023}. Division rule D-STIT in all cases. Top row: lifetime L-STIT, bottom row: lifetime L-AREA resulting in more homogeneous cell area distributions. Left: No modification, middle: RDSSQ resulting in more circular cells, right: ASA avoiding small angles. Parameters were chosen such that all tessellations contain around 500 cells.}
\end{center}
\end{figure}

Leon and Nagel \cite{leon_parameter_2023} extend this idea by listing a framework consisting of 40 model variations in $\R^2$. They consider several models for the directional distribution of the line segments (uniform, discrete, disturbed discrete, and elliptical). The division rules are D-STIT (as in the STIT model) and D-GAUSS (more central splitting). Besides the STIT lifetime L-STIT  which is exponentially distributed with parameter $\Lambda([C])$, they consider an area weighted lifetime L-AREA where the parameter of the exponential distribution is $A(C)$, the area of the cell. Additionally, the authors suggest two modifications that are supposed to generate rounder cells. In D-RDMIN, for given direction of the dividing line, the location of division is chosen such that the minimal isoperimetric shape factor $\min (RD(C_1), RD(C_2))$ of the two resulting cells $C_1$ and $C_2$ is maximized. The isoperimetric shape factor is defined as  $RD(C)= 4 \pi A(C)/P^2(C)$ and equals 1 for a circle. In D-RDSSQ, the  sum of squares $RD(C_1)^2+RD(C_2)^2$ is maximized. Finally, the option ASA (avoid small angles) is introduced which rejects divisions resulting in edge angles smaller than a given threshold $\omega$. Some versions of these models are shown in Fig.~\ref{fig:VariantsSTIT}.
For a Gibbsian extension of STIT tesselltions, see \cite{Georgii}.

\subsubsection{Completely random and Gibbsian T-tessellations}

\label{sec:TTess}
A general theory for so-called T-tessellation, i.e., tessellations with only T-shaped vertices, is developed in \cite{kieu_completely_2013}. The authors only consider tessellations on a bounded, convex, and polygonal window $W$. Similar to the literature for STIT-tessellations, they distinguish between edges (corresponding to the K-segments) and segments (corresponding to the I-segments).

\begin{figure}
\begin{center}
\includegraphics[trim={0 0 3cm 12cm},clip,width=2.7cm]{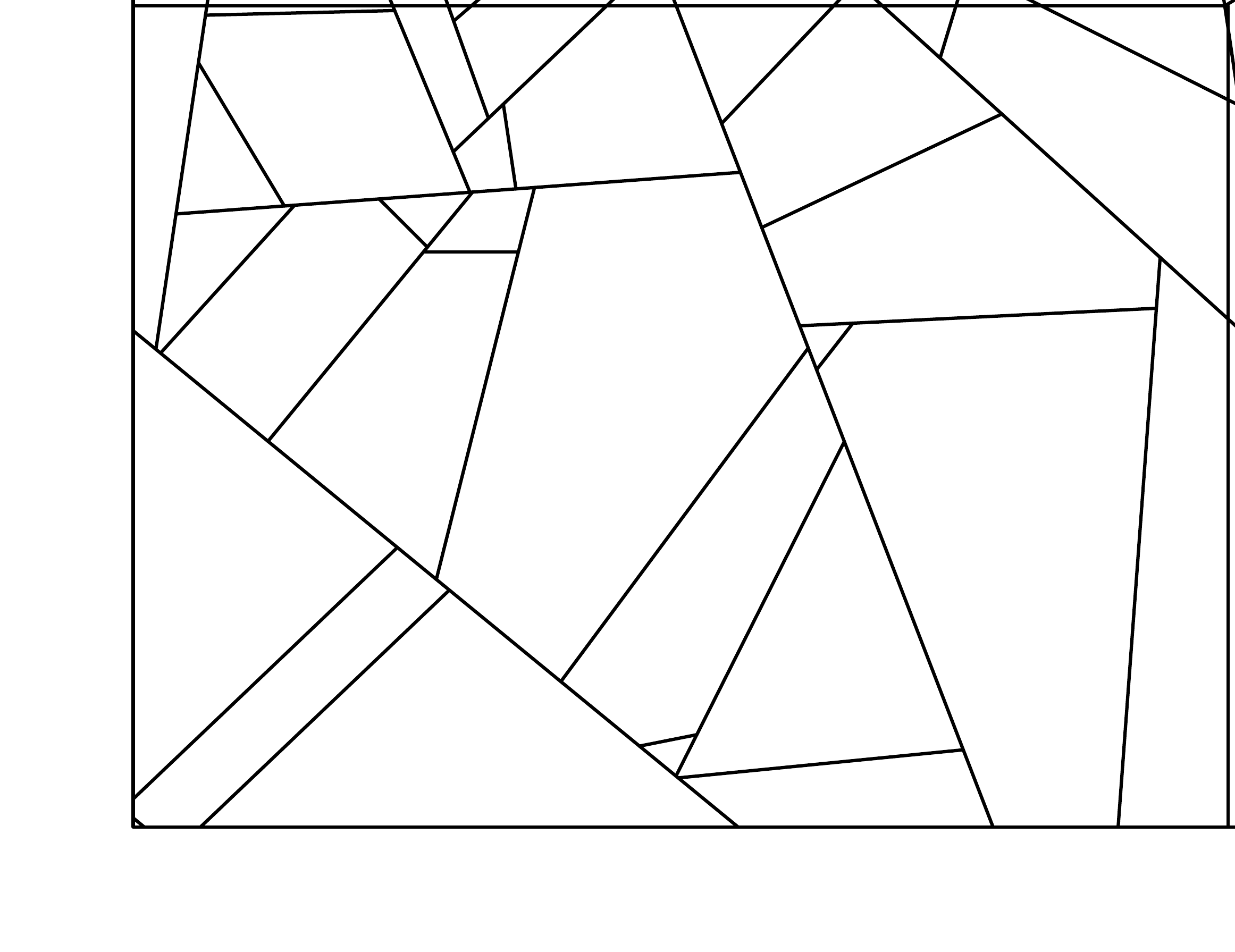}
\includegraphics[trim={0 0 3cm 12cm},clip,width=2.7cm]{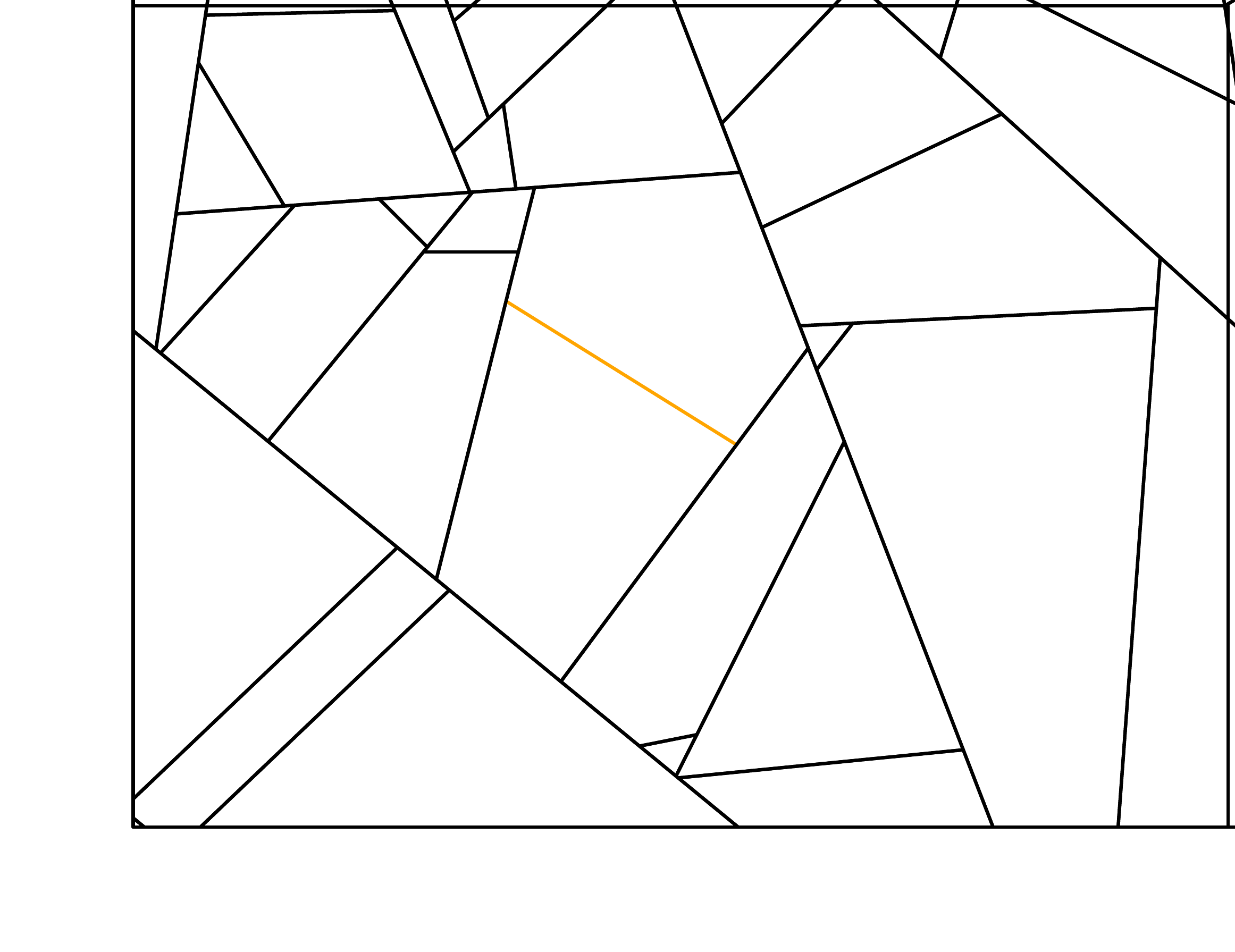}
\includegraphics[trim={0 0 3cm 12cm},clip,width=2.7cm]{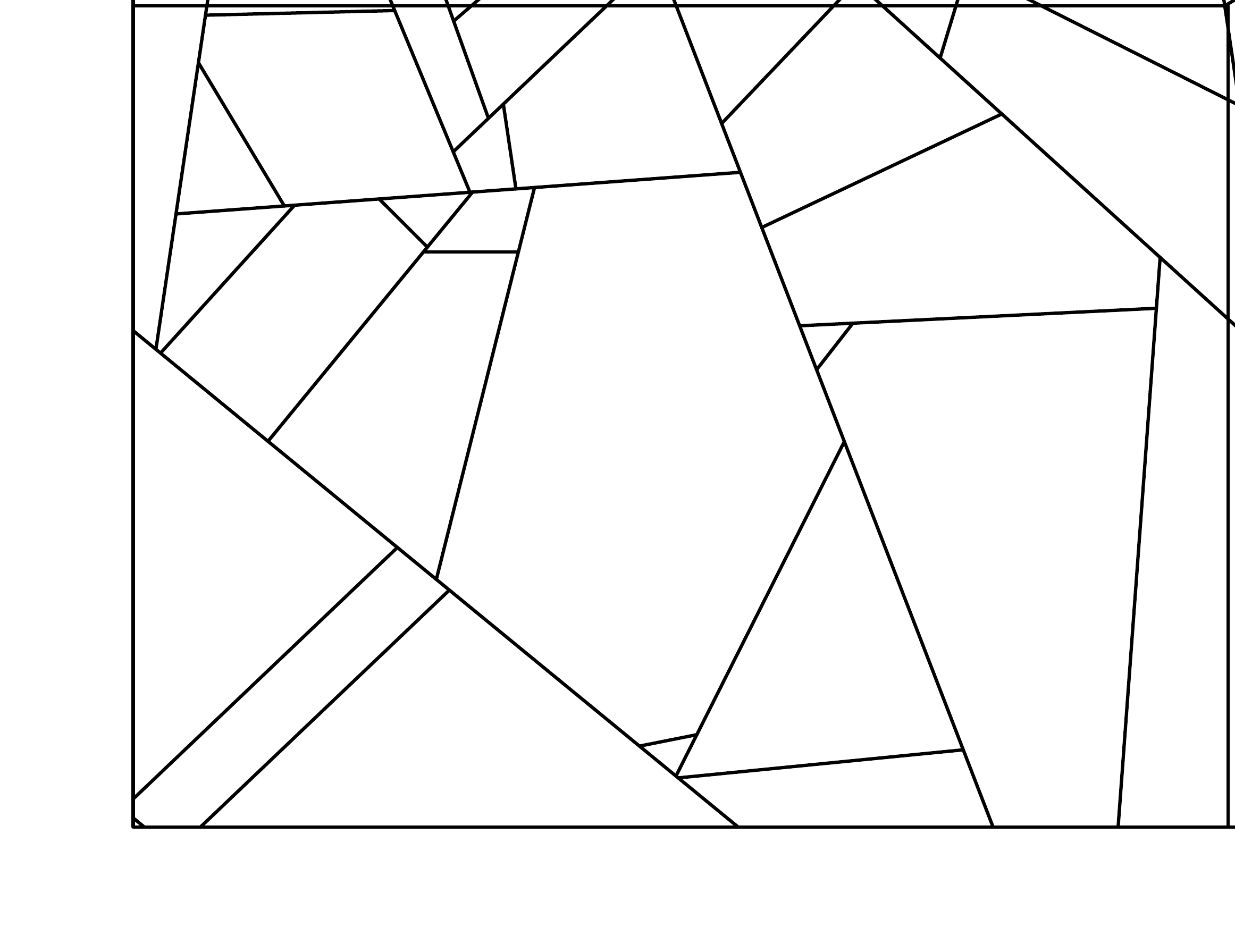}
\includegraphics[trim={0 0 3cm 12cm},clip,width=2.7cm]{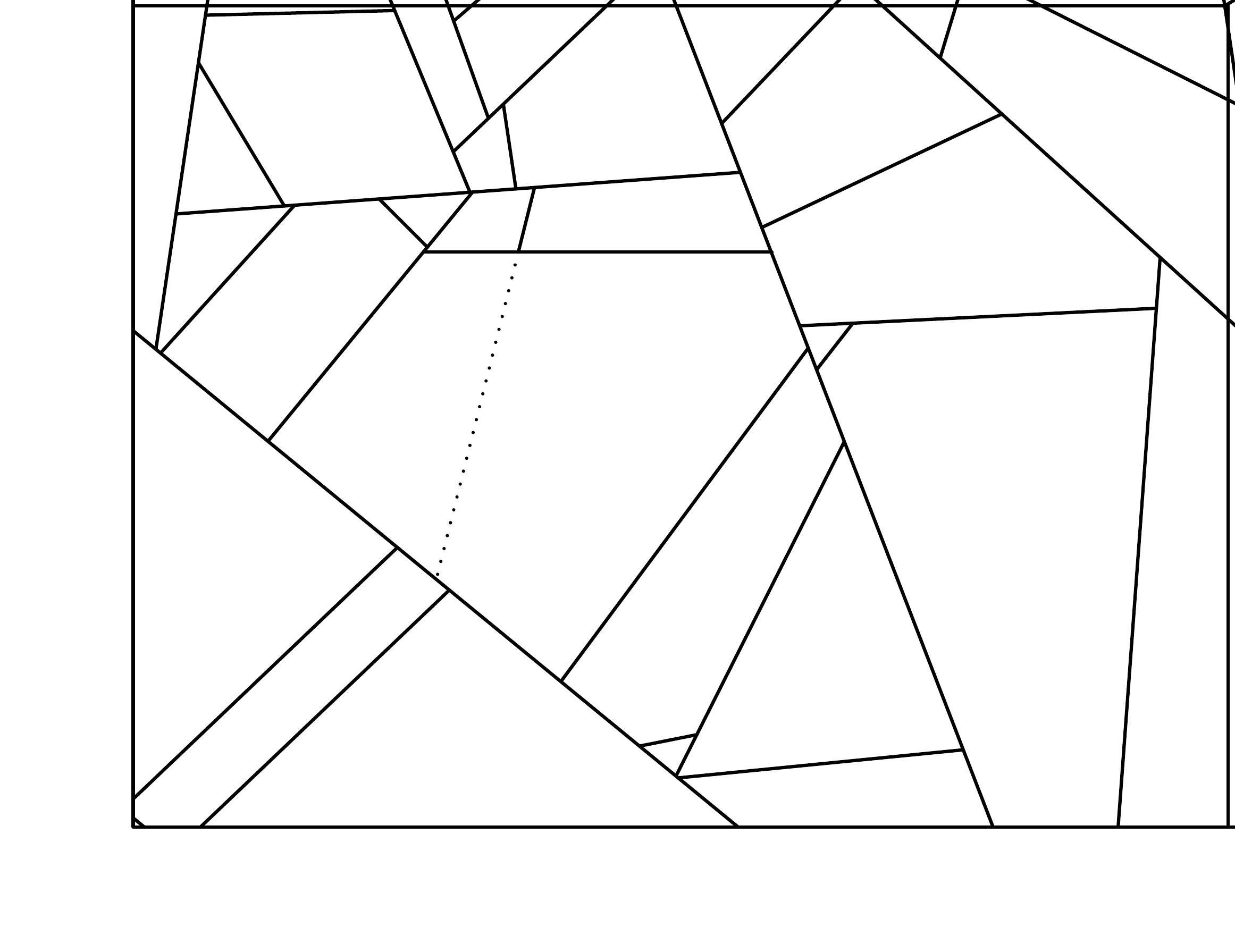}\\
\end{center}
\caption{Operations on T-tessellation: starting configuration, split of a cell by the orange line, merge of the center cells by deletion of an edge, and flip of the dotted edge (from left to right).}\label{Fig:operationsT}
\end{figure}

They show that every T-tessellation can be generated from the empty window $W$ by a finite sequence of \emph{splits}, \emph{merges}, and \emph{flips}.
A split divides one of the cells by a new line segment. A merge deletes a segment consisting of a single edge to merge the two neighbouring cells. Such segments are called nonblocking, while all others are blocking segments. A flip deletes an edge at the end of a blocking segment while simultaneously adding a new edge. This edge is chosen such that it extends a segment that was blocked by the deleted edge. An illustration of these operations is given in Fig.~\ref{Fig:operationsT}.

For point processes, the Poisson process serves as a model for complete spatial randomness. Similarly, a completely random T-tessellation model can be defined.

\begin{definition}[Completely random T-tessellation (CRTT)]
Let $\Tb_T$ be the space of T-tessellations on $W$ equipped with the $\sigma$-algebra $\sigma(\Tb_T)$ defined in analogy to Equation \eqref{Eq:sigmaAlgebra}. For a T-tessellation $T$ let $L(T)$ denote the minimal line system in $W$ containing all edges of $T$. For given $L$, we denote by $\Tb_T(L)$ the set of all T-tessellations $T$ such that $L(T)=L$.

Let $L$ be the restriction of a Poisson line process with unit intensity $\lambda=1$ to the window $W$. A completely random T-tessellation is distributed according to the probability measure $\mu$ given by
\begin{equation}
\label{EQ:DefCRTT}
\mu(A)= Z^{-1} \E \sum_{T \in \Tb_T(L)} \one_A(T), \quad A \in \sigma(\Tb_T).
\end{equation}
\end{definition}
That is, conditioned on the realization of $L$, a completely random T-tessellation follows a uniform distribution on the set $\Tb_T(L)$. For Equation \eqref{EQ:DefCRTT} to be well-defined, it is shown that $\Tb_T(L)$ is finite. The normalizing constant $Z$ is not known.

In \cite{kieu_completely_2013}, results similar to the Slivnyak theorem for the Poisson point process are shown for the model with distribution \eqref{EQ:DefCRTT}. This justifies calling this model a completely random T-tessellation. 

To increase the flexibility of the model, Gibbsian extensions of the completely random T-tessellation are considered in \cite{adamczyk-chauvat_gibbsian_2020} and \cite{kieu_completely_2013}.  

\begin{definition}[Gibbsian T-tessellation]
For a stable, nonnegative functional $h$ on $\Tb_T$, the Gibbs random T-tessellation with unnormalized density $h$ follows a distribution
$$
P(dT) \;\propto \;  h(T) \mu(dT).
$$
\end{definition}

The constant of proportionality has to be chosen such that $P$ has a total mass of $1$. Therefore, the unnormalized density must have a nonzero and finite integral over $\Tb_T$ with respect to $\mu$. As in the case of Gibbs point processes, this can be ensured by a stability condition on $h$, see \cite{kieu_completely_2013} for details. An algorithm for simulating Gibbsian T-tessellations that is based on the Metropolis-Hastings-Green algorithm is also given in \cite{kieu_completely_2013}.

Typical choices of densities are formulated with the aim of controlling the number of cells as well as the distributions of their sizes and shapes. To consider a few examples, we let $\mathring{n}_s(T)$ and $\mathring{n}_v(T)$ denote the number of segments and vertices, respectively, of the tessellation which do not lie on the boundary of $W$.  Additionally, let $l(T)$ denote the total edge length of $T$.

Initially, the completely random T-tessellation was defined with respect to a unit intensity process. A scaled version of the model is obtained when using a Gibbsian density based on
$$
- \log h(T)= -\mathring{n}_s(T) \log (\tau),
$$ where $\tau$ is a scaling parameter that allows to control the number of cells.

The Arak-Clifford-Surgailis tessellation is a Gibbsian T-tessellation with the choice
$$
- \log h(T)= \frac{\tau}{\pi} l(T) + \mathring{n}_v(T) \log 2 - \mathring{n}_s(T) \log (\tau).
$$

Kieu et al. \cite{kieu_completely_2013} conclude that the ACS model tends to produce more small cells than the CRTT model and confirm this by simulation.

\begin{figure}[t]
\begin{center}
\includegraphics[trim={2cm 2cm 1.5cm 1.5cm},clip,width=3.7cm]{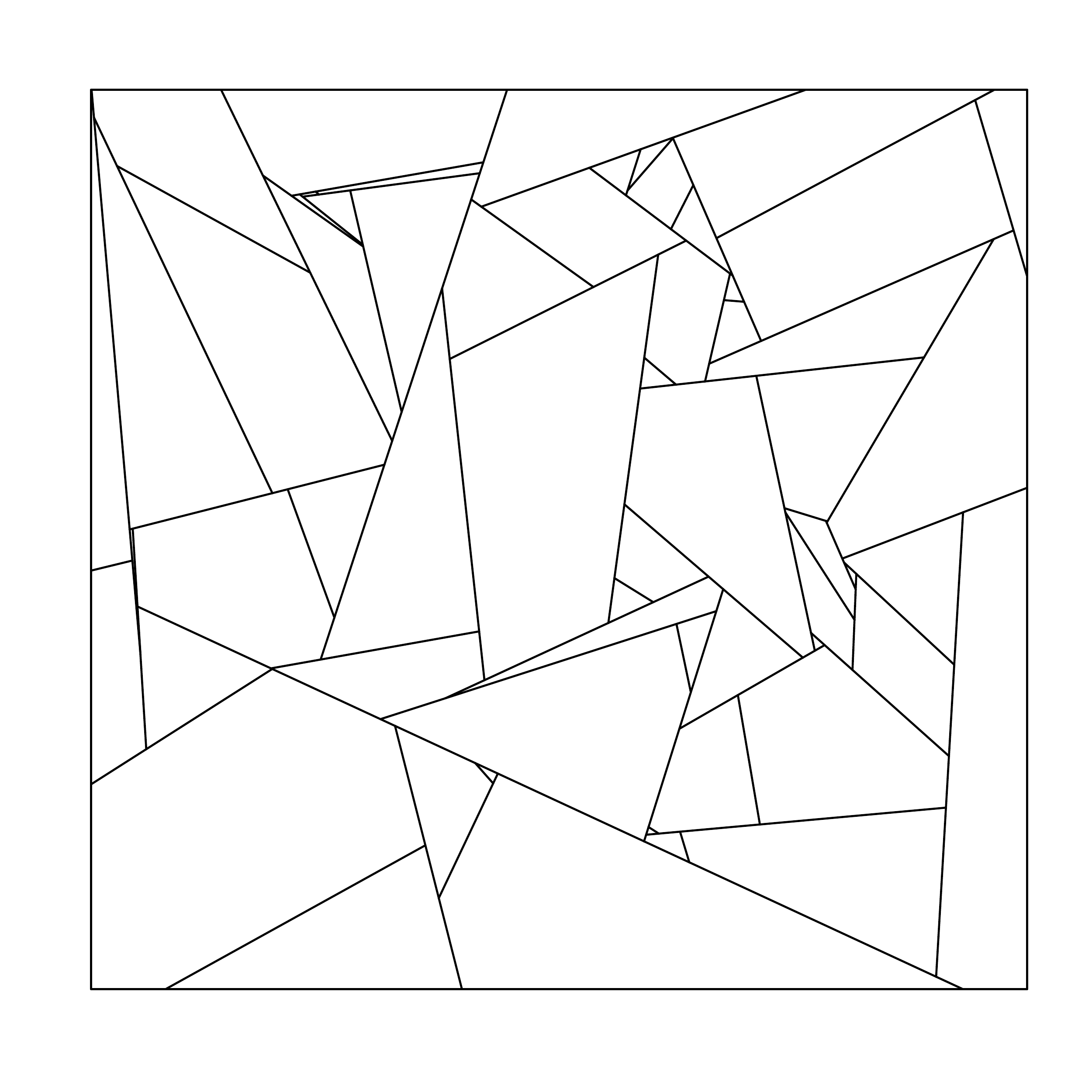} 
\includegraphics[trim={2cm 2cm 1.5cm 1.5cm},clip,width=3.7cm]{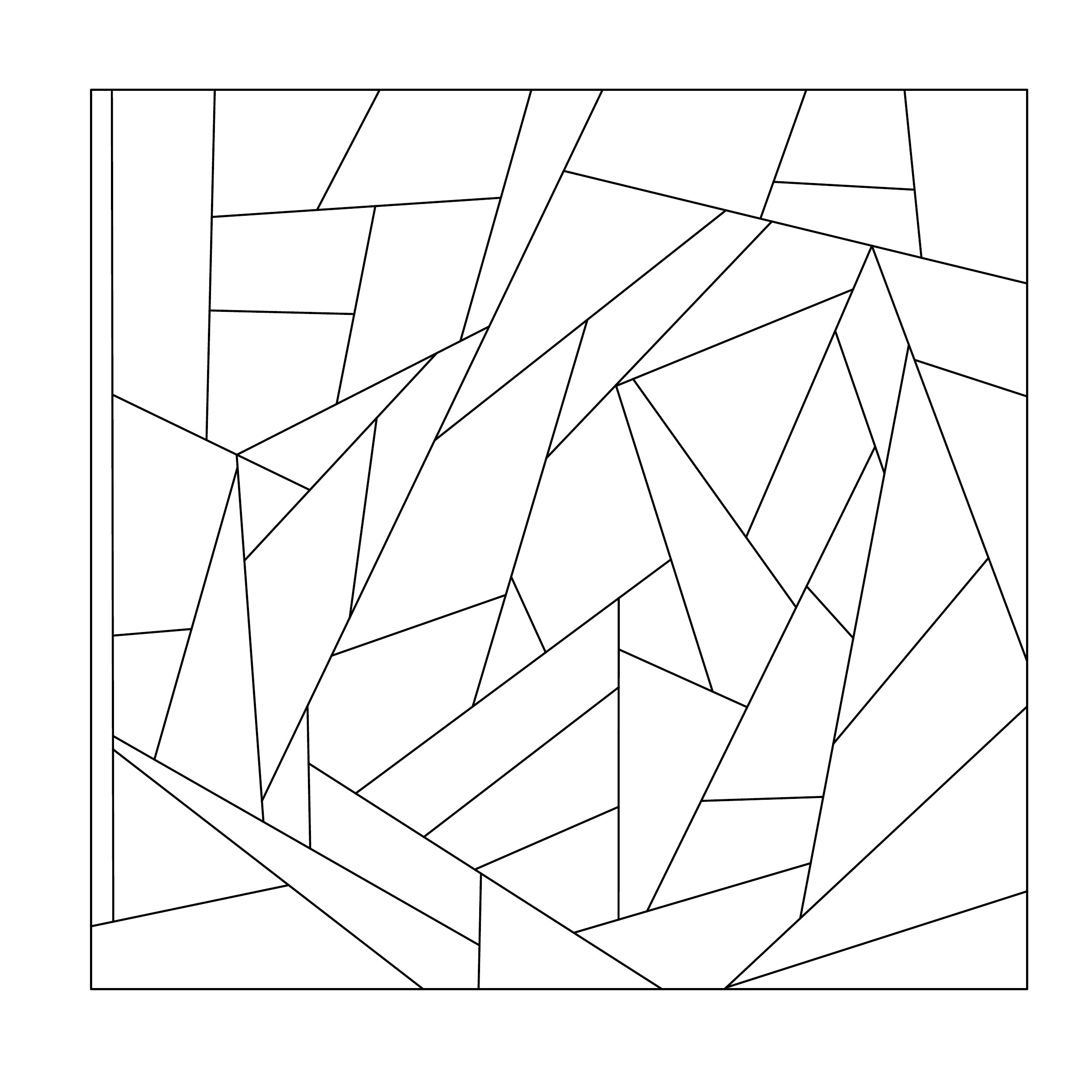} 
\includegraphics[trim={2cm 2cm 1.5cm 1.5cm},clip,width=3.7cm]{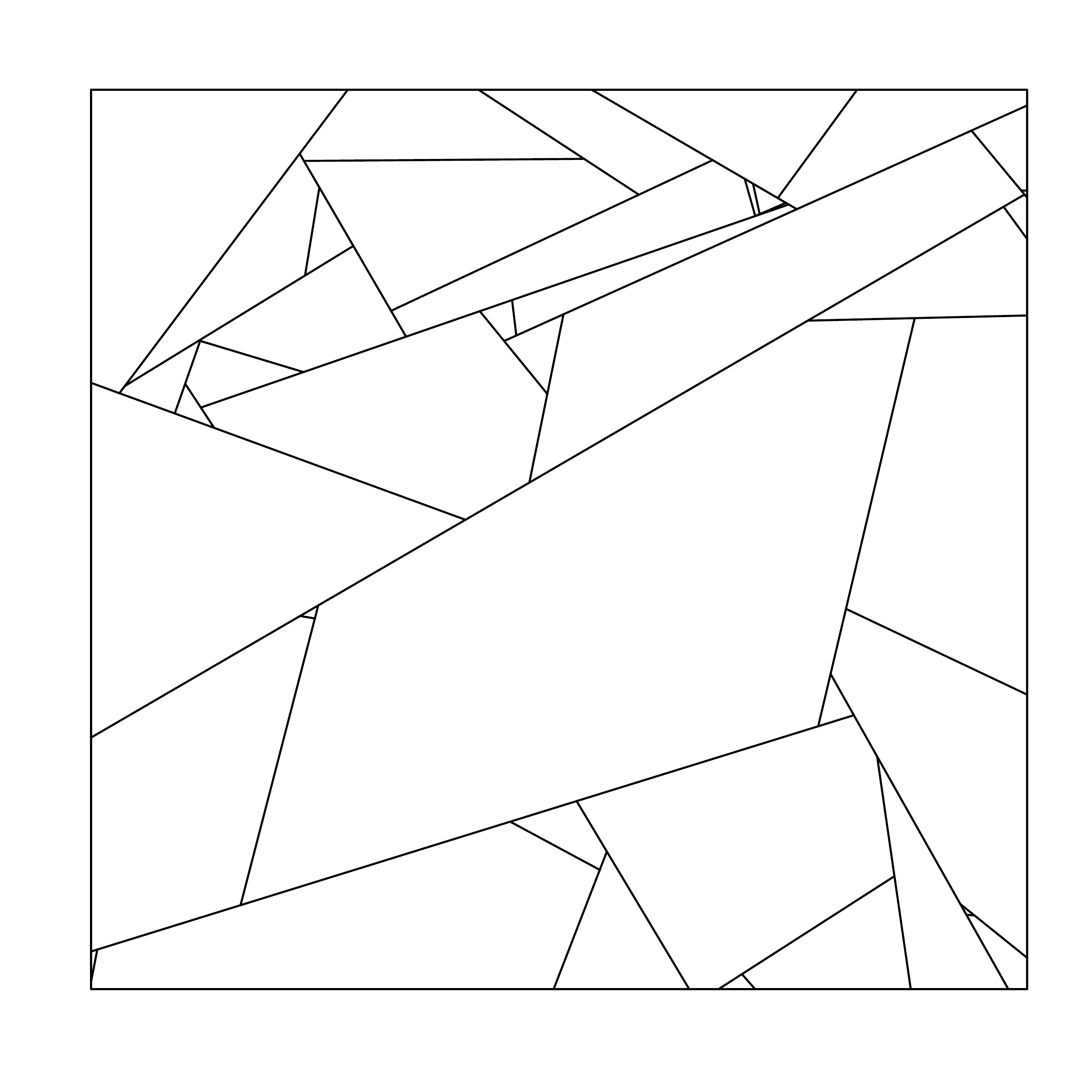} 
\caption{\label{fig:vis_Gibbs_T} Examples of Gibbsian T-tessellations. Left: Completely random T-tessellation, middle: penalizing cell area variability ($\theta_1= 3, \theta_2=10000$), right: favouring cell area variability ($\theta_1= -0.5, \theta_2=-5$)}
\end{center}
\end{figure}

Additionally, they define densities that yield more homogeneous cell area distributions and less small angles between segments than the CRTT model. Let $a^2(T)$ be the sum of squared cell areas and consider the density based on
$$
- \log h(T)= \theta_1\mathring{n}_s(T) + \theta_2 a^2 (T).
$$
Here, $\theta_1$, $\theta_2\in\mathbb{R}$ are parameters allowing to control the number of segments and the homogeneity of the cell areas. Some model realizations are shown in Fig.~\ref{fig:vis_Gibbs_T}.

Inference of model parameters of Gibbsian T-tessellations by Monte Carlo maximum likelihood methods is discussed in \cite{adamczyk-chauvat_gibbsian_2020}. As an application example, they consider the modelling of agricultural landscapes. Goodness-of-fit testing of the models is also discussed. 

\subsubsection{Column tessellations}
A tessellation model in $\R^3$ called column tessellation is introduced in \cite{nguyen_column_2014}. The construction is based on a stationary planar tessellation. Its cells form the base of infinite three-dimensional cylindrical columns which are cut at random heights to obtain the cells of the 3D tessellation. Mean value formulas for this non-face-to-face tessellation are derived in \cite{nguyen_column_2014}.

\section{Iterated tessellations}

\begin{figure}
\begin{center}
\includegraphics[width=3.7cm]{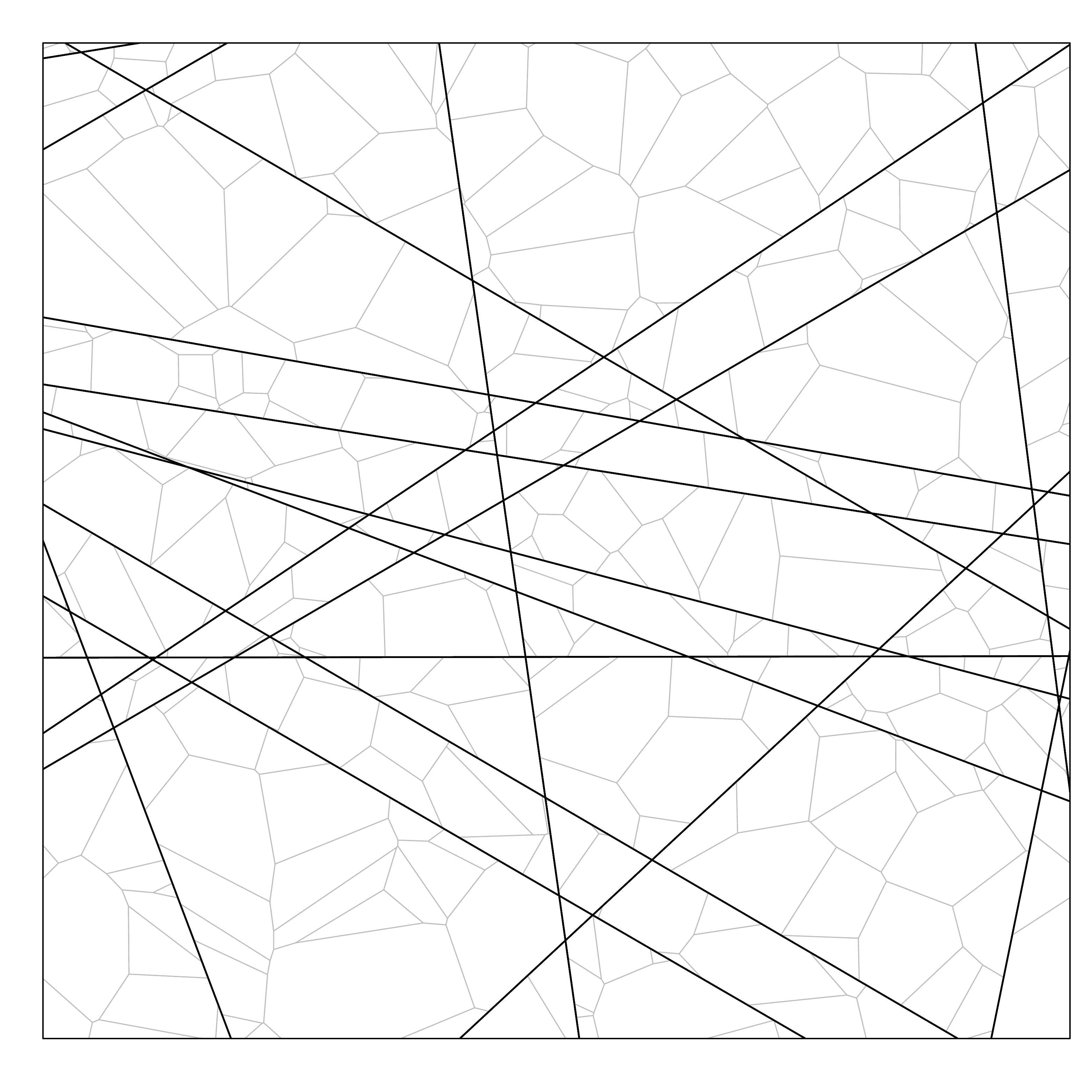}
\includegraphics[width=3.7cm]{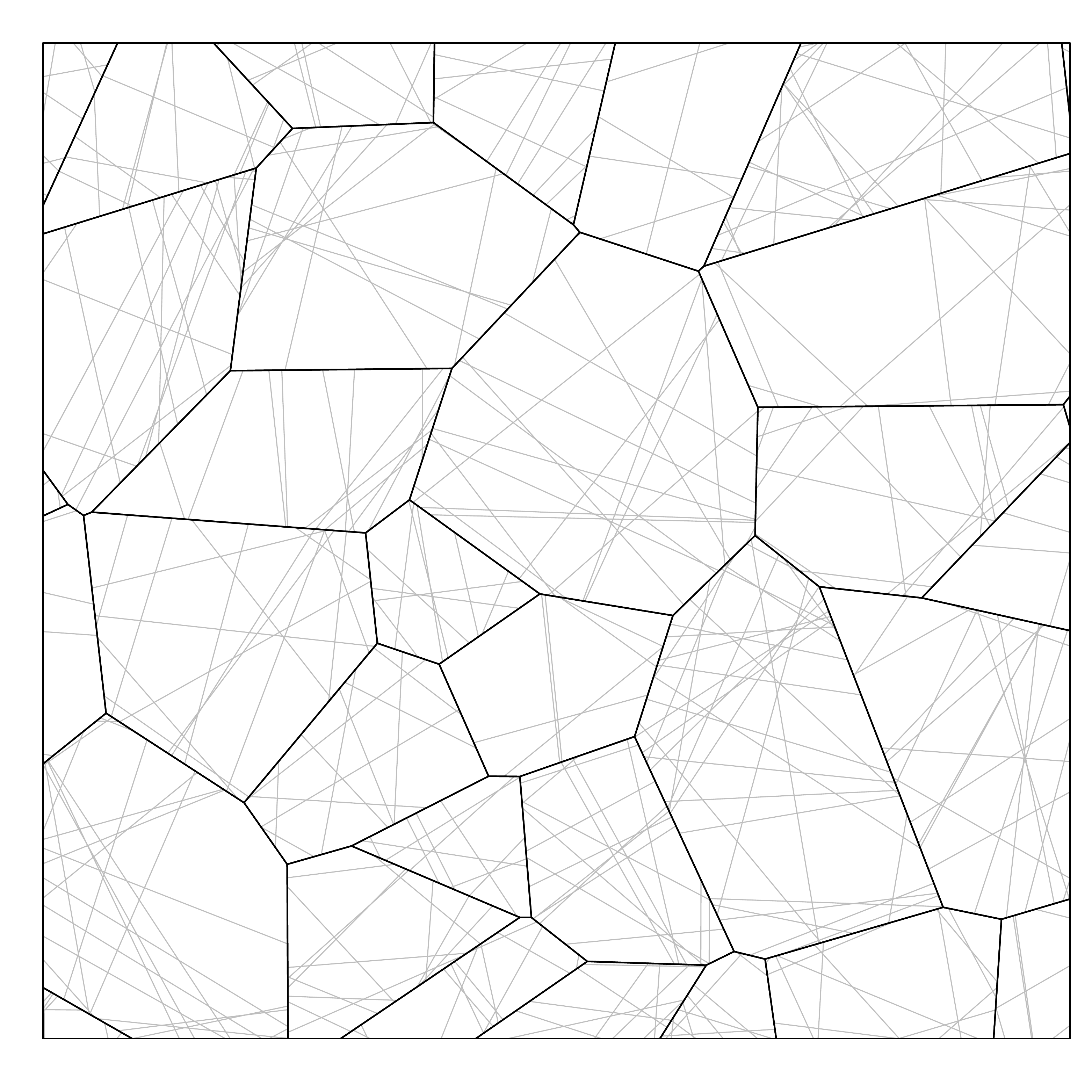}\includegraphics[width=3.7cm]{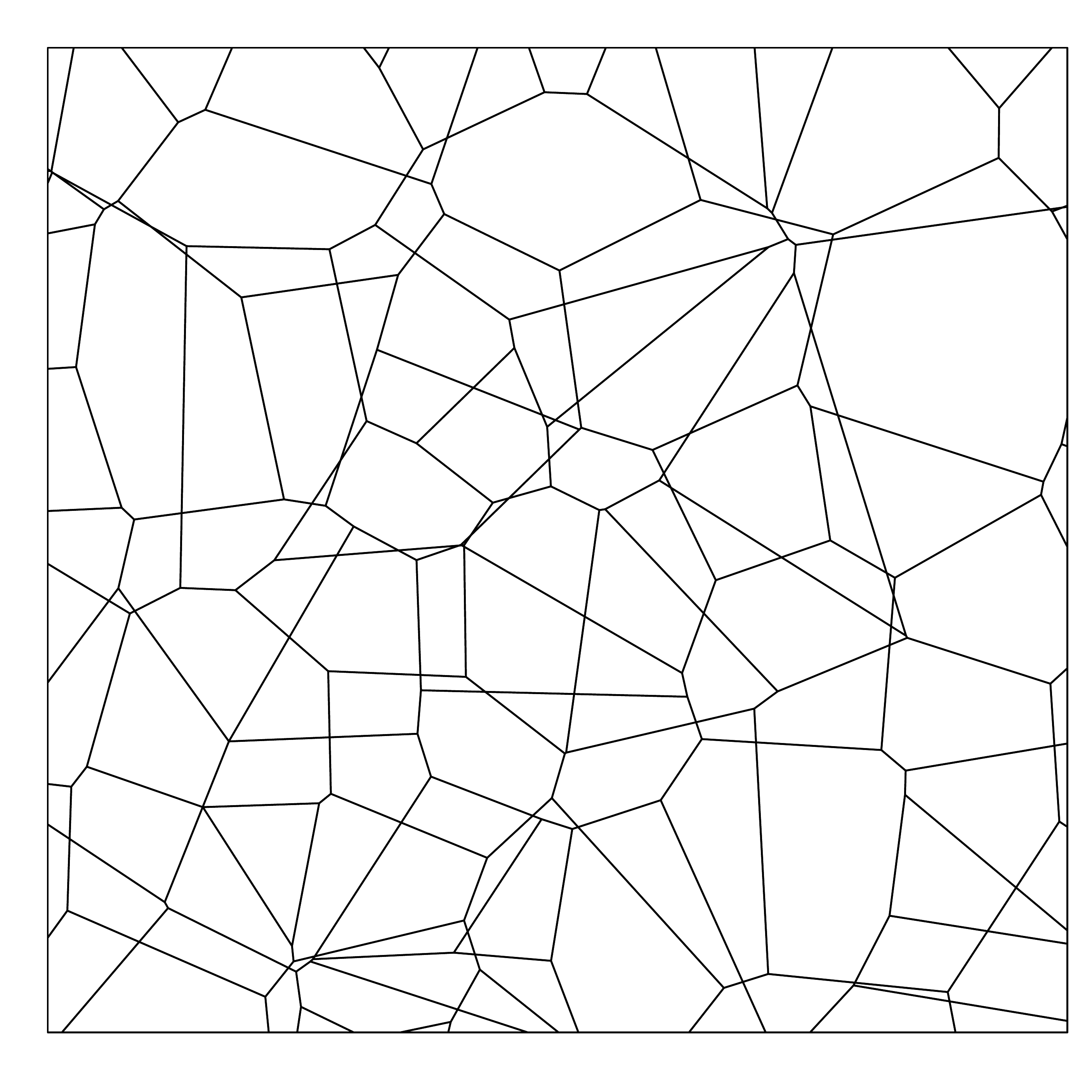}
\end{center}
\caption{Examples of iterated tessellations. Left: Nesting: Poisson-Voronoi tessellations nested into the cells of a Poisson line tessellation. Middle: Nesting: Poisson line tessellations nested into the cells of a Poisson-Voronoi tessellation. Right: Superposition of two Poisson-Voronoi tessellations.}\label{Fig:Nested}
\end{figure}

Iteration of tessellations has already been discussed in the context of STIT tessellations. In \cite{MaierIterated}, a theory for iterated tessellations is developed where the cells of an initial tessellation are further subdivided by component tessellations.  
In \emph{nested tessellations}, an i.i.d. sequence $X_1, X_2, \ldots$ of realizations of the component tessellation is sampled. Then the part of tessellation $X_i$ that intersects the i-th cell of the initial tessellation is copied into that cell. 
If all component tessellations are identical, a \emph{superposition} of two tessellations is obtained. 
Examples for both types of operations are shown in Fig.~\ref{Fig:Nested}.

In \cite{MaierIterated}, conditions ensuring stationarity and isotropy of iterated tessellations in $\R^d$ are formulated. Additionally, formulas for the $k$-face intensities and the expected intrinsic volumes of the $k$-faces are derived. The authors also consider Bernoulli nesting, where only a selection of the cells of the initial tessellation are subdivided.

The distribution of the typical cell of stationary iterated tessellations obtained by choosing the initial and the component tessellations as Poisson-Voronoi and Poisson line tessellations is studied in \cite{maier_distributional_2004}. To this end, the authors formulate a simulation algorithm for the typical cell and derive the distributions of various cell characteristics in a simulation study.
Superpositions of planar Poisson-Voronoi tessellations are studied in \cite{Baccelli2000}.

\section{Dead leaves tessellation}

\begin{figure}
\begin{center}
\includegraphics[width=3.7cm]{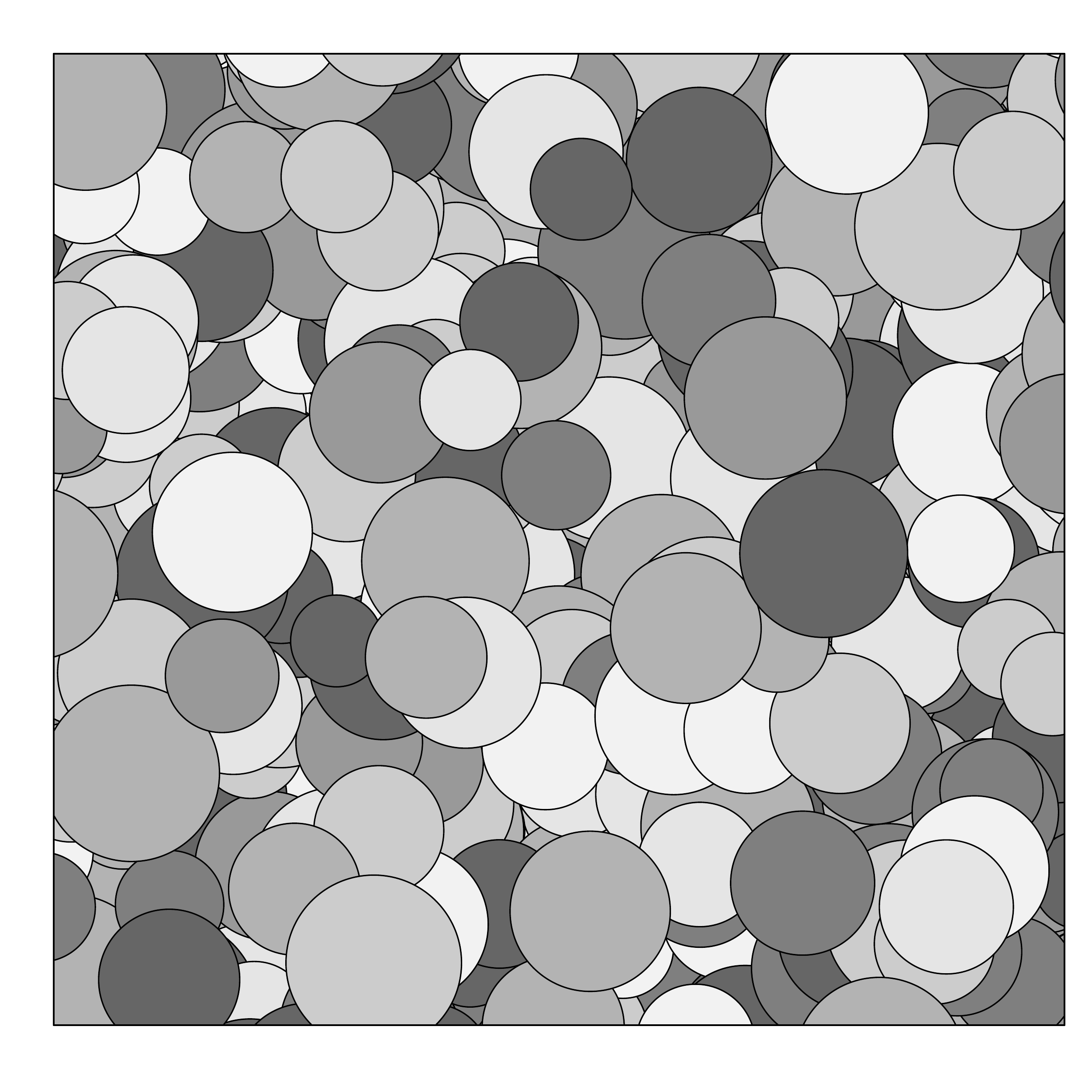}
\includegraphics[width=3.7cm]{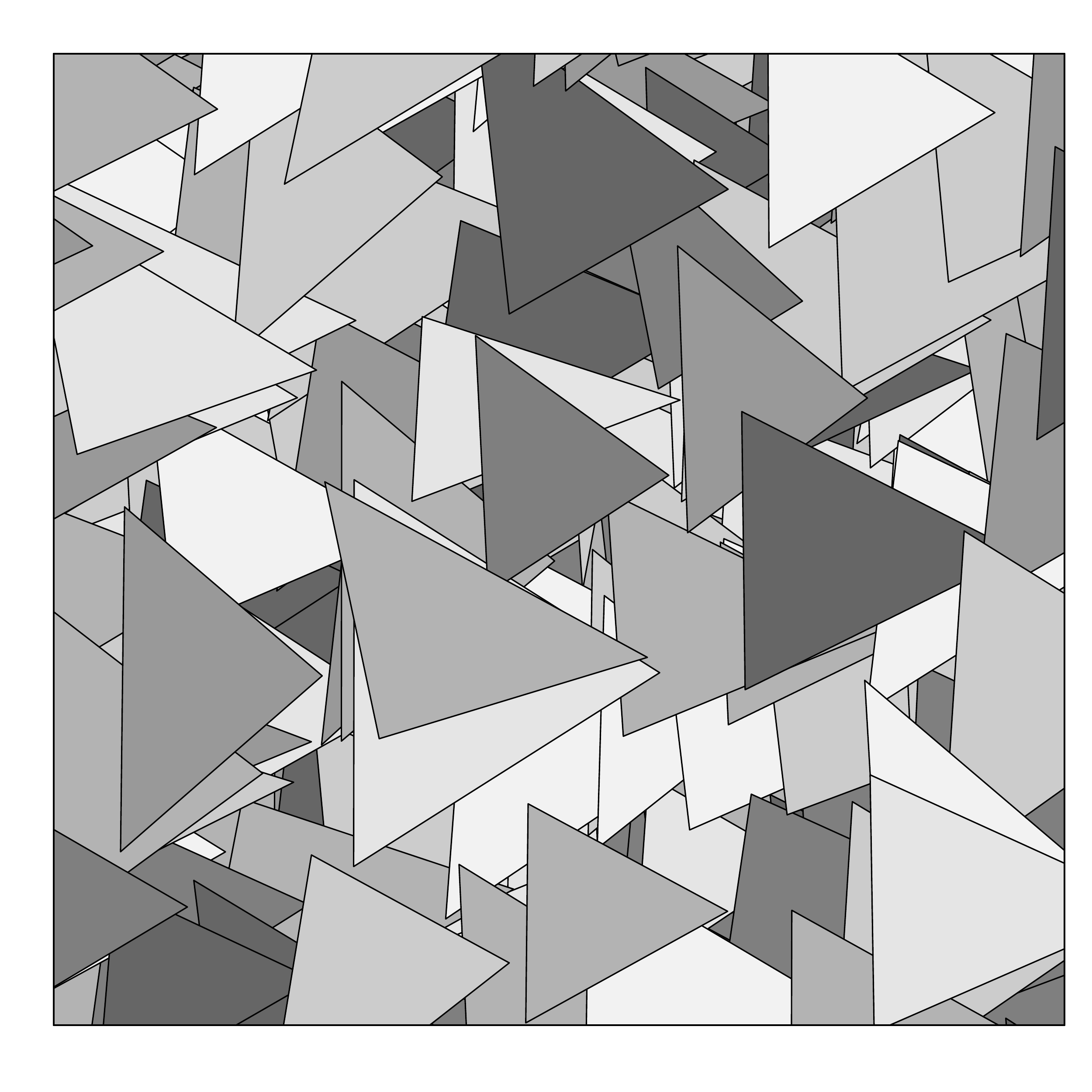}
\end{center}
\caption{Dead leaves tessellations with circular and triangular leaves.}\label{Fig:DeadLeaves}
\end{figure}

The \emph{dead leaves tessellation}, also called \emph{falling leaves tessellation} was introduced in \cite{MatheronDeadLeaves}. As a tessellation of $\R^2$, it is generated by a process that mimics the pattern that fallen leaves generate on the ground. Leaves of given (deterministic or random) shape are placed on the plane such that their positions, orientations and arrival times form a Poisson point process in $\R^2 \times [0, 2\pi) \times [0, \infty)$. Leaves falling later will cover earlier ones. At some point, the ground will completely be covered. The cells of the dead leaves tessellation are then given by all visible leave parts. Some realizations are shown in Fig.~\ref{Fig:DeadLeaves}. Mean value results for dead leaves tessellations with polygonal leaves are given in \cite{CowanFallingLeaves}. Tessellations of leaves with curved boundaries are also briefly discussed, see \cite{penrose2020leaves} for a more general discussion. 

The dead leaves construction can also be generalized to $\R^d$ with $d>2$. Analytical results for this general case are derived in \cite{BordenaveDeadLeaves} and \cite[Ch. 11]{Jeulin2021}.
Perfect simulation of the dead leaves tessellation is discussed in \cite{kendall_perfect_1999}.
Dead leaves models were found to be statistically similar to natural images \cite{GousseauDeadLeaves, JEULIN1989403, Jeulin2021, LeeDeadLeaves2001}.

\section{Fitting of random tessellation models}
\label{sec:Fitting}
\begin{figure}
\begin{center}
\includegraphics[width=3.7cm]{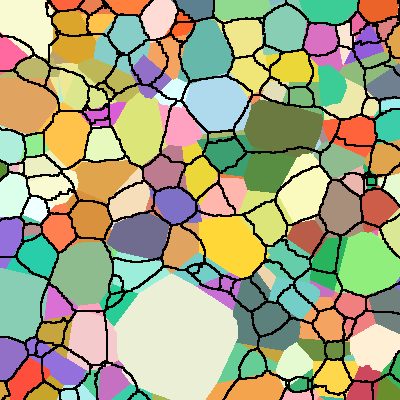}
\includegraphics[width=3.7cm]{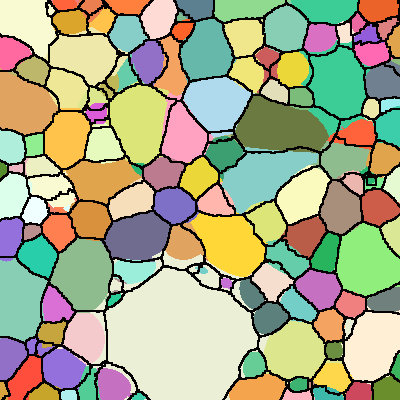}\\
\end{center}
\caption{Sectional images of the 3D approximation of an aluminum foam by a Laguerre tessellation (left) and a GBPD (right). The microstructure's cell boundaries are black while the tessellation cells are colored. }\label{Fig:rekoX2}
\end{figure}

Modelling an observed cellular structure aims at finding realizations of particular model classes that accurately represent the observed cell system. This task can be interpreted in various ways. We will use the following terms:
\begin{enumerate}
    \item 
    \textbf{parametric stochastic modelling}: fitting of a parametric model, e.g. a Poisson-Voronoi tessellation or a Laguerre tessellation generated by a given parametric model of a marked point process, 
    \item 
    \textbf{stochastic reconstruction}: use a stochastic, nonparametric optimization approach to construct a random tessellation that reproduces selected statistics, e.g., the cell size distribution, 
    \item 
    \textbf{approximation or inversion}: for a given model class, e.g., a Laguerre tessellation, find a set of generators such that the resulting tessellation approximates the observed cell system, see Fig.~\ref{Fig:rekoX2} for an example. Goodness of the approximation is measured by a suitable discrepancy measure. In the literature, this problem is also called reconstruction. If the observed cell structure belongs to the selected model class, the fit should be perfect. In this case, one also speaks of inversion of the tessellation.
\end{enumerate}

Approaches 1. and 2. can be used to simulate an arbitrary number of model realizations of variable size that are statistically similar to the data. In contrast, approach 3. just yields one representation of the observed data within the selected model class. 

\subsection{Parametric models}

Parametric stochastic modelling is based on the selection of a suitable parametric model class. Classical models that are considered in this context are Poisson-Voronoi tessellations or Laguerre tessellations generated by dense packings of spheres with a given parametric distribution of radii or volumes. Model fitting then consists in finding model parameters such that statistics like moments of cell characteristics are fit, see \cite{redenbach08fitting, redenbach09microstructure} for Laguerre tessellations and \cite{leon_parameter_2023} for variants of the STIT model. In the latter, also multi objective optimization of several cell characteristics is considered.  

An alternative approach is to consider the probability $P(K)$ for a compact set $K$ to be included in one cell of the tessellation. Explicit formulas for $P(K)$ for various generalized Voronoi tessellations are given in \cite{JEULIN2014139}. Tessellation parameters can be estimated by minimizing some distance between the empirical and the theoretical (or simulated) $P(K)$ for a collection of compact sets $K$.

Due to their formulation based on a density, Gibbs models allow for parameter estimation by (pseudo-) maximum likelihood, see \cite{dereudre_practical_2011} for Delaunay and Voronoi tessellations, \cite{seitl_exploration_2021} for Laguerre tessellations, and \cite{adamczyk-chauvat_gibbsian_2020} for T-tessellations.

\subsection{Reconstruction techniques}

\textbf{Stochastic reconstruction} is mostly considered for Voronoi tessellations and their generalizations. One approach is introduced in  \cite{neper-paper}. It is based on nonlinear optimization, minimizing the discrepancy between prescribed, continuous probability density functions and the sample distributions. In the examples, cell size and sphericity distributions were used. 
An algorithm for constructing Laguerre tessellations of given volumes has been proposed in \cite{givenVolumes}. It uses convex optimization and also includes a step for regularization. Stochastic reconstruction based on Laguerre tessellations generated by Gibbs processes has been discussed in \cite{seitl_exploration_2021}. It uses the Metropolis-Hastings Birth-Death-Move algorithm to iteratively propose new generator configurations. These are accepted with a probability based on the Hastings ratio comparing the density values for the current and the suggested configuration. The proposed algorithm is able to reconstruct moments of cell characteristics but also the whole histograms \cite{seitl_exploration_2021}. 

 A reconstruction approach for T-tessellations could be based on modifications of an initial line system by the split, flip, and merge operations introduced in Section \ref{sec:TTess}.

\textbf{Approximation approaches} have been suggested for various types of input data. The problem of approximating an arbitrary given tessellation by a Voronoi tessellation has been discussed in \cite{voropoly}. The work was later on extended to Laguerre tessellations in \cite{andre}. The idea is to minimize the total mismatch area between the cells of the approximating tessellation and the observed cells
via a gradient descent-based method. The inversion problem of 
finding a system of generators of a Laguerre tessellation was solved in \cite{inverting}.

In cases where the input is given as voxel data, the tessellation cells are voxelized, too, instead of being described explicitly. In \cite{linproGBPD}, an approximation technique based on solving a linear integer programming problem is presented. It uses the relation of tessellations and optimal clusterings and formulates the problem as a weight-balanced least-squares assignment problem. Furthermore, in \cite{neper-paper}, a nonlinear optimization method is proposed  that minimizes the distances between the observed cell system and the tessellation. Stochastic optimization methods such as simulated annealing \cite{ulm2} and a cross-entropy method \cite{ulm-cross-entropy} were considered as well. The idea of these approaches is to iteratively modify randomly chosen generators. The new configuration is accepted with a probability depending on how well the current and the suggested configuration fit. 

Some measurement techniques like X-ray diffraction microscopy yield even less information on the observed cell structure. Instead of a discretization of the cells, only their volumes and centers of mass are reported \cite{petrich_efficient_2021}. In this case, Lyckegard et al. \cite{lyckegaard} propose a simple  heuristic for choosing a set of generators of an approximating Laguerre tessellation. Their solution is often chosen as initial configuration for optimization methods further improving the fit to the real data. Such methods include the minimization of the discrepancies of volume-equivalent spheres via nonlinear optimization \cite{neper-paper} and a cross-entropy method \cite{ulm-cross-entropy-2}.

The approximation methods mentioned above are all based on the Laguerre tessellation, assuming that the observed structure is isotropic. However, this assumption does not always hold true and individual cells may indeed exhibit strong anisotropy \cite{altendorf}. Comparably little research has been conducted on capturing these anisotropies. Tessellations with elliptical generators such as the GBPD seem to be superior to models such as the Laguerre tessellation \cite{SedivyModelSelection}, but require higher computational effort. Optimization techniques discussed above, e.g., the linear integer programming approach \cite{linproGBPD}, simulated annealing \cite{ulm2} or gradient-descent-based methods \cite{petrich_efficient_2021}, can also be applied to GBPDs.

\section{Simulation of random tessellations}
Computer simulations of random tessellation models allow for a visualization of the cell geometry and for comparison of tessellation models. For many tessellation models, no analytical formulas for their cell characteristics are known. In these cases, efficient simulation of the tessellation structure is required to study distributions of cell characteristics. Additionally, simulation is required for most model fitting approaches discussed in Section \ref{sec:Fitting}.

Polyhedral tessellation cells can easily be represented in a data structure listing the vertex coordinates of each polyhedron together with the information which vertices are linked by an edge and which edges form a facet/cell. For tessellations with nonconvex cells, this representation is not sufficient. In some cases, approximations of such tessellation by models with convex cells are considered \cite{OkaBooSugChi00,schaller}.

Various software packages offer implementations of tessellation computation. A selection is given in Table~\ref{Tab:Software}.

\subsection{Line based tessellations}
The T-tessellations in Section \ref{sec:TTess} are mostly defined by an explicit simulation protocol. Hence, we only discuss simulation of Poisson hyperplane processes here.

The decomposition \eqref{thm:eqLambda} of the intensity measure  suggests how to sample this model, see \cite{skm13} for details. Writing a hyperplane $H=L+x$, the coordinate $x\in L^{\perp}$ represents the signed distance of $H$ from the origin. The direction of $H$ is determined by $L$. 
Thus, a hyperplane process corresponds to a point process on the \emph{representation space}
$$
C=\{ (H,p) \in \cL_0 \times \R\}.
$$
In $\R^2$, $C$ can be identified with the set $(0, 2\pi] \times \R$, where $\alpha \in (0, 2\pi]$ yields the line direction via $u=(\cos\alpha, \sin \alpha)$ and $p \in \R$ is the signed distance of the line from the origin.
In a stationary Poisson line process with intensity $\lambda$, the number of lines hitting a set $W$ is Poisson distributed with parameter
$$
\mu_W= \frac{\lambda}{\pi} L(W),
$$
where $L(W)$ is the perimeter of $W$. The easiest way of simulating a Poisson line process on a compact window $W$ is to consider a disc $B(0,r)$ of radius $r$ and centered in the origin such that $W \subset B(0,r)$. Then, a Poisson point process with intensity $\mu= \frac{\lambda}{2\pi}$ is simulated on $(0, 2\pi] \times [-r,r]$. Lines of this process that do not intersect the original window $W$ are rejected.

A Poisson hyperplane process in $\R^3$ is obtained by simulating a Poisson point process on a subset of the representation space $S^{d-1} \times \R$. Here, the intensity is given by
$$
\mu_W= \lambda \overline{b}(W), 
$$
where $\overline{b}(W)$ is the mean width of $W$.

\subsection{Voronoi and Delaunay tessellations}
The computation of Voronoi and Delaunay tessellations and their weighted generalizations is more complex. A summary of algorithms can be found in Chapter 4 of the monograph \cite{OkaBooSugChi00}.

\begin{table}[!t]
\caption{Overview of software packages for simulation of random tessellations}
\label{Tab:Software} 

\begin{tabular}{p{1.5cm}p{5cm}p{1.5cm}p{1.5cm}p{1.5cm}}
\hline\noalign{\smallskip}
Names & Tessellation & Dimension &Language & References\\
\hline\noalign{\smallskip}
Neper & Voronoi, Laguerre & 2D, 3D & C &\cite{Neper}\\
Qhull & Voronoi, Delaunay& 2D, 3D, nD & C& \cite{qhull}\\
CGAL&Voronoi, Laguerre, Delaunay&2D, 3D, nD& C++ &\cite{cgal:eb-23b}\\
CGAL&Johnson-Mehl, Segment Voronoi&2D& C++ &\\
Voro++ &Voronoi, Laguerre & 3D& C++&\cite{Voro}\\
Pomelo & Set Voronoi & 3D & C++& \cite{schaller,pomelo}\\
VRONI & Voronoi, Segment Voronoi, Arc Voronoi & 2D & C++&\cite{Vroni}\\
Gibbs& Gibbs Laguerre & 3D & C++ & \cite{SeitlGithub}\\
wevo & multiplicatively weighted Voronoi &2D &C++ & \cite{held2020efficient,WevoGithub}\\
LiTe & Gibbsian T-tessellation & 2D& C++ & \cite{KieuGithub}\\
crackPattern& STIT variants from Section \ref{Sec:STITExtension} & 2D &C++&\cite{LeonGithub}\\
MATLAB & Voronoi, Delaunay (based on Qhull) & 2D, 3D, nD &MATLAB &\\
Tess & Voronoi, Laguerre (based on Voro++)& 3D&Python &\\
deldir &Voronoi, Delaunay, centroidal Voronoi &2D&R \\
spatstat & Voronoi, Delaunay, Poisson line &2D &R& \cite{Spatstat, spatstat12}\\
tessellation&Voronoi, Delaunay (based on qhull) & 2D, 3D &R\\
transport& Laguerre, Johnson-Mehl (requires CGAL) &2D &R &\\
 \hline
\noalign{\smallskip}
\noalign{\smallskip}
\end{tabular}
\end{table}

\subsection{Simulation of the typical cell}

For some tessellations, methods for simulating single realizations of the typical cell are available. For the Poisson-Voronoi tessellation, this is possible by the radial simulation method of \cite{Quine_Watson_1984}, see also \cite{mol94}. Adaptions for the Johnson-Mehl tessellation and the Laguerre tessellation are described in \cite{Moller1995_johnsonMehl} and \cite{lautensack06:_random_laguer}, respectively.
The typical cell of a Poisson line tessellation can be sampled by the construction given in \cite{George1987}. An alternative approach, along with a method for sampling the zero cell in a Poisson line tessellation, is discussed in \cite{michel_empirical_2007}.
Sampling of the typical cell of iterated tessellations is discussed in \cite{maier_distributional_2004}.

\subsection{Edge treatment}

\begin{figure}[t]
\begin{center}
\includegraphics[trim={4cm 4cm 3cm 3.5cm},clip,width=3.7cm]{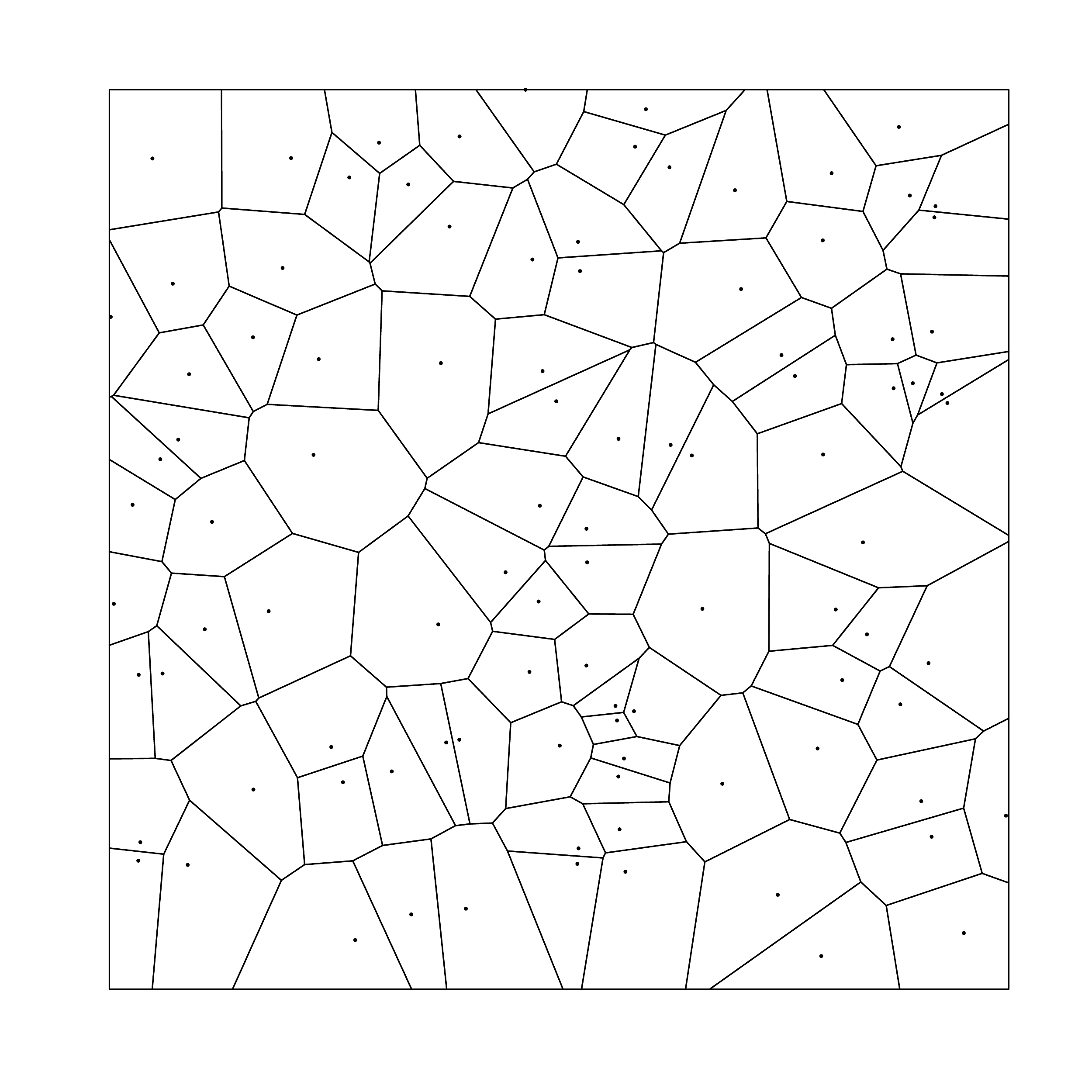} \hspace{0.1cm} \includegraphics[trim={4cm 4cm 3.0cm 3.5cm},clip,
width=3.7cm]{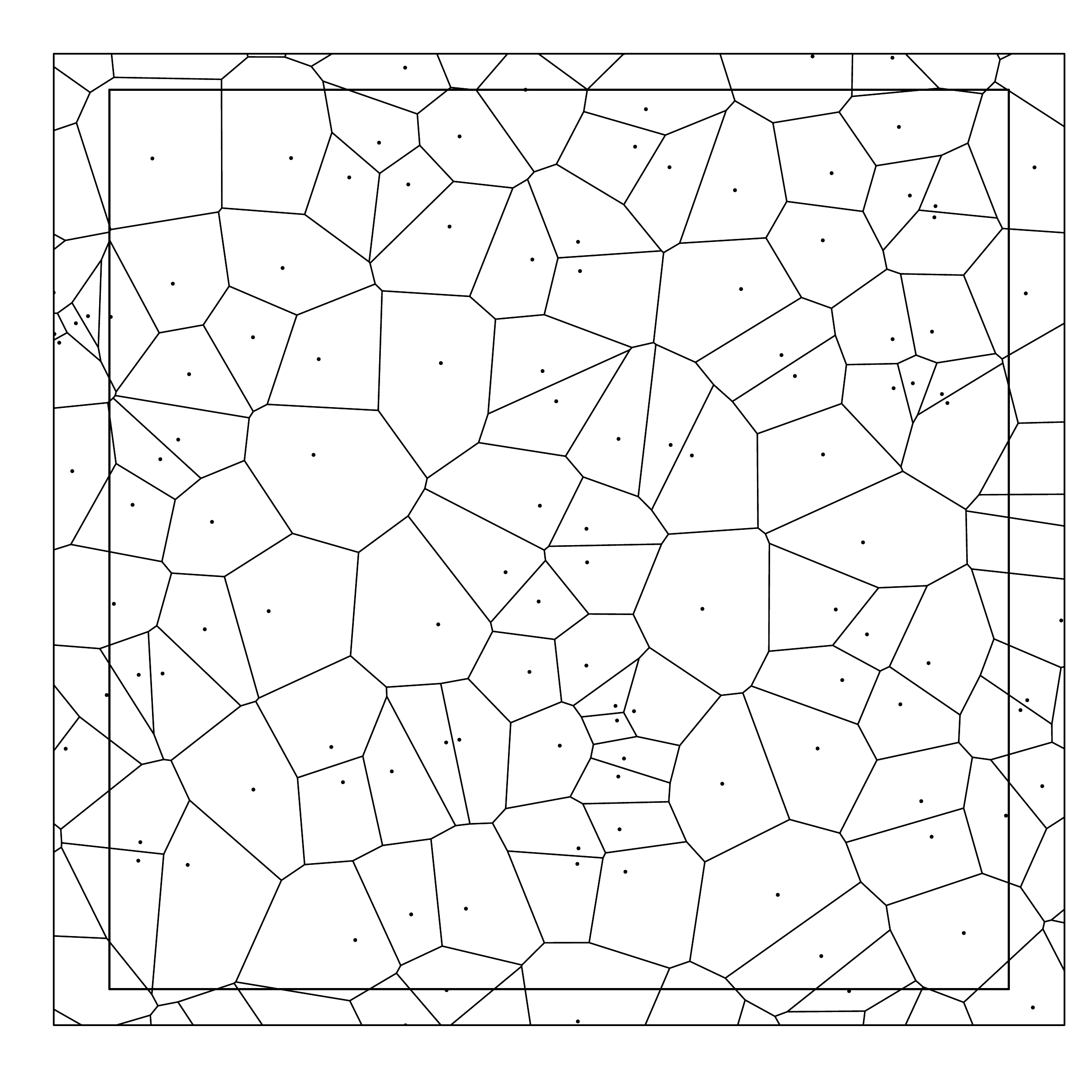} \hspace{0.1cm} \includegraphics[trim={4cm 4cm 3.0cm 3.5cm},clip,width=3.7cm]{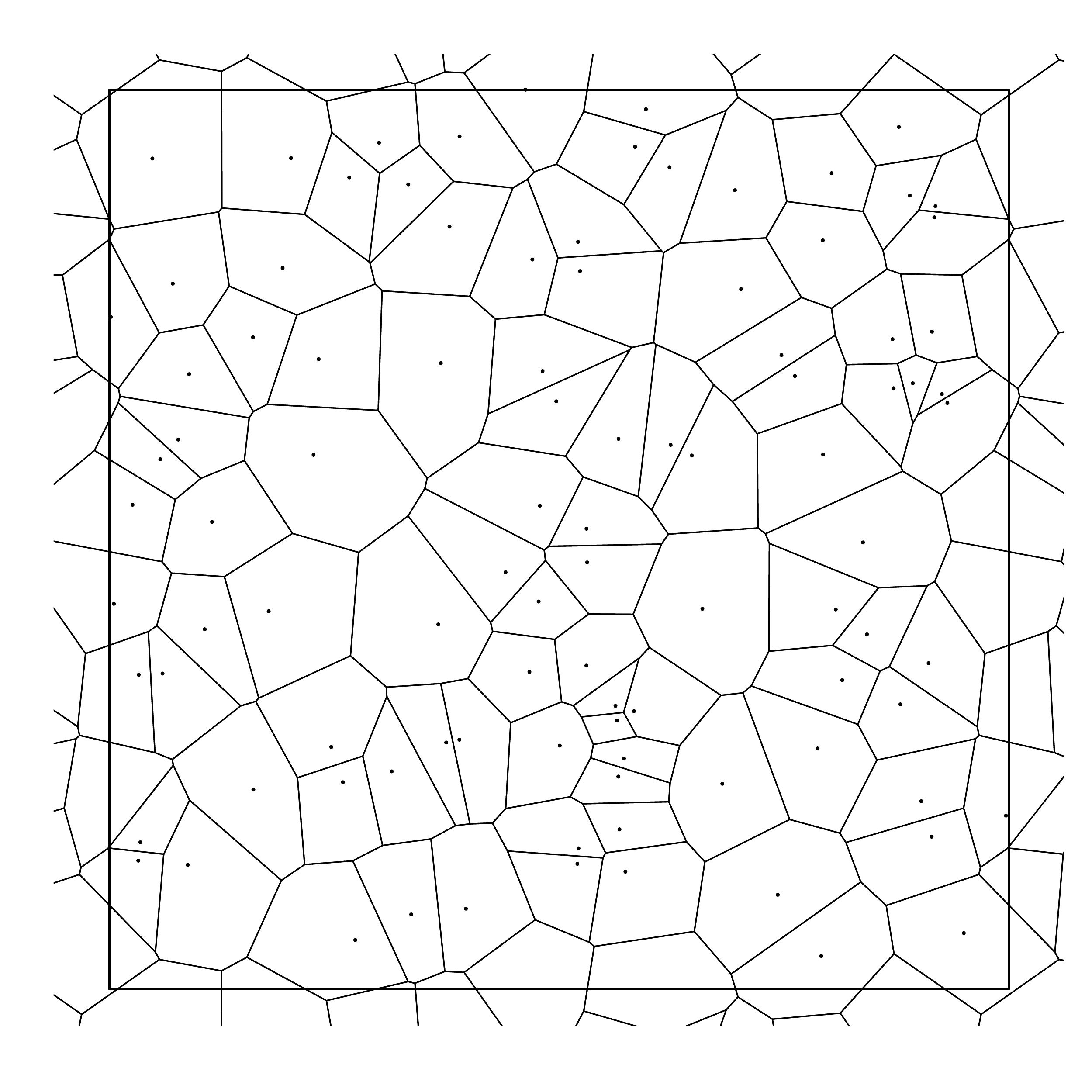}
\caption{\label{fig:vis_Edge} Poisson-Voronoi tessellation on the unit square with different edge treatment. Left: No edge treatment, only points in the unit square are simulated. Note the special shape of boundary cells. Middle: Plussampling. Points are simulated on an extended window. Only the unit square is plotted. Right: Periodic edge treatment.}
\end{center}
\end{figure}

The tessellations discussed in Section~\ref{Sec:Temporal} are mostly defined by explicit simulation processed on bounded windows. Hence, their simulation does not suffer from edge effects. This is different, e.g., for Voronoi and Delaunay tessellations. Here, one has to distinguish between the tessellations (or better diagrams) generated by a finite point set and a cutout of a stationary tessellation.    
This is illustrated in Fig.~\ref{fig:vis_Edge}. In the left panel, a realization of a Poisson process on the unit square is simulated. When computing the Voronoi tessellation of this point pattern, cells intersecting the boundary do not look like parts of the typical cell of a stationary Poisson-Voronoi tessellation. To avoid this effect, the point process can be simulated on a larger window (plussampling). This way, also cells entering the window from the outside can be realized. 

An alternative is periodic edge treatment. That is, the left and the right as well as the top and the bottom edge of the square are identified. Cells that exit the window on the left enter again on the right. In practice, periodic edge treatment is usually realized by working with 9 copies of the original pattern that are arranged in a $3\times 3$ grid. The central square then represents the desired observation window. In 3D, 27 copies in a $3\times 3 \times 3$ grid have to be used.

\bibliographystyle{spmpsci}
\bibliography{litbank,Tessellations}
\end{document}